\documentclass[12pt,reqno]{amsart}

\usepackage{amsmath}
\usepackage{amssymb}
\usepackage{mathrsfs}
\usepackage[OT2,T1]{fontenc}
\usepackage{calligra}
\usepackage[all,cmtip]{xy}
\usepackage{cancel}
\usepackage{bbm}
\usepackage{empheq}
\usepackage{bm}
\usepackage[toc]{appendix}
\usepackage{chngcntr}
\counterwithin{figure}{section}

\usepackage[pdftex,     
          plainpages=false,   
           breaklinks=true,    
           colorlinks=true,
           pdftitle=My Document
           pdfauthor=My Good Self
           colorlinks=true,
	    urlcolor=blue,
	    citecolor=red,
	    linkcolor=blue
          ]{hyperref}
          
\usepackage[a4paper,top=4cm,bottom=4cm,left=2.6cm,right=2.6cm]{geometry}

\theoremstyle{plain}
\newtheorem{theorem}{Theorem}[section]

\theoremstyle{definition}
\newtheorem{defi}[theorem]{Definition}
\newtheorem{regola}[theorem]{Rule}

\newtheorem{es}[theorem]{Example}
\theoremstyle{plain}

\newtheorem{cor}[theorem]{Corollary}
\newtheorem{lemma}[theorem]{Lemma}
\newtheorem{prop}[theorem]{Proposition}
\theoremstyle{remark}
\newtheorem{oss}[theorem]{Remark}

\newlanguage\fakelanguage
\newcommand\cyr{\fontencoding{OT2}\fontfamily{wncyr}\selectfont
   \language\fakelanguage}
\DeclareTextFontCommand{\textcyr}{\cyr}

\numberwithin{equation}{section}

\DeclareMathOperator{\Hom}{Hom}

\usepackage{calligra}

\DeclareMathOperator{\HOM}{\mathscr{H}\text{\kern -3pt {\calligra\large om}}\,}

\makeatletter
\newsavebox{\@brx}
\newcommand{\llangle}[1][]{\savebox{\@brx}{\(\m@th{#1\langle}\)}%
  \mathopen{\copy\@brx\kern-0.5\wd\@brx\usebox{\@brx}}}
\newcommand{\rrangle}[1][]{\savebox{\@brx}{\(\m@th{#1\rangle}\)}%
  \mathclose{\copy\@brx\kern-0.5\wd\@brx\usebox{\@brx}}}
\makeatother

\usepackage{framed,color}
\definecolor{shadecolor}{rgb}{0.92, 0.92, 0.92}
\definecolor{ambra}{rgb}{1.0, 0.75, 0.0}
\definecolor{ametista}{rgb}{0.6, 0.4, 0.8}
\definecolor{auburn}{rgb}{0.43, 0.21, 0.1}
\definecolor{ballblue}{rgb}{0.13, 0.67, 0.8}
\definecolor{cadmiumgreen}{rgb}{0.0, 0.42, 0.24}
\definecolor{candypink}{rgb}{0.89, 0.44, 0.48}
\definecolor{caribbeangreen}{rgb}{0.0, 0.8, 0.6}
\definecolor{airforceblue}{rgb}{0.36, 0.54, 0.66}

\newcommand{\bmh}{\begin{shaded}}
\newcommand{\esh}{\end{shaded}}

\newcommand{\beq}{\begin{equation}}
\newcommand{\eneq}{\end{equation}}
\newcommand{\clgr}[1]{\text{\footnotesize\calligra #1}}

\begin{document}

\title[$qDE$ and $qKZ$ equations for a Projective space]{Equivariant quantum differential equation and qKZ equations for a projective space: Stokes bases as exceptional collections, Stokes matrices as Gram matrices, and \textcyr{B}-Theorem}
\author[Giordano Cotti and Alexander Varchenko]{Giordano Cotti$\>^\circ$ and Alexander Varchenko$\>^\star$}
\maketitle
\begin{center}
\textit{ $^\circ\>$Max-Planck-Institut f\"{u}r Mathematik\\ Vivatsgasse 7, 53111 Bonn, Germany\/}

\vskip4pt
\textit{ $^{\star}\>$Department of Mathematics, University
of North Carolina at Chapel Hill\\ Chapel Hill, NC 27599-3250, USA\/}

\vskip4pt
\textit{ $^{\star}\>$Faculty of Mathematics and Mechanics, Lomonosov Moscow State
University\\ Leninskiye Gory 1, 119991 Moscow GSP-1, Russia\/}

\end{center}

{\let\thefootnote\relax
\footnotetext{\vskip5pt 
\noindent
$^\circ\>$\textit{ E-mail}:  gcotti@mpim-bonn.mpg.de, gcotti@sissa.it,
supported by Max-Planck-Institut f\"{u}r Mathematik, Bonn, Germany
\\
$^\star\>$\textit{ E-mail}:  anv@email.unc.edu,
supported in part by NSF grant DMS-1665239}}

\begin{abstract}
 In \cite{tarvar}
the equivariant quantum differential 
equation ($qDE$)  for a projective space 
 was considered and a compatible system of 
difference $qKZ$ equations was introduced;
 the  space of solutions to 
the joint system of the $qDE$ and $qKZ$ equations
was identified with the space 
of the equivariant $K$-theory algebra of the projective space; 
 Stokes bases in the space of solutions were identified with exceptional
bases in the equivariant $K$-theory algebra. This paper is a continuation of \cite{tarvar}.

We describe the relation between solutions to the joint system
of the $qDE$ and $qKZ$ equations and the topological-enumerative solution to the $qDE$
only, defined as a  generating function of equivariant descendant Gromov-Witten
invariants. The relation is  in terms of the equivariant graded Chern character on
the equivariant $K$-theory algebra,  the equivariant Gamma class of the projective space, 
and the equivariant first Chern class of the tangent bundle of 
the projective space.

We consider a Stokes basis, the associated exceptional basis in the equivariant $K$-theory
algebra, and the associated Stokes matrix.  We show that the Stokes matrix 
equals the Gram matrix of the equivariant Grothendieck-Euler-Poincar\'{e} pairing wrt 
to the  basis, which is the left dual to the associated exceptional basis.

We identify the Stokes bases  in the space of solutions with  explicit
full exceptional collections in the equivariant derived 
category of coherent sheaves on the projective space, where the elements of those
exceptional collections 
are just line bundles on the projective space  and exterior powers of the tangent bundle of
the projective space. 

These statements  are equivariant analogs  of results of G.\,Cotti, B.\,Dubrovin, D.\,Guz-zetti, 
and S.\,Galkin, V.\,Golyshev, H.\,Iritani.

\end{abstract}

\newpage
\begin{flushright}
\emph{In memory of Boris Dubrovin (1950-2019)}
\end{flushright}
\tableofcontents

\section{Introduction}
\addtocontents{toc}{\protect\setcounter{tocdepth}{1}}

\subsection{} We consider the equivariant quantum differential equation ($qDE$) of a complex projective space $\mathbb P^{n-1}$
with the diagonal action of the torus $\mathbb T=(\mathbb C^*)^n$. This equation is the ordinary differential equation
\beq
\left(q\frac{d}{dq}-x*_{q,\bm z}\right)I(q,\bm z)=0,
\eneq
where the unknown function $I(q,\bm z)$ takes values in the equivariant cohomology algebra $H^\bullet_\mathbb T(\mathbb P^{n-1},\mathbb C)$, and $x*_{q,\bm z}\colon H^\bullet_\mathbb T(\mathbb P^{n-1},\mathbb C)\to H^\bullet_\mathbb T(\mathbb P^{n-1},\mathbb C)$ is the operator of \emph{quantum} multiplication by the equivariant first Chern class of the tautological line bundle  on $\mathbb P^{n-1}$. 
The  $qDE$ depends on the equivariant parameters $\bm z=(z_1,\dots, z_n)$ 
 corresponding to the factors of the torus $\mathbb T$. The  $qDE$ has two singularities: a regular singularity at $q=0$ and an irregular singularity at $q=\infty$.

In \cite{tarvar} a compatible system of difference equations, called the $qKZ$ equations, 
was introduced,
\beq
I(q,z_1,\dots, z_i-1,\dots, z_n)=K_i(q,\bm z) I(q,\bm z),\quad i=1,\dots, n,
\eneq
where $K_i$'s are suitable linear operators. In \cite{tar-var} solutions to the joint system of the
 $qDE$ and $qKZ$ equations 
were constructed in the form of $q$-hypergeometric integrals. 
In \cite{tarvar} the space of solutions was  identified with
the space of the equivariant $K$-theory algebra $K_0^{\mathbb T}(\mathbb P^{n-1})$. The  Stokes bases of
the $qDE$ at its irregular singular point were described in terms of the exceptional bases of
the equivariant $K$-theory  and a suitable braid group action on the set of
exceptional bases.   In this paper we continue this study.

\subsection{}
We establish relations between the monodromy data of the joint system of the $qDE$ and $qKZ$ equations for $\mathbb P^{n-1}$ and characteristic classes of objects of the derived category 
$\mathcal D^b_\mathbb T(\mathbb P^{n-1})$ 
of equivariant coherent sheaves on
$\mathbb P^{n-1}$. 

\smallskip
Our first result is on the relation between solutions to the joint system
of the $qDE$ and $qKZ$ equations and the topological-enumerative morphism.
 
 The topological-enumerative morphism is the map 
 $\mathcal S^{o}$, which assigns   a solution of the $qDE$ (only) to 
  every element of $H^\bullet_\mathbb T(\mathbb P^{n-1},\mathbb C)$
  and  whose expansion at $q=0$ is the generating function for  the {equivariant} descendant  Gromov-Witten invariants of $\mathbb P^{n-1}$.

 For $E\in K_0^{\mathbb T}(\mathbb P^{n-1})$ let $\theta(E)$ be
the solution to the joint system of
the $qDE$ and $qKZ$ equations, assigned to $E$ in \cite{tarvar}.  Our \textcyr{B}-Theorem  \ref{bteo}
says that 
$$
\theta(E) = \mathcal S^o
 \big( e^{\pi\sqrt{-1}c_1(\mathbb P^{n-1})} \widehat\Gamma^+_{\mathbb P^{n-1}} {\rm Ch}_{\mathbb T}(E)\big),
$$
where  $c_1(\mathbb P^{n-1})$ is  the equivariant first Chern class of the tangent bundle of $P^{n-1}$, 
$\widehat\Gamma^+_{\mathbb P^{n-1}} $ is the equivariant Gamma-class of $P^{n-1}$,
${\rm Ch}_{\mathbb T}(E)$ is the equivariant graded Chern character of $E$.
In other words, we have the following commutative diagram:
$$
\xymatrix{
\text{$K$-Theory}\ar^{\textcyr{B}\quad\quad\quad}[rr]\ar_{\theta}[dr]&&
        \text{Equiv. Cohomology}
\ar^{\mathcal S^o}[dl]\\
&\text{Solutions of $qDE$}&
}
$$
where $\textcyr{B}(E) := e^{\pi\sqrt{-1}c_1(\mathbb P^{n-1})} \widehat\Gamma^+_{\mathbb P^{n-1}} 
{\rm Ch}_{\mathbb T}(E)$.

Notice that the \textcyr{B}-Theorem is an equivariant analog of results of \cite[Section 5]{gamma1} and \cite[Section 6]{CDG1} for projective spaces. Moreover it is a refinement of the Gamma Theorem in  \cite{tarvar, tar-var}.

\smallskip
Our second result is the identification of the Stokes bases in the space of solutions to the joint
system of  the $qDE$ and $qKZ$ equations with explicit  $\mathbb T$-full exceptional collections 
in the derived category $\mathcal D^b_{\mathbb T}(\mathbb P^{n-1})$
of $\mathbb T$-equivariant coherent sheaves on $\mathbb P^{n-1}$.
We show that the elements of these $\mathbb T$-full 
exceptional collections 
are just line bundles $\mathcal O(i)$
on $\mathbb P^{n-1}$  and exterior powers $\bigwedge^j \mathcal T \otimes \mathcal O(i)$
of the tangent bundle $\mathcal T$ of $\mathbb P^{n-1}$ multiplied by line bundles, see Theorem \ref{exobj}, Corollary \ref{corexcol}, Corollary \ref{corexcol-t} and Theorem \ref{stokexcbas}.
This result is an equivariant version of \cite[Corollary 6.11]{CDG1}.

\smallskip
Our third result is on the relation between the Stokes matrices and Gram matrices of the
 Grothendieck-Euler-Poincar\'{e} pairing on $K_0^{\mathbb T}(\mathbb P^{n-1})$.  
 
 Consider a Stokes sector $\mathcal V$ and the complementary Stokes sector
 $e^{\pi \sqrt{-1}}\mathcal V$. Consider the two exceptional bases in 
 $K_0^{\mathbb T}(\mathbb P^{n-1})$ assigned $\mathcal V$ and
 $e^{\pi \sqrt{-1}}\mathcal V$ in \cite{tarvar}.
  The matrix expressing the  second exceptional basis in terms 
 of the first exceptional basis is called the Stokes matrix associated with the Stokes sector $\mathcal V$.
 We show that the second exceptional basis is  left dual to the first exceptional basis 
wrt  the Grothendieck-Euler-Poincar\'{e} pairing. This fact implies that the Stokes matrix equals
  the Gram
matrix of the Grothendieck-Euler-Poincar\'{e} pairing wrt the second exceptional basis, see Theorem \ref{teostokgram}.
This is an equivariant analog of \cite{guzzetti1} (see also \cite[Section 6]{CDG1} for some refinements of results
in \cite{guzzetti1}).

\subsection{} 

This paper is related to the general theory of D.\,Maulik and A.\,Okounkov in \cite{mo}
 connecting quantum groups and equivariant quantum cohomology of Nakajima quiver varieties. In that context, it was realized that the  $qDE$s of Nakajima quiver varieties admit
some compatible difference equations, called the $qKZ$ equations.

A special case of  Nakajima varieties, namely, the case of the cotangent bundles
$T^*\mathcal F_{\bm \lambda}$   of partial flag varieties $\mathcal F_{\bm \lambda}$ was considered
in \cite{GRTV} and \cite{RTV}
\footnote{Note that the partial flag varieties themselves
are not Nakajima varieties}.  In those papers the $qDE$s and $qKZ$ equations for cotangent bundles
were identified with the dynamical differential equations and $qKZ$ difference equations,
associated in representation theory
with the evaluation module $\Bbb C^N(z_1)\otimes\dots\otimes\Bbb C^N(z_n)$ of the Yangian
$Y(\frak{gl}_N)$.
This identification leads to two constructions of solutions to
 the joint system of the $qDE$ and $qKZ$ equations
for the cotangent bundles.
One  construction in \cite{tarvarhyp} gave solutions in the form of multidimensional
hypergeometric  integrals and
another construction in \cite{tar-var} gave solutions in the form of multidimensional $q$-hypergeometric
integrals.

Also in \cite{tar-var} a suitable limit of the $qDE$s  for cotangent bundles 
of partial flag varieties 
was considered. In that limit
the $qDE$s for cotangent bundles
turn into the $qDE$s for the  partial flag varieties themselves.
Moreover, in that limit the $qKZ$ equations  for
cotangent bundles
survive  and turn into new systems of difference equations
compatible with the $qDE$s  for  partial flag varieties.
These new systems of difference equations were also called the $qKZ$ equations.
Furthermore, it was shown in \cite{tar-var}
 that the $q$-hypergeometric solutions  to  the joint systems of the $qDE$s and $qKZ$ equations
for  cotangent bundles
 have a limit when the $qDE$s and $qKZ$ equations turn 
into the $qDE$s and $qKZ$ equations for partial flag varieties.

The special case of that  limit was considered in \cite{tarvar}
for the partial flag variety $\mathbb P^{n-1}$.
In \cite{tarvar} the 
q-hypergeometric solutions to the joint system of the $qDE$ and $qKZ$ equations 
for $\mathbb  P^{n-1}$
were applied  to
study the monodromy properties of the $qDE$ for $\mathbb P^{n-1}$.

\subsection{}
The paper is organized as follows. 
The basic notions of the derived category of equivariant coherent sheaves and equivariant Helix theory are collected 
in Section \ref{sec2}.
 In Section \ref{sec3} we describe the equivariant derived category and $K$-theory of $\mathbb P^{n-1}$. In Section \ref{sec4} we introduce the equivariant cohomology of $\mathbb P^{n-1}$. In Section \ref{sec5} we discuss the equivariant Gromov-Witten theory of $\mathbb P^{n-1}$. We 
 introduce the $qDE$ and $qKZ$ difference equations, and the topological-enumerative morphism $\mathcal S^o$.
 
  In Section \ref{sec6} we introduce two fundamental systems of solutions of the $qDE$ (only): 
  the Levelt solution and the topological-enumerative solution. We study how they are related, and we describe their monodromy. 
  
  In Section \ref{sec7} we study solutions to the joint system of the 
  $qDE$ and $qKZ$ equations, their integral representations, their asymptotics. 
We describe the corresponding objects and exceptional collections
 in the derived category. In Section \ref{sec9} we prove our \textcyr{B}-Theorem. 
 
 In Section \ref{sec10} we describe the 
  structure of formal solutions to the joint system of the $qDE$ and $qKZ$ 
  difference equations, see Theorem \ref{formred}.
  
  In Section \ref{sec11} we study the Stokes bases of the space of solutions, their normalizations.
 We show that the Stokes bases correspond to $\mathbb T$-full exceptional collections
 in  $\mathcal D^b_{\mathbb T}(\mathbb P^{n-1})$.
  In Section \ref{sec12} we prove that the Stokes matrices coincide 
  with the Gram matrices of the equivariant Grothendieck-Euler-Poincar\'{e} 
  pairing. 
  
  In Section \ref{sec13} we study the specialization of the 
  $qDE$ at points $\bm z$ such that $(e^{2\pi\sqrt{-1}z_1}$, \dots, $e^{2\pi\sqrt{-1}z_n})$ 
  are roots of unity. We show that the monodromy group of the  $qDE$ is $\mathbb Z_n$ 
  only for this specialization of the equivariant parameters $\bm z$. 
 
  In Appendix \ref{redjoisyst} we prove Theorem \ref{teoappb} on the formal normal form
  for a compatible system of a differential equation and a system of difference equations.
  
  In Appendix \ref{qdeDubr} we discuss the relation between the equivariant $qDE$ and the isomonodromic system of differential equations attached to the quantum cohomology of $\mathbb P^{n-1}$. Such a system plays a central role in Dubrovin's theory of Frobenius manifolds \cite{dubro1,dubro2,dubro0,CDG1}.

\subsection{}
The authors thank A.\,Givental, R.\,Rim\'anyi, M.\,Smirnov, and V.\,Tarasov for useful discussions, and also the referee for suggestions improving 
        the exposition of the paper.
The authors are grateful
to the Max-Planck-Institut f\"{u}r Mathematik in Bonn, where 
this project  was developed, for support and 
excellent working conditions. 
The first author is thankful to his teacher, Boris Dubrovin, for his interest in this project.
He will remember with admiration and gratitude his encouraging guidance.

\addtocontents{toc}{\protect\setcounter{tocdepth}{2}}
\section{Equivariant exceptional collections and bases}\label{sec2}
General references for this Section are \cite{gelman,chrissginz,helix}.
\subsection{Basic notions.}
Let $G$ be a linear algebraic reductive group over $\mathbb C$. We denote by
\begin{itemize}
\item $Rep(G)$ the category of finite dimensional complex representations of $G$,
\item $R(G):=K_0(Rep(G))$ (resp. $R(G)_{\mathbb C}:=R(G)\otimes_\mathbb Z	\mathbb C$) the ring of finite dimensional complex representations of $G$ with integer (resp. complex) coefficients.
\end{itemize}
In particular, for a complex torus $\mathbb T:=(\mathbb C^*)^n$ we have $R(\mathbb T)_{\mathbb C}=\mathbb C[Z_1^{\pm 1},\dots,Z_n^{\pm 1}]$. For short, we set $\bm Z:=(Z_1,\dots, Z_n)$ and $\mathbb C[\bm Z^{\pm 1}]:=\mathbb C[Z_1^{\pm 1},\dots, Z_n^{\pm 1}]$. \newline

Let $X$ be a smooth complex projective variety equipped with the action of $G$. We denote by
\begin{enumerate}
\item $\mathcal D^b(X)$ its derived category of coherent sheaves,
\item $\mathcal D^b_{G}(X)$ its derived category of $G$-equivariant coherent sheaves,
\item $ K_0(X)$ (resp. $K_0(X)_{\mathbb C}$) its Grothendieck group (resp. complexified),
\item $ K_0^{G}(X)$ (resp. $K_0^G(X)_{\mathbb C}$) its $G$-equivariant Grothendieck group (resp. complexified).
\end{enumerate}
Any complex of $G$-equivariant quasi-coherent complexes admit  flat and injective resolutions. From this one can deduce that on $\mathcal D^b_G(X)$  all standard derived functors are well defined. In particular, we have a well defined left derived tensor product $\otimes \colon \mathcal D^b_G(X)\times \mathcal D^b_G(X)\to \mathcal D^b_G(X)$, and any  $f\colon X\to Y$, morphism of smooth projective $G$-varieties, induces left and right derived functors $Lf^*\colon \mathcal D^b_G(Y)\to \mathcal D^b_G(X)$ and $Rf_*\colon \mathcal D^b_G(X)\to \mathcal D^b_G(Y)$. It is possible to show that all the standard properties of the derived tensor product, the derived pull-back and push-forward functors are valid in the equivariant setting. Moreover, all these equivariant derived functors are compatible with their non-equivariant versions via the forgetful functor \footnote{For the translation of the theory of derived functors from the non-equivariant setting to the equivariant one, the reader may consult \cite[Chapter 5]{chrissginz}, \cite{blunts} for the topological setting, 
\cite[Section 1.5]{varaves},  and also \cite{lipman}.}.

The structural morphism $\pi\colon X\to {\rm Spec}(\mathbb C)$ endows $K_0(X)$ and $K_0^G(X)$ with a $\mathbb C$-algebra and an $R(G)$-algebra structures, respectively. In addition, it induces serveral push-forward morphisms 
\[\pi_*\colon K_0(X)\to  K_0({\rm Spec\ }\mathbb C)\cong\mathbb Z,\quad \pi_*\colon K_0^G(X)\to R(G).
\]and functors
\[R\pi_*\colon\mathcal D^b(X)\to\mathcal D^b(\mathbb C),\quad R\pi_*^G\colon\mathcal D^b_G(X)\to\mathcal D^b(Rep(G)),
\]
which fit into the diagram
\[\xymatrix{
\mathcal D^b_{G}(X)\ar[r]^{R\pi_*^G\quad}\ar[d]_{\frak F_X}&\mathcal D^b(Rep(G))\ar[d]^{\frak F_{{\rm pt}}}\\
\mathcal D^b(X)\ar[r]^{R\pi_*}&\mathcal D^b(\mathbb C)
}
\]where $\frak F_X,\frak F_{{\rm pt}}$ denote the forgetful functors.
If $V\in{\rm Ob}(\mathcal D^b(X))$ we call a \emph{$G$-equivariant structure} on $V$ any object $V'\in{\rm Ob}(\mathcal D^b_G(X))$ such that $\frak F_X(V')=V$.

\subsection{Equivariant Grothendieck-Euler-Poincar\'{e} characteristic}
The push-forward morphisms 
\[\pi_*\colon K_0(X)\to  K_0({\rm Spec\ }\mathbb C)\cong\mathbb Z,\quad \pi_*^G\colon K_0^G(X)\to R(G),
\]are respectively given by
\[\pi_*(V):=\sum(-1)^i{\rm rk}\ H^i(X,V)\in\mathbb Z,
\]
\[\pi_*^G(V):=\sum(-1)^i[H^i(X,V)]\in R(G),
\]where $[H^i(X,V)]$ denotes the $R(G)$-class of the cohomology space $H^i(X,V)$ seen as a representation of $G$. These morpshims define the Grothendieck-Euler-Poincar\'{e} characteristic of (the isomorphism class of) an object $V$, and its equivariant version respectively. They will be denoted by $\chi,\chi^G$:
\[\chi(V):=\pi_*(V),\quad \chi^G(V):=\pi_*^G(V).
\]
In both cases, using the duality involutions
\beq\label{dualneq}(-)^*\colon K_0(X)\to K_0(X),\quad E\mapsto E^*,
\eneq
\beq\label{dualeq}(-)^*\colon K_0^G(X)\to K_0^G(X),\quad E\mapsto E^*,
\eneq
we can define a non-symmetric paring, called the \emph{Grothendieck-Euler-Poincar\'{e} pairing}
 (or also the \emph{Mukai pairing}):
\beq\chi(E,F):=\chi(E^*\otimes F),\quad \chi^G(E,F):=\chi^G(E^*\otimes F).
\eneq These pairings naturally extend to the complexified algebras $K_0(X)_{\mathbb C}$ and $K_0^G(X)_{\mathbb C}$. In the non-equivariant case, the pairing $\chi$ is $\mathbb C$-bilinear, whereas in the equivariant case the pairing $\chi^G$ is $R(G)_{\mathbb C}$-sesquilinear wrt the duality involution naturally defined on $R(G)_{\mathbb C}$:
\beq\label{dualrap}(-)^*\colon R(G)_{\mathbb C}\to R(G)_{\mathbb C},\quad [V]\mapsto [V^*].
\eneq
That is,
 $\chi^G(\rho_1E_1,\rho_2 E_2)=\rho_1^*\rho_2\,\chi^G(E_1,E_2)$ for $E_1,E_2\in K_0^G(X)$ and $\rho_1,\rho_2\in R(G)$.
 
We consider the involutive operation on $n\times n$-matrices $$(-)^*\colon M_n(R(G)_{\mathbb C})\to M_n(R(G)_{\mathbb C}),$$ which consists in applying \eqref{dualrap} at each entry. 
For $A\in M_n(R(G)_{\mathbb C})$ we  define the matrix $A^\dag\in M_n(R(G)_{\mathbb C})$ as follows:
\beq\label{daga}
(A^\dag)_{\alpha,\beta}:=A^*_{\beta,\alpha},\quad\alpha,\beta=1,\dots,n.
\eneq

If $G=\mathbb T$, then 
 the duality involution acts on $R_G(\mathbb T)_{\mathbb C}\cong \mathbb C[\bm{Z}^{\pm 1}]$
by the formula:
\beq f(\bm{Z})^*=f(\bm{Z}^{-1}),
\eneq
where $f(\bm Z)=f(Z_1,\dots, Z_n)\in\mathbb C[\bm Z^{\pm 1}]$ and
$f(\bm Z^{-1}):=f(Z_1^{-1},\dots, Z_n^{-1})$.

\subsection{Exceptional collections in $\mathcal D^b_G(X)$ and their mutations}
Given two objects $E,F\in{\rm Ob}(\mathcal D^b_G(X))$, we define
\[\Hom^\bullet_G(E,F):=R\pi_*^G\left(E^*\otimes F\right)\in {\rm Ob}(\mathcal D^b(Rep(G))),
\]where $E^*:=R\HOM(E,\mathcal O_X)$ is the ordinary dual sheaf of $E$.
\begin{defi}\label{exccol}
An object $E\in{\rm Ob}(\mathcal D^b_G(X))$ is called an \emph{exceptional object} if and only if
\[\Hom^\bullet_G(E,E)\cong \mathbb C_G,
\]
where $\mathbb C_G$ denotes the object of $\mathcal D^b(Rep(G))$ given by the trivial complex one dimensional representation of $G$, concentrated in degree zero.

An ordered collection $(E_1,\dots, E_n)$ is said to be an \emph{exceptional collection} if and only if
\begin{itemize}
\item all objects $E_i$'s are exceptional objects,
\item and $\Hom^\bullet_G(E_j, E_i)=0$ for $j>i$.
\end{itemize}
\end{defi}
The definitions above are the natural equivariant versions of the standard notions of exceptional objects and collections in $\mathcal D^b(X)$. The following result, due to A. Elagin, gives an insight on the relationships between ordinary exceptional collections in $\mathcal D^b(X)$ and equivariant exceptional collections in $\mathcal D^b_G(X)$. Before stating Elagin's result, let us recall that there is a naturally defined operation of tensor product between objects of $\mathcal D^b_G(X)$ and $\mathcal D^b(Rep(G))$: if $E\in {\rm Ob}(\mathcal D^b_G(X))$ and $V^\bullet\in{\rm Ob}(\mathcal D^b(Rep(G)))$, the tensor product $E\otimes V^\bullet$ is defined as the object of $\mathcal D^b_G(X)$ given by
\beq\label{tensprod}
\bigoplus_{i}E[-i]\otimes V^i.
\eneq

This extends the obvious operation of tensor product between objects of $Coh_G(X)$ and $Rep(G)$. 

If $\mathcal A_1,\dots,\mathcal A_n$ are subcategories of $\mathcal D^b_G(X)$, we denote by $\left\langle \mathcal A_1,\dots,\mathcal A_n\right\rangle$ the smallest full triangulated subcategory of $\mathcal D^b_G(X)$ containing $\mathcal A_1,\dots, \mathcal A_n$.

\begin{defi}
Let $\frak E:=(E_1,\dots, E_n)$ be an exceptional collection in $\mathcal D^b_G(X)$. We  say that $\frak E$ is $G$-\emph{full} if
\beq\label{Gfull}\mathcal D^b_G(X)=\left\langle E_1\otimes\mathcal D^b(Rep(G)),\dots,  E_n\otimes\mathcal D^b(Rep(G))\right\rangle.
\eneq
\end{defi}

\begin{oss}Thus the exceptional collection $(E_1,\dots, E_n)$ is  $G$-full if and only if the collection $\left( E_1\otimes\mathcal D^b(Rep(G)),\dots,  E_n\otimes\mathcal D^b(Rep(G))\right)$ realizes a so-called \emph{semi-orthogonal decomposition} of $\mathcal D^b_G(X)$, see e.g. \cite[Chapter 1]{huy}.
\end{oss}

\begin{oss}Our definition of $G$-fullness is different from the definition of fullness of exceptional collections in triangulated categories. In the paper \cite{bororl}, L.\,Borisov and D.\,Orlov studied bounded derived category of
 $\mathbb T$-equivariant coherent sheaves on smooth toric varieties and Deligne-Mumford stacks. In particular, they described and explicitly constructed full exceptional collections in these categories. Notice that their exceptional 
 collections consist of infinite sets of objects, while we collect an infinite set of objects in one symbol $E_i\otimes \mathcal D^b(Rep(G))$.
\end{oss}

\begin{theorem}[{\cite[Theorem 2.6]{elagin}}]
\label{elagin}
Assume that $(E_1,\dots, E_n)$ is a full exceptional collection of $\mathcal D^b(X)$, and 
each object $E_i$ admits a $G$-equivariant structure $\mathcal E_i$. Then, $(\mathcal  E_1,\dots,\mathcal E_n)$ is a $G$-full exceptional collection in $\mathcal D^b_G(X)$.
\end{theorem}

Being thus important to know under which conditions an exceptional object of $\mathcal D^b(X)$ admits a $G$-equivariant structure, we recall the following result of A.\,Polishchuk.

\begin{theorem}[{\cite[Lemma 2.2]{poli}}]
\label{pol1}
Let $X$ be a smooth projective complex variety equipped with the action of a linear algebraic connected reductive group $G$ with $\pi_1(G)$ torsion free. If $E\in\mathcal D^b(X)$ is an exceptional object, then $E$ admits a $G$-equivariant structure, which is unique up to tensoring by a character of $G$.
\end{theorem}

In the present paper we focus on the case $G=\mathbb T$, and the assumption of Theorem \ref{pol1} applies.

\begin{defi}[Mutations of objects]
Let $E\in{\rm Ob}(\mathcal D^b_G(X))$ be an exceptional object. For any $F\in{\rm Ob}(\mathcal D^b_G(X))$ we  define two new objects
\[\mathbb L_EF,\quad \mathbb R_EF,
\]called the \emph{left} and \emph{right mutations} of $F$ with respect to $E$. These two objects are defined through the distinguished triangles
\beq\label{left}
\xymatrix{\mathbb L_EF[-1]\ar[r]&\Hom^\bullet_G(E,F)\otimes E\ar[r]^{\quad\quad\quad j^*}&F\ar[r]&\mathbb L_EF,}
\eneq
\beq\label{right}
\xymatrix{\mathbb R_EF\ar[r]&F\ar[r]^{j_*\quad\quad\quad}&\Hom^\bullet_G(F,E)^*\otimes E\ar[r]&\mathbb R_EF[1],}
\eneq
where $j^*,j_*$ denote the natural evaluation and coevaluation morphisms.
\end{defi}

\begin{oss}
As in the non-equivariant case, it can be shown that the objects $\mathbb L_EF,\mathbb R_EF$ are uniquely defined (up to unique isomorphism) by the distinguished triangles above. The key property is the exceptionality of $E$. We leave the details to the reader, see  \cite[Section 3.3]{CDG1}.
\end{oss}

\begin{lemma}
Let $E\in{\rm Ob}(\mathcal D^b_G(X))$ be an exceptional object. We have
\[\Hom^\bullet_G(E,\mathbb L_EF)=0,\quad \Hom_G^\bullet(\mathbb R_EF,E)=0,
\]for all objects $F\in{\rm Ob}(\mathcal D^b_G(X))$.
\end{lemma}

\proof
Apply the functor $\Hom^\bullet_G(E,-)$ (resp. $\Hom^\bullet_G(-,E)$) to the distinguished triangle \eqref{left} (resp. \eqref{right}), and use the exceptionality of $E$.
\endproof

\begin{defi}
Let $\frak E:=(E_i)_{i=1}^n$ be an exceptional collection in $\mathcal D^b_G(X)$. For any integer $i$, with $0< i<n$, we define two new collections
\begin{align*}\mathbb L_i(\frak E):&=(E_1,\dots, \mathbb L_{E_i}E_{i+1}, E_i,\dots, E_n),\\
\mathbb R_i(\frak E):&=(E_1,\dots, E_{i+1}, \mathbb R_{E_{i+1}}E_i, \dots, E_n).
\end{align*}
\end{defi}

\begin{prop}
For any $i$, with $0<i<n$, the
 collections $\mathbb L_i(\frak E), \mathbb R_i(\frak E)$ are exceptional.
  Moreover, the mutation operators $\mathbb L_i,\mathbb R_i$ satisfy the following identities:
\beq\label{br1}
\mathbb L_i\mathbb R_i=\mathbb R_i\mathbb L_i={\rm Id},
\eneq
\beq\label{br2}
\mathbb R_i\mathbb R_j=\mathbb R_j\mathbb R_i,\quad \text{if}\ \
|i-j|>1,\quad \mathbb R_{i+1}\mathbb R_i\mathbb R_{i+1}=\mathbb R_i\mathbb R_{i+1}\mathbb R_i.
\eneq
\end{prop}

\proof
The same as in the non-equivariant case, see \cite{helix}, \cite[Section 3.3]{CDG1}.
\endproof

Denote by $\tau_1,\dots, \tau_{n-1}$ the generators of the braid group $\mathcal B_n$, satisfying the relations 
\[\tau_i\tau_{i+1}\tau_i=\tau_{i+1}\tau_i\tau_{i+1},\quad \tau_i\tau_j=\tau_j\tau_i,\quad
\text{if}\ \ |i-j|>1.
\]
 We define the left action of $\mathcal B_n$ on the set of exceptional collections of length $n$ by identifying the action of  $\tau _i$ with $\mathbb R_i$,
 see identities \eqref{br1}-\eqref{br2}.
 
  For our purposes, we  modify this action, by setting 
\beq\label{actbraid}
\tau_i(\frak E):=\mathbb R_{n-i}(\frak E),\quad i=1,\dots,n-1,
\eneq
for any exceptional collection $\frak E=(E_1,\dots, E_n)$.

\begin{oss}
Formula \eqref{actbraid} is in agreement with the notations of \cite{tarvar}, see Remark \ref{chi-A}.
\end{oss}

\subsection{Dual exceptional collections and helices}
\begin{defi}[Dual exceptional collections]Let $\frak E=(E_1,\dots, E_n)$ be an exceptional collection. 
Define the \emph{left} and \emph{right dual} exceptional collections $^\vee\frak E$ and $\frak E^\vee$ as the collections 
\beq\label{braiddualcat}
^\vee\frak E:=\beta(\frak E),\quad \beta:=\tau_1(\tau_{2}\tau_{1})\dots(\tau_{n-2}\dots \tau_{1})(\tau_{n-1}\tau_{n-2}\dots \tau_{1}),
\eneq
\beq
\frak E^\vee:=\beta^{-1}(\frak E).
\eneq 
\end{defi}

\begin{prop}
Let $\frak E=(E_1,\dots, E_n)$ be an exceptional collection, $^\vee\frak E=(^\vee E_1,\dots ^\vee E_n)$ and $\frak E^\vee=(E_1^\vee,\dots, E_n^\vee)$ its left and right dual exceptional collections, respectively. The following orthogonality relations hold true:
\begin{empheq}[left={\Hom^\bullet_G(E_h,E_k^\vee)=}\empheqlbrace]{align*} 
\mathbb C_G,\quad & h=n-k+1,\\
\quad &\\
 0,\quad &\text{otherwise},
\end{empheq}

\begin{empheq}[left={\Hom^\bullet_G(^\vee E_k, E_h)=}\empheqlbrace]{align*} 
\mathbb C_G,\quad & h=n-k+1,\\
\quad &\\
0,\quad &\text{otherwise}.
\end{empheq}
Moreover, for any $F\in {\rm Ob}\left(\mathcal D^b_G(X)\right)$ we have the functorial isomorphism
\beq\label{serredualright}
\Hom_G^\bullet(^\vee E_k,F)\cong\Hom_G^\bullet(F,E_k^\vee)^*.
\eneq
\end{prop}

\proof
The argument is the same as in the non-equivariant case, see  \cite[Section 3.6]{CDG1}.
\endproof

Given an exceptional collection $\frak E$, we  introduce the infinite family of exceptional objects called the \emph{helix} generated by $\frak E$.
\begin{defi}[Helix]
Let $\frak E=(E_1,\dots, E_n)$ be an exceptional collection. Define the \emph{helix} generated by $\frak E$ to be the infinite family of objects $(E_i)_{i\in\mathbb Z}$ defined by the iterated mutations
\[E_{i+n}:=\mathbb R_{E_{i+n-1}}\dots\mathbb R_{E_{i+1}}E_i,\quad E_{i-n}:=\mathbb L_{E_{i-n+1}}\dots\mathbb L_{E_{i-1}}E_i,\quad i\in\mathbb Z.
\]Such a helix is said to be of \emph{period $n$}. Any family of $n$ consecutive objects $(E_i,\dots, E_{i+n})$ is called a \emph{foundation} of the helix.
\end{defi}

\subsection{Exceptional bases in equivariant $K$-theory}In this Section we focus on the $K$-theoretical counterpart of the notion of exceptional collections introduced in Definition \ref{exccol} and of the action of the braid group on them.

\begin{defi}\label{excbas}An element $e\in K_0^G(X)$ is \emph{exceptional} if 
 $$\chi^G(e,e)=\mathbb C_G.$$ A basis $\varepsilon:=(e_i)_{i=1}^n$ of $K_0^G(X)$ as an
  $R(G)$-module, is  \emph{exceptional} if 
\beq
\chi^G(e_i,e_i)=\mathbb C_G,\quad \chi^G(e_j,e_i)=0,\quad \text{for } j>i.
\eneq 
\end{defi}

The following result is a $K$-theoretical analogue of Theorem \ref{elagin}.

\begin{theorem}[{\cite[ Lemma 2.1]{poli}}]
\label{pol2}
Let $(E_1,\dots, E_n)$ be a full exceptional collection in $\mathcal D^b(X)$. If each object $E_i$ admits a $G$-equivariant structure, then the classes $([E_i])_{i=1}^n$ form an exceptional basis of $K_0^G(X)$ as an
$R(G)$-module.
\end{theorem}

\begin{prop}
Let $E\in{\rm Ob}(\mathcal D^b_G(X))$ be an exceptional object. For any $F\in {\rm Ob}(\mathcal D^b_G(X))$ we have
\beq
[\mathbb L_EF]=[F]-\chi^G(E,F)\cdot [E],\quad [\mathbb R_EF]=[F]-\chi^G(F,E)^*\cdot [E].
\eneq
\end{prop}

\proof
From the distinguished triangle \eqref{left}, and  equation \eqref{tensprod}, we deduce 
\begin{align*}
[L_EF]&=[F]-[\Hom^\bullet_G(E,F)\otimes E]=[F]-\left[\bigoplus_i E[-i]\otimes H^i(X,E^*\otimes F)\right]\\
&=[F]-\left(\sum_i(-1)^iH^i(X,E^*\otimes F)\right)[E].
\end{align*}
Analogously, from the distinguished triangle \eqref{right}, we deduce 
\begin{align*}
[R_EF]&=[F]-[\Hom^\bullet_G(F,E)^*\otimes E]=[F]-\left[\bigoplus_i E[-i]\otimes H^i(X,F^*\otimes E)^*\right]\\
&=[F]-\left(\sum_i(-1)^iH^i(X,F^*\otimes E)\right)^*[E].
\end{align*}This completes the proof.
\endproof

\begin{defi}
Let $e\in K_0^G(X)$ be an exceptional element. Given $f\in K_0^G(X)$, we define its \emph{left} and \emph{right mutations} wrt $e$ as the elements
\beq
\mathbb L_{e}f:=f-\chi^G(e,f)\cdot e,\quad \mathbb R_{e}f:=f-\chi^G(f,e)^*\cdot e.
\eneq
\end{defi}

\begin{lemma}\label{lem}
Let $e\in K_0^G(X)$ be an exceptional element. We have
\[\chi^G(e,\mathbb L_ef)=0,\quad \chi^G(\mathbb R_ef,e)=0,
\]for any $f\in K_0^G(X)$.
\qed
\end{lemma}

\begin{defi}
Let $K^G_0(X)$ be a free $R(G)$-module of finite rank and  $\varepsilon:=(e_i)_{i=1}^n$ an exceptional basis of
 $K^G_0(X)$. For any $0<i<n$ define the two new exceptional bases 
\beq
\mathbb L_i\varepsilon:=(e_1,\dots,e_{i-1},\mathbb L_{e_i}e_{i+1},e_i,e_{i+2},\dots, e_n),
\eneq
\beq
\mathbb R_i\varepsilon:=(e_1,\dots,e_{i-1}, e_{i+1}, \mathbb R_{e_{i+1}}e_i,e_{i+2},\dots, e_n).
\eneq

This construction defines the action of the braid group $\mathcal B_n$ on the set of exceptional bases of $K^G_0(X)$,
in which the action of a generator $\tau_i$,  $i=1,\dots, n-1$, is identified with the action of the mutation
 $\mathbb R_{n-i}$.
\end{defi}

\subsection{Dual exceptional bases} Let $\varepsilon:=(e_i)_{i=1}^n$ be an exceptional basis of $K_0^G(X)$. Define the \emph{left} and \emph{right dual exceptional bases}, $^\vee\varepsilon$ and $\varepsilon^\vee$, through the mutations
\beq\label{braiddual}^\vee\varepsilon:=\beta(\varepsilon),\quad \beta:=\tau_1(\tau_{2}\tau_{1})\dots(\tau_{n-2}\dots \tau_{1})(\tau_{n-1}\tau_{n-2}\dots \tau_{1}),
\eneq

\beq
\varepsilon^\vee:=\beta^{-1}\varepsilon.
\eneq

\begin{prop}\label{propconndual}
Let $\varepsilon=(e_i)_{i=1}^n$ be an exceptional basis of $K_0^G(X)$, $^\vee\varepsilon=(^\vee e_i)_{i=1}^n$ and $\varepsilon^\vee=(e_i^\vee)_{i=1}^n$ its left and right dual exceptional basis, respectively. The following orthogonality relations hold true
\beq\label{dual0}	\chi^G(e_h,e_k^\vee)=\delta_{h+k,n+1},\quad \chi^G(^\vee e_k,e_h)=\delta_{h+k,n+1},
\eneq
for $k=1,\dots,n$. In particular, for any $v\in K_0^G(X)$ we have
\beq
\label{dual1}v=\sum_{h=1}^n\chi^G(v,e_h^\vee)^*\ e_{n+1-h},\quad
v=\sum_{h=1}^n\chi^G(^\vee e_h, v)\ e_{n+1-h}.
\eneq
\end{prop}

\proof
We prove the first identity in (\ref{dual0}), the proof of the second is analogous. We have $$\chi^G(e_h,e^\vee_k)=0,\quad\text{for }h=1,\dots, n-k,$$ 
by Lemma \ref{lem}.

If $e,f\in K^G_0(X)$,  $e$ is exceptional and $\chi^G(f,e)=0$, then  $\chi^G(f,v)=\chi^G(f,\mathbb L_ev)$
for any $v\in K^G_0(X)$.  By iteration of this identity, we deduce 
\[\chi^G(e_h,e^\vee_k)=0,\quad\text{for }h=n-k+2,\dots, n,
\]and 
\[\chi^G(e_{n-k+1},e^\vee_{k})=\chi^G(e_{n-k+1},e_{n-k+1})=1.
\]
Identities \eqref{dual1} follow from the sesquilinearity of $\chi^G$.
\endproof

\begin{cor}\label{gramdual}
Let $\varepsilon=(e_i)_{i=1}^n$ be an exceptional basis of $K_0^G(X)$, and 
$\mathcal G$ the Gram matrix of $\chi^G$ wrt $\varepsilon$. 
Then the Gram matrix of $\chi^G$ wrt  $^\vee\varepsilon$ equals the Gram matrix of $\chi^G$ wrt
$\varepsilon^\vee$ 
and
equals \[J\cdot (\mathcal G^{\dag})^{-1}\cdot J,\quad\text{where }J_{\alpha,\beta}=\delta_{\alpha+\beta,n+1}.
\]
\end{cor}

\proof
Let $X=(X^j_k)_{j,k=1}^n$ be the matrix defined by $e^\vee_k:=\sum_{j=1}^nX^j_ke_j$.
Then $X$ satisfies the equation 
$\mathcal GX=J$ by formula \eqref{dual0}. 
Hence the Gram matrix of $\chi^G$ wrt $\varepsilon^\vee$ equals
\[
X^\dag\cdot \mathcal G\cdot X=J\cdot (\mathcal G^{\dag})^{-1}\cdot J. 
\]
The case of $^\vee\varepsilon$ is analogous.
\endproof

\subsection{Serre functor and canonical operator} A Serre functor $\mathcal K\colon\mathcal D^b(X)\to\mathcal D^b(X)$ is a functor defined (uniquely up to canonical isomorphism) by the condition
\beq\label{serredualnon-eq}
\Hom^\bullet(E,F)^*\cong \Hom^\bullet(F,\mathcal K(E)),\quad E,F\in{\rm Ob}(\mathcal D^b(X)).
\eneq
We can take
\beq\label{serneq}\mathcal K=(\omega_X\otimes-)[\dim_\mathbb  CX]\colon\mathcal D^b(X)\to\mathcal D^b(X),
\eneq where $\omega_X$ denotes the canonical sheaf of $X$. 
Analogously, in the equivariant case a Serre functor $\mathcal K\colon\mathcal D^b_G(X)\to\mathcal D^b_G(X)$ is defined by the condition
\beq\label{serredual}
\Hom_G^\bullet(E,F)^*\cong \Hom_G^\bullet(F,\mathcal K(E)),\quad E,F\in{\rm Ob}(\mathcal D^b_G(X)).
\eneq
We can take
\beq\label{ser}\mathcal K=(\omega^G_X\otimes-)[\dim_\mathbb  CX]\colon\mathcal D^b_G(X)\to\mathcal D^b_G(X),
\eneq where $\omega_X^G$ is the $G$-equivariant canonical sheaf of $X$. By abuse of language, we will call \eqref{serneq} (and its equivariant version \eqref{ser}) \emph{the} Serre functor in $\mathcal D^b(X)$ (and $\mathcal D^b_G(X)$, respectively).

The Serre duality \eqref{serredual} implies the \emph{Serre periodicity},
\beq\label{serreper}
\Hom^\bullet_G(E,F)\cong\Hom^\bullet_G(\mathcal K(E),\mathcal K(F)),\quad E,F\in{\rm Ob}(\mathcal D^b_G(X)).
\eneq

\begin{prop}\label{serreexc}
Let $\frak E=(E_i)_{i=1}^n$ be an exceptional collection of length $n$ of $\mathcal D^b_G(X)$. The following operations are equivalent, i.e. produce the same exceptional collection when applied to $\frak E$:
\begin{enumerate}
\item to act on $\frak E$ with the braid $(\tau_1\dots\tau_{n-1})^{-n}$,
\item to take the double right-dual exceptional collection $(\frak E^\vee)^\vee$,
\item to apply the Serre functor to each object of $\frak E$. 
\end{enumerate}
\end{prop}

\proof
The equivalence of points (1) and (2) follows from the well-known identity of braids in $\mathcal B_n$
\beq\label{center}
(\tau_1\dots\tau_{n-1})^n=\beta^2,
\eneq
where $\beta:=\tau_1(\tau_{2}\tau_{1})\dots(\tau_{n-2}\dots \tau_{1})(\tau_{n-1}\tau_{n-2}\dots \tau_{1})$ is the braid that appears in \eqref{braiddualcat} and \eqref{braiddual}, see \cite[Theorem 1.24]{braids}.  The equivalence of (2) and (3) follows from the functorial isomorphism \eqref{serredualright}.
\endproof

\begin{oss}
Note that the element \eqref{center} of $\mathcal B_n$ is the generator of its center $Z(\mathcal B_n)$, see 
\cite[Theorem 1.24]{braids}.
\end{oss}

The $K$-theoretical version of the Serre functor is the so-called \emph{canonical operator} $k\colon K^G_0(X)\to K^G_0(X)$, defined through the identity
\beq\label{cano}
\chi^G(e,f)^*=\chi^G(f,k(e)),\quad e,f\in K^G_0(X).
\eneq
The $K$-theoretical analog of \eqref{serreper}, i.e.
\beq
\chi^G(e,f)=\chi^G(k(e),k(f)),\quad e,f\in K^G_0(X),
\eneq 
shows that the canonical operator $k$ is a $\chi^G$-isometry.

\begin{prop}\label{matrk}
Let $\varepsilon=(e_i)_{i=1}^n$ be a basis of $K_0^G(X)$, and  $\mathcal G$  the Gram matrix of $\chi^G$ wrt $\varepsilon$. Then the matrix of the canonical operator $k\colon K_0^G(X)\to K_0^G(X)$ wrt the basis $\varepsilon$ is equal to
\beq
\mathcal G^{-1}\mathcal G^\dag.
\eneq
\end{prop}

\proof
It follows from identity \eqref{cano}, written in matrix notation.
\endproof

\section{Equivariant derived category, exceptional collections and $K$-theory of $\mathbb P^{n-1}$}\label{sec3}

\subsection{Symmetric functions} Consider the algebra $\mathbb C[\bm Z^{\pm 1}]=\mathbb C[ Z_1^{\pm 1},\dots,Z_n^{\pm 1}]$ of Laurent polynomials in $n$ indeterminates. The \emph{elementary} and \emph{complete} symmetric functions are defined as the elements
\beq
s_k(\bm Z):=\sum_{1\leq i_1<\dots<i_k\leq n}\prod_{j=1}^kZ_{i_j},\quad k=1,\dots,n,
\eneq
\beq
m_k(\bm Z):=\sum_{\substack{i_1\geq 0,\dots i_n\geq 0\\ i_1+\dots+i_n=k}}Z_1^{i_1}\dots Z_n^{i_n},\quad k\in\mathbb Z_{>0}.
\eneq
Put  $s_0=1$, $m_0=1$. We have 
\beq
\sum_{i=0}^k(-1)^im_i(\bm Z)s_{k-i}(\bm Z)=0,\quad k\in\mathbb Z_{>0}.
\eneq

\subsection{Torus action}\label{torusaction}Let $n\geq 2$. 
Consider the diagonal action of $\mathbb T=(\mathbb C^*)^n$ on the space $\mathbb C^n$. Such an action  induces an action of $\mathbb T$ on $	\mathbb P^{n-1}$, the projective space parametrizing the one dimensional subspaces $F\subset \mathbb C^n$. If $(u_1,\dots, u_n)$ denote the standard basis of $\mathbb C^n$, denote by $pt_I\in\mathbb P^{n-1}$, with $I=1,\dots, n$, the point corresponding to the coordinate line spanned by $u_I$. The points $pt_I$, $I=1,\dots, n$ are the fixed points of the $\mathbb T$-action.

\subsection{Derived category}The action of $\mathbb T$ on $\mathbb C^n$ induces naturally a $\mathbb T$-structure on the structural sheaf $\mathcal O_{\mathbb P^{n-1}}$ and the tautological line bundle $\mathcal O(-1)$ on $\mathbb P^{n-1}$. Any vector bundle obtained from $\mathcal O_{\mathbb P^{n-1}}$ and $\mathcal O(-1)$ through tensorial operations inherits a ``natural'' $\mathbb T$-structure. 

The derived category $\mathcal D^b(\mathbb P^{n-1})$ admits a well-known full exceptional collection, the \emph{Beilinson exceptional collection}
\[\frak B:=(\mathcal O,\mathcal O(1),\dots,\mathcal O(n-1)).
\]
Such an exceptional collection, with its natural $\mathbb T$-structure, is an exceptional collection in $\mathcal D^b_\mathbb T(\mathbb P^{n-1})$.
Its $K$-theoretical counterpart $([\mathcal O(i-1)])_{i=1}^n$ defines an exceptional basis of $K_0^{\mathbb T}(\mathbb P^{n-1})$ (in accordance with Theorem \ref{elagin}, Theorem \ref{pol1} and Theorem \ref{pol2}). 

\subsection{Equivariant $K$-theory }The equivariant $K$-theory algebra 
$K_0^{\mathbb T}(\mathbb P^{n-1})_{\mathbb C}$ admits the following presentation
\beq\label{Kpn}
K_0^{\mathbb T}(\mathbb P^{n-1})_{\mathbb C}\cong{\mathbb C[X^{\pm 1},\bm{Z}^{\pm 1}]}\Bigg/{\langle\prod_{j=1}^n(X-Z_j)\rangle},
\eneq
where the variable $X$ corresponds to the tautological line bundle $\mathcal O(-1)$ over $\mathbb P^{n-1}$, and the variables $Z_1,\dots, Z_n$ are the equivariant parameters corresponding to the factors of the torus $\mathbb T=(\mathbb C^*)^n$. 

Under the presentation \eqref{Kpn}, the duality involution \eqref{dualeq} is given by
\[f(X,\bm Z)^*:=f(X^{-1},\bm Z^{-1}),\quad f\in K_0^{\mathbb T}(\mathbb P^{n-1})_{\mathbb C}.
\]

The equivariant Grothendieck-Euler-Poincar\'{e} pairing $\chi^\mathbb T$ on
$K_0^{\mathbb T}(\mathbb P^{n-1})_{\mathbb C}$ is given by the formula
\beq\label{chitproj}
\chi^{\mathbb T}(f,g)=\sum_{a=1}^n\frac{f(Z_a^{-1},{\bm Z}^{-1})g(Z_a,{\bm Z})}{\prod_{j\neq a}(1-Z_a/Z_j)}=-\sum_{a=1}^n{\rm Res}_{X=Z_a}\frac{f(X^{-1},{\bm Z}^{-1})g(X,{\bm Z})}{X\prod_{j=1}^n(1-X/Z_j)},
\eneq
by the Atiyah-Bott equivariant localization theorem \cite{atbot}.

\begin{oss}
By putting $Z_i=1$, for $i=1,\dots,n$, in \eqref{Kpn} and \eqref{chitproj}, we obtain the presentation of the non-equivariant $K$-theory of $\mathbb P^{n-1}$ and its non-equivariant Grothendieck-Euler-Poincar\'e pairing.
\end{oss}

The class of the $\mathbb T$-equivariant canonical sheaf $[\omega^{\mathbb T}_{\mathbb P^{n-1}}]$ is obtained by twisting the class $X^n=[\mathcal O(-n)]$ with a character of $\mathbb T$:
\beq\label{cansh}[\omega^{\mathbb T}_{\mathbb P^{n-1}}]=\frac{X^n}{\prod_{j=1}^nZ_j}\quad \text{in }K^\mathbb T_{0}(\mathbb P^{n-1}).
\eneq

\begin{lemma}\label{chiO}
For $i,j\in\mathbb Z$ we have 
\begin{empheq}[left={\chi^{\mathbb T}([\mathcal O(i)],[\mathcal O(j)])=}\empheqlbrace]{align*} 
m_{j-i}(\bm Z^{- 1}),\quad & i\leq j,\\
0,\quad &j<i<j+n,\\
(-1)^{n-1}m_{i-j-n}(\bm Z)\prod_{j=1}^nZ_j,\quad &i\geq j+n.
\end{empheq}
\qed
\end{lemma}

\begin{oss}\label{chi-A}
In  \cite{tarvar}, instead of the  pairing $\chi^{\mathbb T}$ on $K^\mathbb T_{0}(\mathbb P^{n-1})_{\mathbb C}$, it is studied another non-symmetric pairing $A$ defined by the formula
\beq
A(f,g):=\pi_*\left(f^*\cdot g\cdot (-1)^{n-1}\frac{X^n}{\prod_{j=1}^nZ_j}\right),
\eneq 
where $f,g\in K^\mathbb T_{0}(\mathbb P^{n-1})_{\mathbb C}$ and 
$\pi\colon \mathbb P^{n-1}\to {\rm Spec}(\mathbb C)$. 
In \cite[Section 6]{tarvar}, a notion of exceptional bases of 
$K^\mathbb T_{0}(\mathbb P^{n-1})_{\mathbb C}$ wrt the pairing $A$, analogous to Definition \ref{excbas}, 
is given. 
From \eqref{ser}, \eqref{cano} and \eqref{cansh}
we deduce the following relationships between $A$ and $\chi^G$:
\beq
A(f,g)=\chi^{\mathbb T}(g,f)^*.
\eneq
This implies, in particular, that $A$-exceptional bases of $K^\mathbb T_{0}(\mathbb P^{n-1})_{\mathbb C}$ are exactly $\chi^{\mathbb T}$-exceptional bases, although ordered \emph{in the opposite order}. Moreover, if we denote by $\mathbb L^A,\mathbb R^A$ (resp. $\mathbb L^{\chi^\mathbb T},\mathbb R^{\chi^\mathbb T}$) the morphisms of left/right mutations wrt $A$ (resp. $\chi^\mathbb T$), then  
\beq
\mathbb L^A=\mathbb R^{\chi^\mathbb T},\quad \mathbb R^A=\mathbb L^{\chi^\mathbb T}.
\eneq 
\end{oss}

\subsection{Diophantine constraints on Gram matrices.} In this section, we show that the Gram matrices of 
 $\chi^\mathbb T$ wrt exceptional bases of $K_0^{\mathbb T}(\mathbb P^{n-1})$ satisfy certain
 Diophantine constraints.

Given $\mathcal G\in GL(n,\mathbb Z[\bm Z^{\pm 1}])$, denote
\beq\label{pg}
p_{\mathcal G}(\lambda):=\det\left(\lambda\cdot\mathbbm 1-\mathcal G^{-1}\mathcal G^\dag\right)\in\mathbb Z[\bm Z^{\pm 1},\lambda].
\eneq

\begin{lemma}We have
\beq
p_{\mathcal G^{-1}}(\lambda)=p_{\mathcal G}(\lambda)^*,
\eneq
where for any $f(\bm Z,\lambda)\in\mathbb Z[\bm Z^{\pm 1},\lambda]$ we define $f(\bm Z,\lambda)^*:=f\left(\bm Z^{-1},\lambda\right)$. 
\end{lemma}
\proof
Notice that $p_{\mathcal G^{-1}}(\lambda)=\det\left(\left(\lambda\cdot\mathbbm 1-\mathcal G^{-1}\mathcal G^\dag\right)^\dag\right)=p_{\mathcal G}(\lambda)^*$.
\endproof

\begin{theorem}\label{propdioph}
Let $\varepsilon:=(e_i)_{i=1}^n$ be a basis of $K_0^\mathbb T(\mathbb P^{n-1})$, and let $\mathcal G$ be the Gram matrix of $\chi^\mathbb T$ wrt $\varepsilon$. The following identity holds true:
\beq
\label{dioph1}
p_{\mathcal G}(\lambda)=\sum_{j=0}^n(-1)^j\lambda^{n-j}s_j\left((-1)^{n-1}\frac{Z_1^n}{s_n(\bm Z)},\dots,(-1)^{n-1}\frac{Z_n^n}{s_n(\bm Z)}\right).
\eneq
\end{theorem}

\proof
From presentation \eqref{Kpn} and 
equation \eqref{cansh}, it is readily seen that the eigenvalues of the canonical operator $k$ are 
\[(-1)^{n-1}\frac{Z_1^n}{s_n(\bm Z)},\dots,(-1)^{n-1}\frac{Z_n^n}{s_n(\bm Z)}.
\]Then, identity \eqref{dioph1} follows from Proposition \ref{matrk}.
\endproof

If we expand \eqref{pg} in powers of $\lambda$, i.e.
\beq
p_\mathcal G(\lambda)=\sum_{j=0}^n(-1)^j\lambda^{n-j}p_j(\mathcal G),
\eneq
for suitable polynomial functions $p_j(\mathcal G)$ of the entries of $\mathcal G$ and $\mathcal G^\dag$, from the identity \eqref{dioph1} we deduce the validity of $n$ constraints:
\beq\label{dioph1.2}
p_j(\mathcal G)=s_j\left((-1)^{n+1}\frac{Z_1^n}{s_n(\bm Z)},\dots,(-1)^{n+1}\frac{Z_n^n}{s_n(\bm Z)}\right),\quad j=1,\dots,n.
\eneq
If $\mathcal G$ is a Gram matrix of $\chi^{\mathbb T}$, then $\det \mathcal G=1$, and we have
\beq
p_\mathcal G\left(\frac{1}{\lambda}\right)=\frac{(-1)^n}{\lambda^n}p_{\mathcal G}(\lambda)^*,
\eneq
so that 
\beq
p_{n-j}(\mathcal G)=p_j(\mathcal G)^*,\quad j=0,\dots, n.
\eneq
Thus, we are left with $[\frac{n}{2}]$ constraints  for the entries of $\mathcal G$. Let us write the constraints for $n=2,3,4$.

\begin{prop}
Let $\varepsilon:=(e_1,e_2)$ be an exceptional basis of $K_0^{\mathbb T}(\mathbb P^1)$, and let
\[\mathcal G=\begin{pmatrix}
1&g\\
0&1
\end{pmatrix},\quad g\in\mathbb Z[Z_1^{\pm1},Z_2^{\pm 1}],
\]be the Gram matrix of $\chi^\mathbb T$ wrt $\varepsilon$. Then, the Laurent polynomial $g$ is a solution of the equation
\beq\label{diophP1}
gg^*=\frac{(Z_1+Z_2)^2}{Z_1Z_2}.
\eneq
All the solutions of \eqref{diophP1} are of the form
\beq\label{solp1}
g(Z_1,Z_2)=Z_1^\alpha Z_2^\beta(Z_1+Z_2)\in\mathbb Z[Z_1^{\pm1},Z_2^{\pm 1}],
\eneq
where $\alpha,\beta\in\mathbb Z$.
\end{prop}

\proof
We have $p_{\mathcal G}(\lambda)=\lambda^2+\left(gg^*-2\right)\lambda+1$, and the only non-trivial constraint \eqref{dioph1.2} is
\beq
gg^*-2=\frac{Z_1^2+Z_2^2}{Z_1Z_2},
\eneq
which coincides with \eqref{diophP1}. Notice that $g$ is a solution of \eqref{diophP1} if and only if $\gamma:={g}\cdot {s_1(Z_1,Z_2)}^{-1}$ is a solution of 
$\gamma\gamma^*=1,$
whose solutions are $\gamma(Z_1,Z_2)=Z_1^\alpha Z_2^\beta,$ with $\alpha,\beta\in\mathbb Z.$
\endproof

\begin{oss}
By Lemma \ref{chiO}, the matrix $\mathcal G$ corresponding to the solution \eqref{solp1} coincide with the Gram matrix wrt the exceptional basis $\left([p\otimes \mathcal O],\quad [q\otimes \mathcal O(1)]\right),$ where $p,q\in R(\mathbb T)\cong\mathbb  Z[Z_1^{\pm1},Z_2^{\pm 1}]$ are characters of $\mathbb T$ such that $p^*q=Z_1^{\alpha+1}Z_2^{\beta+1}.$
\end{oss}

\begin{prop}
Let $\varepsilon=(e_1,e_2,e_3)$ be an exceptional basis of $K_0^{\mathbb T}(\mathbb P^2)$, and let
\[\mathcal G=\begin{pmatrix}
1&a&b\\
0&1&c\\
0&0&1
\end{pmatrix},\quad a,b,c\in\mathbb Z[Z_1^{\pm1},Z_2^{\pm1},Z_3^{\pm1}],
\]
be the Gram matrix of $\chi^\mathbb T$ wrt $\varepsilon$. Then, the triple $(a,b,c)$ is a  solution of the Markov-type equations
\begin{align}\label{markov1}
aa^*+bb^*+cc^*-ab^*c&=3-\frac{Z_1^3+Z_2^3+Z_3^3}{Z_1Z_2Z_3},\\
\label{markov2}
aa^*+bb^*+cc^*-a^*bc^*&=3-\frac{Z_2^3Z_3^3+Z_1^3Z_3^3+Z_2^3Z_3^3}{Z_1^2Z_2^2Z_3^2}.
\end{align}
\end{prop}

Notice that the triple $(a,b,c)=(s_1(\bm Z), s_2(\bm Z), s_1(\bm Z))$ gives a solutions of 
\eqref{markov1} and \eqref{markov2}.

The properties of the Markov-type 
equations \eqref{markov1}-\eqref{markov2} and its solutions are discussed  in \cite{cv-inprep}. 

\begin{oss}

In the non-equivariant case, the Gram matrices wrt exceptional bases are upper triangular matrices 
with ones on the diagonal and integer  
entries $(a,b,c)$ satisfying the equation 
\beq
\label{mMa}
a^2+b^2+c^2-abc=0,
\eneq
 see 
\cite{helix, bondalsympl}. This Diophantine equation is equivalent to 
 the famous Markov 
equation,
\beq
\label{Mark}
a^2+b^2+c^2 - 3abc=0,
\eneq
 see \cite{aigner}. A triple of integers  $(a,b,c)$ is a solution of  (\ref{Mark}) if and only if the triple
 of integers 
  $(3a,3b, 3c)$ is a solution of   (\ref{mMa}).
 
 If $a,b,c\in\mathbb Z[Z_1^{\pm1},Z_2^{\pm1},Z_3^{\pm1}]$ is a solution of equations 
 \eqref{markov1} and \eqref{markov2}, then putting  $Z_1=Z_2=Z_3=1$ in the Laurent polynomials
 $a,b,c$ we obtain a triple of integers satisfying the Markov equation (\ref{Mark}). For example, the solution 
 $$(s_1(\bm Z), s_2(\bm Z), s_1(\bm Z))$$   gives the minimal Markov triple $(3,3,3)$.
 Thus equations 
 \eqref{markov1} and \eqref{markov2} may  be considered as a Laurent polynomial deformation of the classical Markov equation.
 \end{oss}

\begin{prop}
Let $\varepsilon=(e_1,e_2,e_3,e_3)$ be an exceptional basis of $K_0^{\mathbb T}(\mathbb P^3)$, and let
\[\mathcal G=\begin{pmatrix}
1&a&b&c\\
0&1&d&e\\
0&0&1&f\\
0&0&0&1
\end{pmatrix},\quad a,b,c,d,e,f\in\mathbb Z[Z_1^{\pm1},Z_2^{\pm1},Z_3^{\pm1},Z_4^{\pm 1}],
\]
be the Gram matrix of $\chi^\mathbb T$ wrt $\varepsilon$. Then, $(a,b,c,d,e,f)$ is a  solution of the equations
\begin{align}\nonumber
aa^*&+bb^*+cc^*+dd^*+ee^*+ff^*\\
\label{mark4n1}
&-a^*bd^*-a^*ce^*-b^*cf^*-d^*ef^*+a^*cd^*f^*\\
\nonumber
&=4+\frac{Z_2 Z_4 Z_3}{Z_1^3}+\frac{Z_1 Z_4 Z_3}{Z_2^3}+\frac{Z_1 Z_2 Z_4}{Z_3^3}+\frac{Z_1 Z_2 Z_3}{Z_4^3},\\
\nonumber\quad& 
\end{align}
\begin{align}
\nonumber -2 a a^*&-2 b b^*-2 c c^*-2 d d^*-2 e e^*-2 f f^*\\
\nonumber &+a b^* d+a^* b d^*+a c^* e+a^* c e^*+b^* c f^*+b c^* f+d e^* f+d^* e f^*\\
\nonumber&-a b^* e f^*-a^* b e^* f-b c^* d^* e-b^* c d e^*\\
\label{mark4n2}
&+a a^* f f^*+b b^* e e^*+c c^* d d^*\\
\nonumber&=-6+\frac{Z_4^2 Z_2^2}{Z_1^2 Z_3^2}+\frac{Z_3^2 Z_2^2}{Z_1^2 Z_4^2}+\frac{Z_1^2 Z_2^2}{Z_3^2 Z_4^2}+\frac{Z_3^2 Z_4^2}{Z_1^2 Z_2^2}+\frac{Z_1^2 Z_4^2}{Z_2^2 Z_3^2}+\frac{Z_1^2 Z_3^2}{Z_2^2 Z_4^2},\\
\nonumber&\quad
\end{align}
\begin{align}
\nonumber aa^*&+bb^*+cc^*+dd^*+ee^*+ff^*\\
\label{mark4n3}
&-ab^*d-ac^*e-bc^*f-de^*f+ac^*df\\
\nonumber&=4+\frac{Z_1^3}{Z_2 Z_3 Z_4}+\frac{Z_2^3}{Z_1 Z_3 Z_4}+\frac{Z_3^3}{Z_1 Z_2 Z_4}+\frac{Z_4^3}{Z_1 Z_2 Z_3}.
\end{align}
\end{prop}

\begin{oss}
In the corresponding 
non-equivariant case, these Diophantine constraints on the Gram matrices wrt exceptional collection reduce to the equations
\begin{align}
a^2+b^2+c^2+d^2+e^2+f^2-abd-ace&-bcf-def+acdf=8,\\
(af-be+cd)^2&=16,
\end{align}
on the integers $(a,b,c,d,e,f)$,  see e.g. \cite{bondalsympl}. These constraints may be re-obtained by putting $Z_1=Z_2=Z_3=Z_4=1$ in 
equations \eqref{mark4n1}-\eqref{mark4n3}.
\end{oss}

\section{Equivariant cohomology of $\mathbb P^{n-1}$}\label{sec4}

\subsection{Equivariant cohomology}  Consider the $\mathbb T$-equivariant cohomology algebra 
$H^\bullet_{\mathbb T}(\mathbb P^{n-1},\mathbb C)$. Denote 
\begin{itemize}
\item by $x$ the first equivariant Chern class of the tautological line bundle $\mathcal O(-1)$ on $\mathbb P^{n-1}$ with its standard $\mathbb T$-structure,
\item by $\bm y=(y_1,\dots, y_{n-1})$ the equivariant Chern roots of the quotient bundle $\mathcal Q$ (if $F\subset \mathbb C^n$ is the line represented by $p\in\mathbb P^{n-1}$, then the fiber $\mathcal Q_p$ is the quotient $\mathbb C^n/F$),
\item by $\bm z=(z_1,\dots, z_n)$ the equivariant parameters corresponding to the factors of the torus $\mathbb T$,
\item by $\Omega$ the complement in $\mathbb C^n$ to the union of the hyperplanes
\[z_i-z_j=m,\quad{i,j=1,\dots,n},\quad i\neq j,\quad m\in\mathbb Z.
\]
\end{itemize}

It is well known that
\begin{align}\label{equivring}H^\bullet_{\mathbb T}(\mathbb P^{n-1},\mathbb C)&\cong{\mathbb C[x,\bm z]}\Bigg/{\left\langle\prod_{i=1}^n(x-z_i)\right\rangle}\\
&\cong \mathbb C[x,\bm y,\bm z]^{\frak S_{n-1}}\Bigg/\left\langle(u-x)\prod_{j=1}^{n-1}(u-y_j)-\prod_{a=1}^n(u-z_a)\right\rangle,
\end{align}
where $\mathbb C[x,\bm y,\bm z]^{\frak S_{n-1}}$ is the algebra of polynomials in $x,\bm y,\bm z$ symmetric in the variables $y_1,\dots, y_{n-1}$. The equivariant cohomology $H^\bullet_{\mathbb T}(\mathbb P^{n-1},\mathbb C)$
 is a module over the ring $H^\bullet_{\mathbb T}({\rm pt},\mathbb C)\cong \mathbb C[\bm z]$. By setting all the equivariant parameter $z_i$'s to zero in \eqref{equivring}, we obtain the presentation of the classical cohomology algebra 
\[H^\bullet(\mathbb P^{n-1},\mathbb C)\cong \mathbb C[x]/\langle x^n\rangle.
\]

\subsection{Extension of scalars}
Denote by $\mathcal O_\Omega$ the ring of holomorphic functions on the domain $\Omega$.
This ring is  a module over the ring $H^\bullet_{\mathbb T}({\rm pt},\mathbb C)\cong \mathbb C[\bm z]$. Set
\beq
H_\mathbb T^\Omega(\mathbb P^{n-1}):=H^\bullet_{\mathbb T}(\mathbb P^{n-1},\mathbb C)\otimes_{H^\bullet_{\mathbb T}({\rm pt},\mathbb C)}\mathcal O_\Omega.
\eneq
A class $\alpha\in H_\mathbb T^\Omega(\mathbb P^{n-1})$ is uniquely 
determined by the restrictions $\alpha|_{pt_I}\in \mathcal O_\Omega$ at fixed points.

Following the notations of \cite{tarvar}, we will use three different bases of $H_\mathbb T^\Omega(\mathbb P^{n-1})$:
\begin{enumerate}
\item the standard basis $(1,x,\dots, x_{n-1})$, where $x_\alpha:=x^\alpha$;
\item the basis $(g_1,\dots, g_n)$ defined by
\beq\label{gbas} g_i:=\prod_{a=i+1}^n(x-z_a),\quad i=1,\dots, n-1,\quad\text{and } g_n:=1;
\eneq
\item the idempotent basis $(\Delta_1,\dots, \Delta_n)$ defined by the \emph{Lagrange inteprolating} polynomials
\beq\Delta_i:=\prod_{j\neq i}\frac{x-z_j}{z_i-z_j},\quad i=1,\dots, n,
\eneq
\end{enumerate}
We have
\begin{equation}\label{eqidem}\Delta_i\cdot \Delta_j=\delta_{ij}\Delta_i.
\end{equation}

\subsection{Poincar\'{e} pairing and $\mathcal D$-matrix} Denote by 
\beq 
\eta\colon H^\bullet_{\mathbb T}(\mathbb P^{n-1},\mathbb C)
\times H^\bullet_{\mathbb T}(\mathbb P^{n-1},\mathbb C)\to 
H^\bullet_{\mathbb T}({\rm pt},\mathbb C)\cong \mathbb C[\bm z]
\eneq
 the \emph{equivariant Poincar\'{e} metric} given by equivariant integration
\beq\label{eqpoi}\eta(v,w):=\int_{\mathbb P^{n-1}}^{eq}v\cdot w=\sum_{a=1}^n\frac{v(z_a,\bm z)w(z_a,\bm z)}{\prod_{j\neq a}(z_a-z_j)},\quad v,w\in H^\bullet_{\mathbb T}(\mathbb P^{n-1},\mathbb C).
\eneq

 The equivariant cohomology $H^\bullet_{\mathbb T}(\mathbb P^{n-1},\mathbb C)$
  with the equivariant Poincar\'{e} metric $\eta$ 
  is a Frobenius algebra over the ring $H^\bullet_{\mathbb T}({\rm pt},\mathbb C)
  \cong \mathbb C[\bm z]$: 
\beq\label{frobalg} \eta(a\cdot b,c)=\eta(a, b\cdot c),\quad a,b,c\in H^\bullet_{\mathbb T}(\mathbb P^{n-1},\mathbb C).
\eneq

By bilinearity, we extend the Poincar\'e pairing to $H^\Omega_\mathbb T(\mathbb P^{n-1})$.
The idempotent vectors are pairwise orthogonal:
\beq\eta(\Delta_i,\Delta_j)=\eta(\Delta_{i}\cdot\Delta_j,1)=\int_{\mathbb P^{n-1}}\Delta_i\delta_{ij}=\delta_{ij}{\chi_i},\quad \chi_i:=\frac{1}{\prod_{j\neq i}(z_i-z_j)}.
\eneq
Define the matrix $\mathcal D=(\mathcal D_{j\alpha})$ as the matrix of the base change
\[x_\alpha=\sum_{j=1}^n\mathcal D_{j\alpha}\Delta_j,\quad \alpha=0,\dots, n-1.
\]

\begin{lemma} We have $$\mathcal D_{j\alpha }=z_j^{\alpha},\quad \alpha=0,\dots, n-1,\quad j=1,\dots, n.$$
Thus $\mathcal D$ is the Vandermonde matrix
\[\mathcal D=
\begin{pmatrix}
1&z_1&z_1^2&\dots&z_1^{n-1}\\
1&z_2&z_2^2&\dots&z_2^{n-1}\\
\vdots&&&\vdots\\
1&z_n&z_n^2&\dots&z_n^{n-1}
\end{pmatrix}.
\]Its inverse $\mathcal D^{-1}$ is 
\begin{empheq}[left={(\mathcal D^{-1})_{\alpha j}=}\empheqlbrace]{align*} 
      (-1)^\alpha\frac{s_{\alpha+1}^j(z)}{\prod_{m\neq j}(z_j-z_m)},\quad &0\leq\alpha<n-1,\\
      \\
      \frac{1}{\prod_{m\neq j}(z_j-z_m)},\quad &\alpha=n-1,
        \end{empheq}
        where
        \[s_{k}^j(z):=\sum_{\substack{1\leq m_1<\dots<m_{n-k}\leq n\\ m_1,\dots, m_{n-k}\neq j}}z_{m_1}\dots z_{m_{n-k}}.
        \]
  \end{lemma}
  
  \proof
  The identity $x=z_1\Delta_1+\dots+z_n\Delta_n$ implies
   the identity $x_\alpha=z_1^\alpha\Delta_1+\dots z_n^\alpha\Delta_n.$
  \endproof

\begin{lemma}
Let $\eta=(\eta_{\alpha\beta})_{\alpha,\beta}$, with $$\eta_{\alpha\beta}:=\eta(x_\alpha,x_\beta),$$ be the Gram matrix of the equivariant Poincar\'{e} metric.  We have
\begin{empheq}[left={\eta_{\alpha\beta}=}\empheqlbrace]{align*} 
        0\quad&\text{if }\alpha+\beta<n-1,\\
        1\quad&\text{if }\alpha+\beta=n-1,\\
        m_{\alpha+\beta-n+1}(\bm z)\quad&\text{if }\alpha+\beta>n-1.
        \end{empheq}
\end{lemma}

\proof
It readily follows from the identity $\mathcal D^{T}\cdot \operatorname{diag}(\chi_1,\dots,\chi_n)\cdot \mathcal D=\eta.$
\endproof

\subsection{Equivariant characteristic classes}\label{charclass}Consider a $\mathbb T$-equivariant vector bundle $V$ of rank $r$ on $\mathbb P^{n-1}$, with equivariant Chern roots $\xi_1,\dots,\xi_r$. 

\begin{defi}\label{chcl1}
Define the \emph{graded equivariant Chern character} of $V$ as the characteristic class
\[{\rm Ch}_{\mathbb T}(V):=\sum_{j=1}^r\exp(2\pi \sqrt{-1}\xi_j)\in {H}^{\Omega}_{\mathbb T}(\mathbb P^{n-1}).
\]
\end{defi}

\begin{es}
For $V=\mathcal O(k)$, $k\in\mathbb Z$, the graded Chern character ${\rm Ch}_{\mathbb T}(V)$ is the class 
\[{\rm Ch}_{\mathbb T}(\mathcal O(k))=\exp(-2\pi\sqrt{-1}kx).
\]
This is the element of ${H}^{\Omega}_{\mathbb T}(\mathbb P^{n-1})$ whose restriction at the fixed point $pt_I$ is 
\linebreak
$\exp(-2\pi\sqrt{-1}kz_I)$, for $I=1,\dots,n$.
\end{es}

\begin{lemma}\label{lemmach}
Let $V\in K^{\mathbb T}_0(\mathbb P^{n-1})_{\mathbb C}$ and $Q(\bm Z)\in\mathbb C[\bm Z^{\pm 1 }]$. We have
\[{\rm Ch}_{\mathbb T}(Q(\bm Z)V)=Q(\acute{\bm Z}){\rm Ch}_{\mathbb T}(V),\quad \acute{\bm Z}:=(e^{2\pi\sqrt{-1}z_1},\dots, e^{2\pi\sqrt{-1}z_n}).
\]
\end{lemma}
\proof
By additivity it is sufficient to prove the lemma for a monomial 
$Q(\bm Z)=Z_1^{\alpha_1}\dots Z_n^{\alpha_n}$. If $(\xi_i)_{i}$ 
are the equivariant Chern roots of $V$, then $(\xi_i+\sum_{j=1}^n\alpha_j z_j)_i$ 
are the equivariant Chern roots of $Q(\bm Z)V$.
\endproof

\begin{defi}\label{chcl2}
Given any meromorphic function $F$ on $\mathbb C$, holomorphic at $0$, with Taylor expansion of the form
\[F(t)=1+\sum_{k=1}^\infty F_kt^k,
\]we define the $\widehat F$-class of $V$ as the characteristic class 
\[\widehat{F}(V):=\prod_{j=1}^r  F(\xi_j).
\]
\end{defi}

\begin{oss}
We consider only the vector bundles $V$ and functions $F(t)$ such that $\widehat{F}(V)$ 
are elements of $H^\Omega_\mathbb T(\mathbb P^{n-1})$.
\end{oss}

\begin{defi}[Gamma classes]\label{chcl3}
The $\widehat\Gamma^\pm$-classes of $V$ are defined as the characteristic classes of $V$ obtained from the Taylor expansions
$$\Gamma(1\pm t)=\exp\left(\mp \gamma t+\sum_{k=2}^\infty(\mp 1)^k\frac{\zeta(k)}{k} t^k \right).$$
\end{defi}

Definitions \ref{chcl1}, \ref{chcl2}, \ref{chcl3}  naturally extend to objects of the equivariant derived category $\mathcal D^b_{\mathbb T}(\mathbb P^{n-1})$.

\begin{oss}
If $V=T\mathbb P^{n-1}$, the $\widehat F$-class of $V$ is called the \emph{$\widehat F$-class 
of $\mathbb P^{n-1}$}. We will denote it by $\widehat F_{\mathbb P^{n-1}}$. 
Since the Chern roots of $T\mathbb P^{n-1}$ are $(y_a-x)_{a=1}^{n-1}$, such a class is given by
\[\widehat F_{\mathbb P^{n-1}}=\prod_{a=1}^{n-1} F(y_a-x).
\]
This is the class 
whose restriction at the fixed point $pt_I$ is 
$\prod_{a\neq I}F(z_a-z_I)$.
It is an element of $H^\Omega_{\mathbb T}(\mathbb P^{n-1})$ if $F(t)$ has poles
only  at points of $\mathbb Z$. This is the case for $\widehat \Gamma^\pm_{\mathbb P^{n-1}}$.

\end{oss}

\section{Equivariant quantum cohomology of $\mathbb P^{n-1}$}\label{sec5}

\subsection{Equivariant Gromov-Witten invariants}\label{seceqgw} For a given 
$d\in H_2(\mathbb P^{n-1},\mathbb Z)$ and given 
integers $g,m\geq 0$, denote by $\overline{\mathcal M}_{g,m}(\mathbb P^{n-1},d)$ 
the moduli stack of genus $g$ stable maps to $\mathbb P^{n-1}$ with degree $d$ and $m$ marked points. We
 assume that either $d>0$ or $2g+m>2$ so that $\overline{\mathcal M}_{g,m}(\mathbb P^{n-1},d)$ 
 is non-empty. The $\mathbb T$-action on $\mathbb P^{n-1}$ induces
  a $\mathbb T$-action on $\overline{\mathcal M}_{g,m}(\mathbb P^{n-1},d)$.
   Given $m$ cohomological classes 
   $$\gamma_1,\dots, \gamma_m\in H^\bullet_{\mathbb T}(\mathbb P^n,\mathbb C),$$ 
   and integers $d_1,\dots, d_n\in\mathbb Z_{\geq 0}$, we define the \emph{genus $g$, 
   degree $d$, $\mathbb T$-equivariant descendant Gromov-Witten invariants} 
   of $\mathbb P^{n-1}$ to be the polynomials
\beq
\label{eqgw}
\langle\tau_{d_1}(\gamma_1),\dots,\tau_{d_m}(\gamma_m)\rangle^{\mathbb P^{n-1},\mathbb T}_{g,n,d}:=\left(\int^{eq}_{[\overline{\mathcal M}_{g,m}(\mathbb P^{n-1},d)]^{\rm vir}_{\mathbb T}}
\prod_{j=1}^m\psi^{d_j}_{j}{\rm ev}^*_j(\gamma_j)\right)\in 
H^\bullet_{\mathbb T}({\rm pt},\mathbb C),
\eneq 
where 
\begin{itemize}
\item $[\overline{\mathcal M}_{g,m}(\mathbb P^{n-1},d)]^{\rm vir}_{\mathbb T}\in A^{\mathbb T}_{D_{\rm vir}}(\overline{\mathcal M}_{g,m}(\mathbb P^{n-1},d))$, with $D_{\rm vir}:=nd+(n-4)(1-g)+m$, is the equivariant virtual fundamental class\footnote{Its existence is ensured by the \emph{properness} of $\overline{\mathcal M}_{g,m}(\mathbb P^{n-1},d)$. From this property, it also follows that equivariant Gromov-Witten invariants are polynomials in $\bm z$, see \cite[Section 3]{liush} and references therein.},
\item the map ${\rm ev}_j\colon \overline{\mathcal M}_{g,m}(\mathbb P^{n-1},d)\to\mathbb P^{n-1}$ is the evaluation at the $j$-th marked point, which is $\mathbb T$-equivariant,
\item the classes $\psi_j\in A^{\mathbb T}_{1}(\overline{\mathcal M}_{g,m}(\mathbb P^{n-1},d))$ denote any equivariant lift of the first Chern classes of the universal cotangent line bundles $\mathcal L_j$ on $\overline{\mathcal M}_{g,m}(\mathbb P^{n-1},d)$.
\end{itemize}
We refer the interested reader to the expository article \cite{liush}, and references therein, for details. If all $d_i$'s are zero, then the polynomials above are called \emph{primary} equivariant Gromov-Witten invariants.

\subsection{Equivariant Gromov-Witten potential}\label{seceqgwp}Consider the standard basis $(x_\alpha)_{\alpha=0}^{n-1}$ of $H^\bullet_{\mathbb T}(\mathbb P^{n-1},\mathbb C)$, seen as a $H^\bullet_{\mathbb T}({\rm pt},\mathbb C)$-module. Denote by $\bm t:=(t^0,\dots, t^{n-1})$ the corresponding dual coordinates on $H^\bullet_{\mathbb T}(\mathbb P^{n-1},\mathbb C)$, so that the generic element of $H^\bullet_{\mathbb T}(\mathbb P^{n-1},\mathbb C)$ is
\[\gamma=\sum_{\alpha=0}^{n-1}t^\alpha x_\alpha.
\]
Consider the generating function $F^{\mathbb P^{n-1},\mathbb T}_0\in H^\bullet_{\mathbb T}({\rm pt},\mathbb C)[\![t^0,\dots, t^n]\!]\cong \mathbb C[\bm z][\![\bm t]\!]$, called  \emph{equivariant Gromov-Witten potential of }$\mathbb P^{n-1}$, defined by
\begin{align}\nonumber
F^{\mathbb P^{n-1},\mathbb T}_0(\bm t):=&\sum_{m=0}^\infty\sum_{d=0}^\infty\frac{1}{m!}\langle\underbrace{\gamma,\dots,\gamma}_{m\text{ times}}\rangle_{0,m,d}^{\mathbb P^{n-1},\mathbb T}\\
\label{GWpot}=&\sum_{m=0}^\infty\sum_{d=0}^\infty\sum_{\alpha_1,\dots\alpha_m=0}^{n-1}\frac{t^{\alpha_1}\dots t^{\alpha_m}}{m!}\langle x_{\alpha_1},\dots, x_{\alpha_m}\rangle_{0,m,d}^{\mathbb P^{n-1},\mathbb T}.
\end{align}

\begin{theorem}[{\cite[Theorem 3.1]{giv1}}]
\label{eqWDVV}
The function $F^{\mathbb P^{n-1},\mathbb T}_0(t)$ satisfies the $WDVV$-equations
\[\frac{\partial^3F^{\mathbb P^{n-1},\mathbb T}_0}{\partial t^\alpha\partial t^\beta\partial t^\lambda}\eta^{\lambda\mu}\frac{\partial^3F^{\mathbb P^{n-1},\mathbb T}_0}{\partial t^\mu\partial t^\delta\partial t^\nu}=\frac{\partial^3F^{\mathbb P^{n-1},\mathbb T}_0}{\partial t^\nu\partial t^\beta\partial t^\lambda}\eta^{\lambda\mu}\frac{\partial^3F^{\mathbb P^{n-1},\mathbb T}_0}{\partial t^\mu\partial t^\delta\partial t^\alpha}.
\]
\end{theorem}

\subsection{Equivariant quantum cohomology}\label{seceqc}

The \emph{big equivariant quantum product} $*$ defined by 
\beq\label{qprod}
x_\alpha* x_\beta=\sum_{\lambda,\mu}\frac{\partial^3F^{\mathbb P^{n-1},\mathbb T}_0}{\partial t^\alpha\partial t^\beta\partial t^\lambda}\eta^{\lambda\mu}x_\mu,
\eneq
defines on $H^\bullet_{\mathbb T}(\mathbb P^{n-1},\mathbb C)[\![\bm{t}]\!]$ a \emph{Frobenius algebra} structure, namely a commutative, associative algebra with unit (the element $1$) whose product is compatible with the equivariant Poincar\'{e} metric \eqref{eqpoi}, that is 
\beq\label{qFrobalg}
\eta(a*b,c)=\eta(a,b*c),\quad a,b,c\in H^\bullet_\mathbb T(\mathbb P^{n-1},\mathbb C).
\eneq
This algebra structure on $H^\bullet_{\mathbb T}(\mathbb P^{n-1},\mathbb C)[\![\bm{t}]\!]$ 
is called the 
\emph{big equivariant quantum cohomology} of $\mathbb P^{n-1}$. 
It gives an example of a \emph{formal Frobenius manifold}
 \cite[Chapter III]{manin}. 
The quantum product \eqref{qprod} is a
 \emph{deformation} of the product in classical cohomology. 
 It is customary to denote the big quantum 
 product also by $*_{\bm t}$  to emphasize its dependence on parameters $t^i$'s. 

\subsection{Quantum connection}\label{seceqqconn} The \emph{quantum connection} of the equivariant quantum cohomology of $\mathbb P^{n-1}$ is defined by the formula
\beq
\nabla_{\alpha,\kappa}^{\rm quant}\colon H^\bullet_{\mathbb T}(\mathbb P^{n-1},\mathbb C)[\![\bm t]\!]\to H^\bullet_{\mathbb T}(\mathbb P^{n-1},\mathbb C)[\![\bm t]\!],\quad \alpha=0,\dots, n-1,
\eneq
\beq\label{qconn}
\nabla^{\rm quant}_{\alpha,\kappa}:=\kappa\frac{\partial}{\partial t^\alpha}-x_\alpha*_{\bm t},\quad 
\eneq 
 where $\kappa\in\mathbb C^*$ is the \emph{spectral parameter}. 
 The associativity of the quantum multiplication $*_{\bm t}$, i.e. Theorem \ref{eqWDVV},
  is equivalent to the flatness condition of the quantum connection 
  $\nabla^{\rm quant}_{\alpha,\kappa}$, for all $\kappa\in\mathbb C^*$:
\beq
\left[\nabla^{\rm quant}_{\alpha,\kappa},\nabla^{\rm quant}_{\beta,\kappa}\right]=0,\quad \alpha,\beta=0,\dots,n-1,\quad \kappa\in\mathbb C^*.
\eneq The system of equations for flat sections of the quantum connection is called the system of \emph{equivariant quantum differential equations}.

\begin{defi}
The $\mathbb T$-equivariant \emph{topological-enumerative morphism} is the element 
$$
\mathcal S(\bm t,\kappa)
\in{\rm End}(H^\bullet_{\mathbb T}(\mathbb P^{n-1},\mathbb C))[\![\bm t]\!][\![\kappa^{-1}]\!]
$$ 
defined by the formula
\begin{align*}\eta(\mathcal S(\bm t,\kappa)u,v)=&\ \eta(u,v)\\
+&\sum_{d=0}^\infty\sum_{m=0}^\infty\sum_{\alpha_1,\dots,\alpha_m=0}^{n-1}\frac{t^{\alpha_1}\dots t^{\alpha_m}}{m!}\Bigg\langle u,x_{\alpha_1},\dots,x_{\alpha_m},\frac{v}{\kappa-\psi}\Bigg\rangle_{0,m+2,d}^{\mathbb P^{n-1},\mathbb T},
\end{align*}where $u,v\in H^\bullet_{\mathbb T}(\mathbb P^{n-1},\mathbb C)$, $\psi:=c_1(\mathcal L_{m+2})$ and the term $\frac{1}{\kappa-\psi}$ has to be expanded in power series $\sum_{j=0}^\infty\psi^j\kappa^{-j-1}$. 
\end{defi}
\begin{defi}
The (big) $\mathbb T$-equivariant $J$-function of $\mathbb P^{n-1}$ is the cohomology-valued function defined by the identity
\beq\label{js}\eta(J(\bm t,\kappa),a)=\eta(1,\mathcal S(\bm t,\kappa)a),\quad a\in H^\bullet_{\mathbb T}(\mathbb P^{n-1},\mathbb C).
\eneq
\end{defi}

\begin{theorem}[\cite{giv1,Giv}]\label{teosoper}
For any $\kappa\in\mathbb C^*$, and any $\alpha\in H^\bullet_\mathbb T(\mathbb P^{n-1},\mathbb C)$, the cohomology class
\[\mathcal S(\bm t,\kappa)\alpha
\] is a flat section of the quantum connection
 $\nabla^{\rm quant}_{\alpha,\kappa}$, namely it satisfies the following system of differential equations
\beq	
\label{flateqs}
\kappa\frac{\partial}{\partial t^\beta}\mathcal S(\bm t,\kappa)\alpha=
x_\beta*_{\bm t}\mathcal S(\bm t,\kappa)\alpha, \quad \beta=0,\dots, n-1.
\eneq
\end{theorem}

\proof
The validity of equations \eqref{flateqs} is equivalent to the topological recursion relations in genus 0 for Gromov-Witten invariants with descendants \cite{witten}. For the proof in the non-equivariant case, see \cite{dubronapoli}, \cite[Lecture 6]{dubro1}, \cite[Lecture 2]{dubro2}, \cite[Chapter 10]{cox}, \cite[Section 7]{CDG}. For the adaptation to the equivariant case see \cite[Section 6]{giv1} and \cite[Sections 1 and 2]{Giv}.
\endproof

\subsection{Small equivariant quantum product for $\mathbb P^{n-1}$}
\begin{defi}The \emph{small} quantum product of $\mathbb P^{n-1}$ is obtained by specializing the parameters $t^i$'s of the big quantum product \eqref{qprod} as follows: $t^i=0$ for $i\neq 1$. 
\end{defi}

It is customary to put $q:=\exp(t^1)$ and to denote by $*_q$ the small quantum product. Following the notations of \cite{tarvar}, we  denote by $*_{q,\bm z}$ the small quantum product, underlining its  dependence on the equivariant parameters $\bm z$. 

A detailed study of the 
equivariant Gromov-Witten invariants of $\mathbb P^{n-1}$ (and more general flag varieties) and its small quantum cohomology can be found for example in \cite{givkim, kim, mihalcea}. For a fixed $q\in\mathbb C^*$, the small quantum product operator
\beq x*_{q,\bm z}\colon H^\bullet_\mathbb T(\mathbb P^{n-1},\mathbb C)\to H^\bullet_\mathbb T(\mathbb P^{n-1},\mathbb C),
\eneq
is the $\mathbb C[\bm z]$-linear morphism defined by the identities
\begin{align}
x*_{q,\bm z}x_j&\ =x_{j+1},\quad j=0,\dots,n-2,\\
x*_{q,\bm z}x_{n-1}
&\ =q+\sum_{i=1}^n(-1)^{i-1}s_i(\bm z)x_{n-i},
\end{align} where $s_i(\bm z)$ are the elementary symmetric polynomials in $\bm z$.

\begin{oss}
In the basis $(g_1,\dots, g_n)$, the operator $x*_{q,\bm z}$ is given by
\begin{align}
x*_{q,\bm z}g_i&\ =z_ig_i+g_{i-1},\quad i=2,\dots, n,\\
x*_{q,\bm z}g_1&\ =z_1g_1+qg_n.
\end{align}
\end{oss}

\subsection{$R$-matrices and $qKZ$ operators}For $a,b\in\left\{1,\dots, n\right\}$, with $a\neq b$, we define a family of $\mathbb C[\bm{z}]$-linear operators, called the $R$-matrices,  
\[R_{ab}(u)\colon H^\bullet_\mathbb T(\mathbb P^{n-1},\mathbb C)\to H^\bullet_\mathbb T(\mathbb P^{n-1},\mathbb C),
\]depending on a parameter $u\in\mathbb C$, and defined by the formulae
\begin{align*}
R_{ab}(u)g_i:=g_i,&\quad i\neq a,b,\\
R_{ab}(u)g_b:=g_a,&\quad R_{ab}(u)g_a:=g_b+ug_a.
\end{align*}
The $R$-matrices satisfy the Yang-Baxter equation
\[R_{ab}(u-v)R_{ac}(u)R_{bc}(v)=R_{bc}(v)R_{ac}(u)R_{ab}(u-v),
\]for $a,b,c$ all distinct, and the inversion relation
\[R_{ab}(u)R_{ba}(-u)=1.
\]Define the operators $E_1,\dots, E_n$ such that
\[E_ig_j:=\delta_{ij}g_j.
\]Define the $qKZ$ operators $K_1,\dots, K_n$ by the formula
\beq
K_i\colon H^\bullet_\mathbb T(\mathbb P^{n-1},\mathbb C)\to H^\bullet_\mathbb T(\mathbb P^{n-1},\mathbb C),
\eneq
\beq\label{qkzop} K_i:=R_{i,i-1}(z_i-z_{i-1}-1)\dots R_{i,1}(z_i-z_1-1)q^{-E_i}R_{i,n}(z_i-z_n)\dots R_{i,i+1}(z_i-z_{i+1}).
\eneq
\subsection{Equivariant $qDE$ and $qKZ$ difference equations}\label{secqkzdiff} Consider the vector bundle $H$ over the base space $\mathbb C^n$, with fiber over $\bm z_0$ given by the equivariant cohomology algebra \eqref{equivring} specialized at $\bm z=\bm z_0$, i.e.
\beq
\left.H^\bullet_{\mathbb T}(\mathbb P^{n-1},\mathbb C)\right|_{\bm z=\bm z_0}.
\eneq
Denote by ${\rm pr}\colon \mathbb C^*\times\mathbb C^n\to\mathbb C^n$ the natural projection.
Consider the pull-back vector bundle ${\rm pr}^*H$.

The quantum connection described in Section \ref{seceqqconn} defines a differential operator
\beq
\nabla_{q\frac{d}{dq},\kappa}:=\kappa q\frac{d}{dq}-x*_{q,\bm z},
\eneq 
acting on sections $I(q,\bm z)$ of the vector bundle ${\rm pr}^*H$. Following \cite{tarvar}, we fix $\kappa=1$. 

The \emph{(small) equivariant quantum differential equation} ($qDE$ for short) of $\mathbb P^{n-1}$ is the differential equation
\beq\label{eqde}
\nabla_{q\frac{d}{dq},\kappa=1}I(q,\bm z)=\left(q\frac{d}{dq}-x*_{q,\bm z}\right)I(q,\bm z)=0,
\eneq where $I$ is a section of the vector bundle ${\rm pr}^*H$. The $qDE$ is thus the equation for flat sections of ${\rm pr}^*H$.

\begin{defi}
Fix $q\in\mathbb C^*$, $\bm z,\bm z'\in\mathbb C^n$. Define the isomoprhism of vector spaces 
\beq\label{isonu}\Pi_{\bm z,\bm z'}\colon {\rm pr}^*H_{q,\bm z}\to {\rm pr}^*H_{q,\bm z'},\quad \left.x_{\alpha}\right|_{q,\bm z}\mapsto \left.x_{\alpha}\right|_{q,\bm z'},
\eneq
for $\alpha=0,\dots, n-1$.
\end{defi}

\begin{defi}
Fix $q\in\mathbb C^*$, $\bm z,\bm z'\in\mathbb C^n$. For $i=1,\dots, n$, define the isomorphisms of vector spaces 
\beq\label{isothetaqkz}
\Theta_{\bm z,\bm z'}\colon {\rm pr}^*H_{q,\bm z}\to {\rm pr}^*H_{q,\bm z'},\quad \left.g_j\right|_{q,\bm z}\mapsto \left.g_j\right|_{q,\bm z'},
\eneq
for $j=1,\dots, n$, where $g_j$'s are the elements of bases \eqref{gbas}. 
\end{defi}

For $\bm z\in\mathbb C^n$, $i=1,\dots, n$, we  use the following notations:
\begin{align}
\bm z_i^\pm:&=(z_1,\dots, z_i\pm 1,\dots, z_n)\in\mathbb C^n,\\
\Theta_{\bm z,i}^\pm :&=\Theta_{\bm z,\bm z_i^\pm}\colon {\rm pr}^*H_{q,\bm z}\to {\rm pr}^*H_{q,\bm z_i^\pm},\\
^\pm\Theta_{\bm z,i} :&=\Theta_{\bm z_i^\pm,\bm z}\colon {\rm pr}^*H_{q,\bm z_i^\pm}\to {\rm pr}^*H_{q,\bm z}.
\end{align}

For every fiber ${\rm pr}^*H_{q, \bm z}$ we have the $qKZ$-operators $K_1,\dots, K_n$ defined by equation \eqref{qkzop}.
\begin{defi}
Fix $q\in\mathbb C^*$, $\bm z\in\mathbb C^n$. The 
{\it $qKZ$-discrete connection}
 on the bundle $H$ is given by the datum of the isomorphisms of vector spaces
\beq\label{qkziso}
\Theta_{\bm z,i}^-\circ K_i(q,\bm z)\colon {\rm pr}^*H_{q,\bm z}\to {\rm pr}^*H_{q,\bm z_i^-}.
\eneq
\end{defi}

The system of difference equations 
\beq 
\label{qkz}
I(q,z_1,\dots,z_i-1,\dots,z_n)=\left[\Theta_{\bm z,i}^- \circ K_i(q,\bm z)\right]I(q,\bm z),\quad i=1,\dots, n,
\eneq
is called the system of the \emph{$qKZ$ difference equations}. These are equations for flat sections for the $qKZ$ discrete connection.

\begin{theorem}[{\cite[Theorem 3.1]{tarvar}}]
\label{compat}
The joint system of equation \eqref{eqde} and \eqref{qkz} is compatible.
\end{theorem}

\begin{oss}\label{equivqkz}
The $qKZ$ difference equations \eqref{qkz}, can be written in the equivalent form
\beq
I(q,z_1,\dots, z_i+1,\dots, z_n)=\left[\Theta_{\bm z,i}^+\circ K'_i(q,\bm z) \right]I(q,\bm z),\quad i=1,\dots, n,
\eneq
where the operators \beq
K'_i(q,\bm z):=\ ^+\Theta_{\bm z,i}\circ K_i(q,\bm z_i^+)^{-1}\circ \Theta_{\bm z,i}^+	
\eneq act on the fiber ${\rm pr}^*H_{q,\bm z}$. In terms of 
the $R$-matrices we have  
\begin{align}
\nonumber
K_i(q,\bm z_i^+)^{-1}=&R_{i+1,i}(z_{i+1}-z_i-1)\dots R_{n,i}(z_n-z_i-1)q^{E_i}\cdot\\
&\cdot R_{1,i}(z_1-z_i)\dots R_{i-1,i}(z_{i-1}-z_i).
\end{align}
\end{oss}

\begin{oss}
The $qKZ$ operators are defined in the $g$-basis \eqref{gbas}. That
 basis is the limit of the stable 
envelope basis of the equivariant cohomology of the cotangent bundle
$T^*\mathbb P^{n-1}$, in the limit, in which
 the equivariant cohomology of the cotangent bundle $T^*\mathbb P^{n-1}$
turns into the equivariant cohomology of the base space $\mathbb P^{n-1}$. 
See \cite{RTV} on the stable envelopes for the 
cotangent bundle $T^*\mathbb P^{n-1}$, see \cite[Section 7]{GRTV} and \cite[Section 11.4]{tar-var} 
on that limit.
\end{oss}

\section{Equivariant $qDE$ of $\mathbb P^{n-1}$ and its topological-enumerative solution}\label{sec6}
\subsection{Equivariant quantum differential equation} We consider the equivariant quantum differential equation \eqref{eqde} written wrt the standard basis $(x_\alpha)_{\alpha=0}^{n-1}$, namely,
\begin{equation}
\label{eqpn}\frac{dY}{dq}=\mathcal A(q,\bm {z})Y,\quad \mathcal A(q,\bm {z}):=\frac{1}{q}
\begin{pmatrix}
0&&\dots&0& q+(-1)^{n-1}s_n(\bm z)\\
1&0&\dots&0&(-1)^{n-2}s_{n-1}(\bm z)\\
0&1&\dots&0&(-1)^{n-3}s_{n-2}(\bm z)\\
&&\ddots&&\vdots\\
&&&1&s_1(\bm z)
\end{pmatrix}.
\end{equation}
We have
\beq	\label{eqpn2}\mathcal A(q,\bm z)=\mathcal A_0+\frac{1}{q}\mathcal A_1(\bm z),
\eneq where 
\[\mathcal A_0:=\begin{pmatrix}
0&&\dots&0& 1\\
0&0&\dots&0&0\\
0&0&\dots&0&0\\
&&\ddots&&\vdots\\
&&&0&0
\end{pmatrix},\quad \mathcal A_1(\bm z):=\begin{pmatrix}
0&&\dots&0& (-1)^{n-1}s_n(\bm z)\\
1&0&\dots&0&(-1)^{n-2}s_{n-1}(\bm z)\\
0&1&\dots&0&(-1)^{n-3}s_{n-2}(\bm z)\\
&&\ddots&&\vdots\\
&&&1&s_1(\bm z)
\end{pmatrix}.
\]
The eigenvalues of the matrix $\mathcal A_1(\bm z)$ are exactly $z_1,\dots, z_n$, as it easily follows from Vi\`ete formulae.
Notice that $\mathcal A_1(\bm z)$ denotes the matrix of equivariant multipilcation $$H^\bullet_{\mathbb T}(\mathbb P^{n-1},\mathbb C)\to H^\bullet_{\mathbb T}(\mathbb P^{n-1},\mathbb C),\quad f\mapsto x\cdot f,$$ whereas $\mathcal A_0$ represents the \emph{quantum correction terms} of the product. Moreover, we have 
\beq
\label{diag}
\mathcal D\cdot \mathcal A_1(\bm z)\cdot \mathcal D^{-1}=\mathcal Z:=\operatorname{diag}(z_1,\dots, z_n),
\eneq
since the classes $\Delta_i$'s are the idempotents of the equivariant cohomology algebra.

The differential system \eqref{eqpn} has a regular singularity at $q=0$ 
and an irregular singularity (of Poincar\'{e} rank 1) at $q=\infty$.

\subsection{Levelt Solution}$\quad$

\begin{theorem}
\label{teolev1}

There exist {unique} $n\times n$-matrix valued functions $(G_k(\bm z))_{k=1}^\infty$, meromorphic on 
$\mathbb C^n$ and regular on $\Omega$, such that the gauge transformation
\beq\label{gaugelevelt}
Y(q,\bm z)=G(q,\bm z)\tilde Y(q,\bm z),\quad G(q,\bm z)
=\mathcal D^{-1}\left(\mathbbm 1+\sum_{k=1}^\infty G_k(\bm z)q^k\right),
\eneq
transforms the differential system \eqref{eqpn} into the differential equation
\beq
\frac{d}{dq}\tilde Y=\frac{1}{q}\mathcal Z\tilde Y.
\eneq
Moreover, the formal power series 
 $G(q,\bm z)$ converges to a meromorphic function on $\mathbb C\times\mathbb
  C^n$, regular on $\mathbb C\times\Omega$.
\end{theorem}

\proof
Let us look for a formal gauge transformation
\[Y=\mathcal D^{-1}G\tilde Y,\quad G(q,\bm z)=\sum_{k=0}^\infty G_k(\bm z)q^k,
\]which puts the system \eqref{eqpn} into the simplest \emph{normal form}
\[\frac{d\tilde Y}{dq}=\frac{1}{q}\mathcal Z\tilde Y.
\]
This requirement
implies the following equation for $G$:
\[\mathcal D \mathcal A_0\mathcal D^{-1}G-\frac{dG}{dq}+\frac{1}{q}[\mathcal Z,G]=0,
\]which reduces to the following recurrence equations  for the coefficients $G_k$'s:
\begin{align}\label{rec1}\mathcal ZG_0=G_0\mathcal Z&,\\
\label{rec2}
\mathcal D \mathcal A_0\mathcal D^{-1}G_k+[\mathcal Z,G_{k+1}]-&(k+1)G_{k+1}=0.
\end{align}
Equation \eqref{rec1} is satisfied if and only if $G_0$ is diagonal.
So we may choose $G_0$ to be $\mathbbm 1$. 

For $k\geq 1$ equation \eqref{rec2}  {uniquely} determines
$G_{k+1}$ in terms  of $G_k$.
Indeed, the linear operator
\[\varphi_k\colon M_n(\mathbb C)\to M_n(\mathbb C),\quad X\mapsto [\mathcal Z,X]-(k+1)X,
\]
has eigenvalues $z_i-z_j-(k+1),\quad i,j=1,\dots, n,$ which are nonzero, since $\bm z\in \Omega$.
 Hence, $\varphi_k$ is invertible, and we deduce
\[G_{k+1}=\varphi_k^{-1}(-\mathcal D \mathcal A_0\mathcal D^{-1}G_k).
\]
The power series $G(q,\bm z)$ is convergent. This follows from the regularity of the singularity $q=0$ of \eqref{eqpn}. The proof is standard, e.g. see \cite{wasow,sibook,div1}.
\endproof

\begin{cor}
\label{solstand}
For $(G_k(\bm z))_{k=1}^\infty$ as in Theorem \ref{teolev1}, the matrix valued function 
\beq
\label{levelt}
Y_o(q,\bm z)=\mathcal D^{-1}\left(\mathbbm 1+\sum_{k=1}^\infty G_k(\bm z)q^k\right)q^\mathcal Z,
\eneq
is a solution of system \eqref{eqpn}. For each fixed $z\in\Omega$, the function $Y_o$ is a fundamental system of solutions.
\end{cor}

We call the fundamental solution $Y_o$ the \emph{Levelt fundamental solution}, following the terminology of \cite[Chapter 2]{ABRH}.

\smallskip
Fix $(q,\bm z)$ and increase the argument of $q$ by $2\pi$. The analytic continuation of the solutions of \eqref{eqpn} along this curve produces the \emph{monodromy operator} $M_0(\bm z)$ on the space of solutions.

\begin{cor}\label{monoeig}
The Levelt fundamental solution $Y_o(q,\bm z)$ is an eigenbasis for the the monodromy operator $M_0(\bm z)$. The matrix of the monodromy operator $M_0(\bm z)$ wrt the solution $Y_o(q,\bm z)$ is 
\beq
M_0(\bm z)=\exp\left(2\pi\sqrt{-1}\,\mathcal Z\right).
\eneq
\end{cor}

\proof
We have 
\[Y_o\left(e^{2\pi\sqrt{-1}}q,\bm z\right)=Y_o(q,\bm z)M_0(\bm z),\quad M_0(\bm z)
=\exp\left(2\pi\sqrt{-1}\,\mathcal Z\right).
\]
\endproof

\subsection{Topological-enumerative solution}
Recall the topological-enumerative morphism
$\mathcal S(\bm t,\kappa)$
of Section \ref{seceqqconn}, where $\bm t=(t^0,t^1, \dots,t^{n-1})$. Denote
\beq
\label{nn}
\mathcal S^o(q):=\mathcal S(0, \log q,0,\dots,0,1),
\eneq
where the last argument is  $\kappa=1$. We call $\mathcal S^o(q)$ the restriction of 
$\mathcal S(\bm t,\kappa)$ to the small equivariant quantum locus. 

Define the equivariant cohomology valued
functions $\Psi_{{\rm top},1}(q,\bm z),\dots, \Psi_{{\rm top},n}(q,\bm z)$ by
the formula:
\beq
\label{6.11}
\Psi_{{\rm top},m}(q,\bm z):=\mathcal S^o(q)x_{m-1},\quad m=1,\dots,n.
\eneq
By Theorem \ref{teosoper}, these functions are solutions of the equivariant quantum differential equation \eqref{eqde}. 
Let $Y_{\rm top}(q,\bm z)$ be the matrix of the operator 
$\mathcal S^o(q)$ wrt the basis $(x_\alpha)_{\alpha=0}^{n-1}$\,: 
\beq\label{psitop}
\Psi_{{\rm top},m}(q,\bm z)=\sum_{\alpha=0}^{n-1}[Y_{\rm top}(q,\bm z)]^\alpha_m x_\alpha,\quad m=1,\dots,n.
\eneq

The matrix $Y_{\rm top}(q,\bm z)$ is a solution of the matrix differential system 
\eqref{eqpn}. We call it the \emph{topological-enumerative solution} of \eqref{eqpn}.

\begin{theorem}\label{cor1}The topological-enumerative solution is the \emph{unique} solution of \eqref{eqpn} of the form
\beq\label{topsol1}Y_{\rm top}(q,\bm z)=\Phi(q,\bm z)q^{\mathcal A_1(\bm z)},
\eneq where 
\[\Phi(q,\bm z)=\mathbbm 1+\sum_{j=1}^\infty\Phi_j(\bm z)q^j.
\]The coefficients $\Phi_j$ are holomorphic on $\Omega$,  they are related to descendant Gromov-Witten invariants through the equation
\beq\Phi_j(\bm z)^\lambda_\alpha=\Bigg\langle x_\mu,\frac{x_\alpha}{1-\psi}\Bigg\rangle_{0,2,j}^{\mathbb P^{n-1},\mathbb T}\eta^{\mu\lambda},\quad j\in\mathbb N_{>0},\quad \alpha,\lambda=0,\dots, n-1.
\eneq
Here $\psi$ is the first Chern class of the universal cotangent line bundle $\mathcal L_2$ on the moduli space $\overline{\mathcal M}_{0,2}(\mathbb P^{n-1},j)$ at the second marking.
Furthermore, we have
\[Y_{\rm top}(q,\bm z)=Y_o(q,\bm z)\cdot \mathcal D,
\]
where $Y_o$ is the Levelt  fundamental solution of \eqref{eqpn} described in Corollary
 \ref{solstand} and $\mathcal D$ is \eqref{diag}. In particular, for each fixed $\bm z\in\Omega$, the matrix
  $Y_{\rm top}(q,\bm z)$ is a fundamental system of solutions of  \eqref{eqpn}.
\end{theorem}

\proof 
For  $a\in H^\bullet_{\mathbb T}(\mathbb P^{n-1},\mathbb C)$ let $\mathcal S^o(q)a$ be the corresponding solution of
the $qDE$. We obtain
\[\mathcal S^o(q)a=q^{x_1} a+\sum_{d=1}^\infty q^{d}\sum_{\lambda,\mu=0}^{n-1}
\Bigg\langle x_\mu,\frac{q^{x_1} a}{1-\psi}\Bigg\rangle_{0,2,d}^{\mathbb P^{n-1},
\mathbb T}\eta^{\mu\lambda}x_\lambda.
\]
by using the divisor axiom for descendant Gromov-Witten invariants, see \cite[Chapter 10]{cox}. Notice that
\begin{align*}
Y_{\rm top}(q,\bm z)=&\left(\mathbbm 1+\sum_{j=1}^\infty\Phi_j(\bm z)q^j\right)q^{\mathcal A_1(\bm z)}\\
=&\left(\mathbbm 1+\sum_{j=1}^\infty\Phi_j(\bm z)q^j\right)\cdot \mathcal D^{-1}\cdot\mathcal D\cdot q^{\mathcal A_1(\bm z)}\cdot\mathcal D^{-1}\cdot\mathcal D\\
=&\underbrace{\left(\mathcal D^{-1}+\sum_{j=1}^\infty\Phi_j(\bm z)\mathcal D^{-1}q^j\right)q^{\mathcal D\cdot \mathcal A_1(\bm z)\cdot \mathcal D^{-1}}}_{Y_o(q,\bm z)}\cdot \mathcal D,
\end{align*}
where in the last line we used \eqref{diag}. The uniqueness of  a solution of the form \eqref{topsol1} thus follows from the uniqueness of the solution $Y_o$ in Corollary \ref{solstand}.
\endproof

\subsection{Scalar equivariant quantum differential equation}\label{secscaleqpn}Let $Y$ be a fundamental solution of the differential system \eqref{eqpn}. Then, the matrix $\hat{Y}:=\eta\cdot Y\cdot \eta^{-1}$ is a solution of the differential system
\begin{equation}\label{eqpndual}\frac{d\hat{Y}}{dq}=\mathcal A(q,\bm z)^T\hat{Y}.
\end{equation}
This follows from the Frobenius algebra property \eqref{qFrobalg}.

Equation \eqref{eqpndual} can be reduced to the scalar differential equation
\begin{equation}\label{scaleqpn}\vartheta_q^n\phi=\left(q+(-1)^{n-1}s_n(\bm z)\right)\phi+\sum_{j=1}^{n-1}(-1)^{n-j-1}s_{n-j}(\bm z)\vartheta_q^j\phi,\quad\vartheta_q:=q\frac{d}{dq},
\end{equation}that will be called the \emph{scalar equivariant quantum differential equation} of $\mathbb P^{n-1}$.

Given $n$ linearly independent solutions $(\phi_1,\dots,\phi_n)$ of \eqref{scaleqpn} one can reconstruct a fundamental matrix solution $\hat{Y}$  of  system \eqref{eqpndual} by setting
\[(\hat{Y})^h_k:=\vartheta^h_q\phi_k,\quad h=0,\dots, n-1,\quad k=1,\dots, n.
\]
\begin{oss}
In the non-equivariant limit $z_1=\dots=z_n=0$,  equation \eqref{scaleqpn} reduces to the equation
\[\vartheta_q^n\phi=q\phi,
 \]
 which coincides with the scalar quantum differential equation of $\mathbb P^{n-1}$,
\[
\vartheta_s^n\phi=(ns)^n\phi,\quad\vartheta_s:=s\frac{d}{ds}\,,
\]
under the change of variables $q=s^n$. The monodromy of this equation has been studied in \cite{guzzetti1,CDG1}.
\end{oss}

\begin{theorem}
The matrix
\beq
\hat{Y}(q,\bm z):=\begin{pmatrix}
a_1&a_2&\dots&a_{n}\\
\vartheta_q a_1&\vartheta_q a_2&\dots&\vartheta_q a_{n}\\
\vdots&&&\vdots\\
\vartheta_q^{n-1} a_1&\vartheta_q^{n-1} a_2&\dots&\vartheta_q^{n-1} a_{n}
\end{pmatrix}(\mathcal D^{-1})^T,
\eneq
where 
\beq
a_j(q,\bm z):=q^{z_j}\left(1+\sum_{d=1}^\infty q^d\frac{1}{\prod_{i=1}^n\prod_{m=1}^d(z_j-z_i+m)}\right),\quad j=1,\dots,n,
\eneq
is a fundamental matrix solution  of the differential system \eqref{eqpndual}. The corresponding solution 
$\eta^{-1}\cdot \hat{Y}\cdot \eta$ of the equivariant differential system \eqref{eqpn} 
is the topological-enumerative solution, 
\[
Y_{\rm top}(q,\bm z)=\eta^{-1}\cdot \hat{Y}\cdot \eta.
\]
\end{theorem}

\proof
Equation \eqref{js} implies that the components, wrt to the standard basis 
$(x_\alpha)_{\alpha=1}^n$, of the $J$-function, restricted to the small equivariant quantum locus, are solutions of the scalar equivariant quantum differential equation \eqref{scaleqpn}.
The small equivariant $J$-function of $\mathbb P^{n-1}$, computed by A.\,Givental \cite{giv1}, B.J.\,Lian, K.\,Liu, S.-T.\,Yau \cite{LLY} is given by the formula
\begin{align*}
J(q,\bm z)&=q^x\left(1+\sum_{d=1}^\infty q^d\frac{1}{\prod_{i=1}^n\prod_{m=1}^d(x-z_i+m)}\right).
\end{align*}
We have 
\[x=\sum_{j=1}^nz_j\Delta_j,\quad 1=\sum_{j=1}^n\Delta_j,\quad \frac{1}{\sum_i\alpha_i\Delta_i}=\sum_i\frac{1}{\alpha_i}\Delta_i,\]
for any $\alpha_i\in\mathbb C^*$. We deduce 
\begin{align*}
J(q,\bm z)&=q^x\left(1+\sum_{d=1}^\infty q^d\frac{1}{\prod_{i=1}^n\prod_{m=1}^d\sum_{j=1}^n(z_j-z_i+m)\Delta_j}\right)\\
&=q^x\left(1+\sum_{d=1}^\infty q^d\frac{1}{\sum_{j=1}^n\prod_{i=1}^n\prod_{m=1}^d(z_j-z_i+m)\Delta_j}\right)\\
&=q^x\left(1+\sum_{d=1}^\infty\sum_{j=1}^n q^d\frac{\Delta_j}{\prod_{i=1}^n\prod_{m=1}^d(z_j-z_i+m)}\right)\\
&=q^x\sum_{j=1}^n\left(1+\sum_{d=1}^\infty q^d\frac{1}{\prod_{i=1}^n\prod_{m=1}^d(z_j-z_i+m)}\right)\Delta_j\\
&=\sum_{j=1}^nq^{z_j}\left(1+\sum_{d=1}^\infty q^d\frac{1}{\prod_{i=1}^n\prod_{m=1}^d(z_j-z_i+m)}\right)\Delta_j\\
&=\sum_{\alpha=0}^{n-1}\left\{\sum_{j=1}^nq^{z_j}\left(1+\sum_{d=1}^\infty q^d\frac{1}{\prod_{i=1}^n\prod_{m=1}^d(z_j-z_i+m)}\right)(\mathcal D^{-1})_{\alpha j}\right\}x_\alpha.
\end{align*}
If we define
\[a_j(q,\bm z):=q^{z_j}\left(1+\sum_{d=1}^\infty q^d\frac{1}{\prod_{i=1}^n\prod_{m=1}^d(z_j-z_i+m)}\right),\quad j=1,\dots,n,
\]then the matrix
\[\hat{Y}(q,\bm z):=\begin{pmatrix}
a_1&a_2&\dots&a_{n}\\
\vartheta_q a_1&\vartheta_q a_2&\dots&\vartheta_q a_{n}\\
\vdots&&&\vdots\\
\vartheta_q^{n-1} a_1&\vartheta_q^{n-1} a_2&\dots&\vartheta_q^{n-1} a_{n}
\end{pmatrix}( \mathcal D^{-1})^T
\]is a solution of the differential system \eqref{eqpndual}, and the corresponding solution $$\eta^{-1}\cdot \hat{Y}\cdot \eta$$ is the topological solution of system \eqref{eqpn}. 
\endproof

\section{Solutions of the equivariant $qDE$ and $qKZ$ difference equations}\label{sec7}
\subsection{$q$-Hypergeometric Solutions} In this section we define a fundamental system of solutions of the joint system of equations \eqref{eqde} and \eqref{qkz} described in \cite{tarvar}. 

\begin{defi}[Master and Weight function]
Define the \emph{master function} $\Phi$ and the $H^\bullet_{\mathbb T}(\mathbb P^{n-1},\mathbb C)$-valued \emph{weight function} $W$ by the formulae
\beq\label{maswei}
\Phi(t,q,\bm z):=e^{\pi\sqrt{-1}\sum_{i=1}^nz_i}\left(e^{-\pi\sqrt{-1}n}q\right)^t\prod_{a=1}^n\Gamma(z_a-t),\quad W(t,\bm y):=\prod_{j=1}^{n-1}(y_j-t).
\eneq
Recall that $y_1,\dots, y_{n-1}$ denote the equivariant Chern roots of the natural quotient bundle $\mathcal Q$ on $\mathbb P^{n-1}$.
\end{defi}

\begin{oss}\label{diff1}
Notice the difference in the definition (\ref{maswei}) of the master 
function $\Phi$ with respect to \cite{tar-var} and \cite{tarvar}. In  \cite[Section 4.1]{tarvar} the master function is defined as
\[\Phi(t,q,\bm z):=\left(e^{\pi\sqrt{-1}(2-n)}q\right)^t\prod_{a=1}^n\Gamma(z_a-t),
\]
i.e. differing from \eqref{maswei} by the
 extra factor $\exp(\pi\sqrt{-1}(2t+\sum_{i=1}^nz_i))$. 
 In \cite{tar-var} the general case of partial flag varieties 
 is considered. The master function in 
 \cite[Section 11.4]{tar-var} (see formula (11.16)), specialized to the case of projective spaces, is 
\[
\Phi(t,q,\bm z):=e^{-\pi\sqrt{-1}\sum_{i=1}^nz_i}\left(e^{-\pi\sqrt{-1}n}q\right)^t\prod_{a=1}^n\Gamma(z_a-t).
\]
Thus it differs from the function $\Phi$ in \eqref{maswei} by the factor $\exp(2\pi\sqrt{-1}\sum_{i=1}^nz_i)$.
\end{oss}

\begin{defi}[Jackson Integrals] Define the Jackson integrals $\Psi_J$, 
 $J=1,\dots, n$, to be the $H^\bullet_{\mathbb T}(\mathbb P^{n-1},\mathbb C)$-valued 
 functions defined on $\widetilde{\mathbb C^*}\times\Omega$ by the formula
\beq\label{jackint}
\Psi_J(q,\bm y,\bm z):=-\sum_{r=0}^\infty\underset{t=z_J+r}{\rm Res}\Phi(t,q,\bm z)W(t,\bm y).
\eneq
Here $\widetilde{\mathbb C^*}$ is the universal cover of $\mathbb C^*$.
\end{defi}

\begin{theorem}[\cite{tar-var}]
The functions $\Psi_J(q,\bm y,\bm z)$ with $J=1,\dots, n$ are holomorphic on $\widetilde{\mathbb C^*}\times\Omega$. Each of them is a solution of the equivariant quantum differential equation \eqref{eqde} and of 
the qKZ difference equations \eqref{qkz}. These functions form a basis of solutions of this joint system of equations. 
\end{theorem}

We will call the solutions $\Psi_J$ the \emph{$q$-hypergeometric solutions}.

\begin{oss}
Notice that in \cite{tarvar} the system of $qKZ$ equations differs from the one considered in \cite{tar-var} and in this article by a sign in the rhs of \eqref{qkz}, due to the
 different normalization of the master function $\Phi$.
\end{oss}

\begin{cor}[{\cite[Formula (11.19)]{tar-var}}]\label{cor2}
The $q$-hypergeometric solutions $\Psi_J$ admit the following expansion
\[\Psi_J(q,\bm y,\bm z)=e^{\pi\sqrt{-1}\sum_{i=1}^nz_i}\left(e^{-\pi\sqrt{-1}n}q\right)^{z_J}\prod_{a\neq J}\Gamma(1+z_a-z_J)\left(\Delta_J+\sum_{k=1}^\infty\Psi_{J,k}(\bm z)q^k\right),
\]where the classes $\Psi_{J,k}(\bm z)$ are rational functions in $(z_1,\dots,z_n)$, regular on $\Omega$.
\end{cor}

Define the matrix $Y_{\text{$q$-hyp}}=\left([Y_{\text{$q$-hyp}}]^\alpha_J\right)_{\alpha,J}$ by the formula
\[
\Psi_J=\sum_{\alpha=0}^{n-1}
[Y_{\text{$q$-hyp}}]^\alpha_J\ x_\alpha,
\quad J=1,\dots, n,
\]
then $Y_{\text{$q$-hyp}}$ is  a fundamental matrix solution of the differential system \eqref{eqpn}.

\begin{theorem}\label{teoconn1}
The connection matrix $C$ relating the topological-enumerative solution with the $q$-hypergeometric solution, 
\[
Y_{\textnormal{$q$-hyp}}(q,\bm z)=Y_{\rm top}(q, \bm z)\cdot C,
\]
is given by the formula
\beq
\label{conn1}
C=\mathcal D^{-1}\cdot{\rm diag}\left(e^{\pi\sqrt{-1}(-nz_j+\sum_{i=1}^nz_i)}
\prod_{a\neq j}\Gamma(1+z_a-z_j)\right)_{j=1}^n.
\eneq
This matric $C$ is the matrix attached to the morphism
\beq \rho\colon H^\bullet_{\mathbb T}(\mathbb P^{n-1},\mathbb C)\to H^{\Omega}_{\mathbb T}(\mathbb P^{n-1}),\quad v\mapsto  e^{\pi\sqrt{-1}c_1(\mathbb P^{n-1})}
\cdot\widehat\Gamma^+_{\mathbb P^{n-1}}\cdot v,
\eneq where we fix 
\begin{itemize}
\item the basis $(\Delta_j)_{j=1}^n$ in the domain of $\rho$,
\item the basis $(x_\alpha)_{\alpha=0}^{n-1}$ in the target space of $\rho$.
\end{itemize} 
\end{theorem}

\proof The proof follows from Theorem \ref{cor1} and Corollary \ref{cor2}.
Notice that
\begin{align*}
c_1(\mathbb P^{n-1})=\sum_{i=1}^nz_i-nx,
\qquad
\Gamma^+_{\mathbb P^{n-1}}=\prod_{a=1}^{n-1}\Gamma(1+y_a-x).
\end{align*}
Each term of the entries of the diagonal matrix in \eqref{conn1} can be indentified
 with the multiplication by these classes wrt the basis $(\Delta_i)_{i}$.
\endproof

\begin{oss}

The functions $\Psi_{{\rm top},m}(q,\bm z)$ defined in \eqref{6.11} are not solutions of the $qKZ$ difference equations \eqref{qkz}, since the matrix $C$ given by \eqref{conn1} is not 1-periodic in the equivariant parameters $z_1,\dots,z_n$.
\end{oss}

\subsection{Identification of solutions with $K$-theoretical classes}\label{solKclas}
Following \cite{tarvar}, we introduce the symbols
\beq
\acute{T}:=\exp(2\pi\sqrt{-1}t),\quad \acute{Z_J}:=\exp(2\pi\sqrt{-1}z_J),\quad J=1,\dots, n.
\eneq
\begin{defi}
Let $Q(X,\bm{Z})\in\mathbb C[X^{\pm1}, \bm{Z}^{\pm 1}]$ be a Laurent polynomial. Define 
\[\Psi_Q(q,\bm y,\bm z):=\sum_{J=1}^nQ(\acute{Z_J},\acute{\bm Z})\Psi_J(q,\bm y,\bm z).
\]The function $\Psi_Q$ is a solution of the joint system of equations \eqref{eqde} and \eqref{qkz}. If $Q(X,\bm{Z})=X^{m}$, we denote the corresponding solution $\Psi_Q$ by $\Psi^m$, i.e.
\[\Psi^m=\sum_{J=1}^n\acute{Z}_J^m\Psi_J.
\]
\end{defi}

\begin{oss}\label{diff2}
Notice that $\Psi^m$ in \cite{tarvar} equals $\exp\left(-\pi\sqrt{-1}\sum_{i=1}^nz_i\right)\Psi^{m-2}$ of this paper. This is due to the difference of normalizations of the master function, see Remark \ref{diff1}. 

\end{oss}

\begin{theorem}[\cite{tarvar}]\label{teosolkclass}
There is a well-defined morphism from $K_0^{\mathbb T}(\mathbb P^{n-1})_{\mathbb C}$ to the space of solutions of the joint system of equations \eqref{eqde} and \eqref{qkz}, defined by the association 
\[Q\mapsto\Psi_Q,
\]under the isomorphism \eqref{Kpn}.
\end{theorem}

\begin{cor}[{\cite[Corollary 4.4]{tarvar}}]\label{corsol}
For any $k\in\mathbb Z$, we have
\beq
\sum_{i=0}^n(-1)^{n-i}s_{n-i}(\acute{\bm  Z})\Psi^{k+i}(q,\bm y, \bm z)=0,
\eneq
where $s_i(\acute{\bm Z})$ are the elementary symmetric polynomials in $\acute{\bm Z}$.
\end{cor}

\begin{theorem}[{\cite[Theorem 11.3]{tar-var}}]\label{solkth}
For any $k$, $(\Psi^{k+i}(q,\bm y,\bm z))_{i=0}^{n-1}$ is a basis of the space of solutions 
of the joint system \eqref{eqde} and \eqref{qkz}.
\end{theorem}

\begin{oss}
The idea that space of solutions of the 
$qDE$ and $qKZ$  equations
is naturally identified with the space of the $K$-algebra can be observed in \cite{tarvar971,tarvar972} and was implicitly discussed  there.
\end{oss}

\subsection{Module $\mathcal S_n$ of solutions}

\begin{defi}
Define the space $\mathcal S_n$ of solutions to the joint system \eqref{eqde} and \eqref{qkz} of the form
\beq
\sum_{m=1}^nQ_m(\acute{\bm Z})\Psi^m(q,\bm y,\bm z),\quad Q_m\in\mathbb C[\bm Z^{\pm 1}].
\eneq 
The space $\mathcal S_n$ admits a structure of a $\mathbb C[\bm Z^{\pm 1}]$-module, the multiplication by $Q(\bm Z)$ being defines as the multiplication by $Q(\acute{\bm Z})$.
\end{defi}

By Corollary \ref{corsol}, the module $\mathcal S_n$ contains all the solutions $\Psi^m(q,\bm y,\bm z)$, $m\in\mathbb Z$.

\begin{cor}[{\cite[Corollary 4.6]{tarvar}}]\label{corsn}
The module $\mathcal S_n$ contains a basis of solutions of the joint system of equations \eqref{eqde},
\eqref{qkz}. Moreover, the map 
$\theta\colon K_0^{\mathbb T}(\mathbb P^{n-1})_{\mathbb C}\to\mathcal S_n$ defined by
the formula
\beq
\label{isoth}
\theta(X^m):=\Psi^m(q,\bm y, \bm z),\quad m\in\mathbb Z,
\eneq 
defines an isomorphism of $\mathbb C[\bm Z^{\pm 1}]$-modules.
\end{cor}

Using the isomorphism $\theta$ we define a sesquilinear form on the module $\mathcal S_n$ as the image of the $\chi^\mathbb T$-form on $K_0^\mathbb T(\mathbb P^{n-1})_{\mathbb C}$. The notions of exceptional bases and the action of the braid group on them can be lifted to $\mathcal S_n$.

\subsection{Integral representations for solutions}For $p\in\mathbb C$, let us denote by $C(p)$ the parabola in $\mathbb C$ defined by the equation
\[C(p):=\left\{p+t^2+t\sqrt{-1}\colon t\in\mathbb R\right\}.
\]Given a point $\bm z\in\Omega$, take $p$ such that all the points $z_1,\dots,z_n$ line inside $C(p)$. The value of the integral \eqref{int} below does not depend on a particular choice of $p$, so we will simply denote $C(p)$ by $C(\bm z)$.

\begin{lemma}[{\cite[Lemma 11.5]{tar-var}}]For any Laurent polynomial $Q(X,\bm Z)$ we have
\beq\label{int}
\Psi_Q(q,\bm y,\bm z)=\frac{1}{2\pi\sqrt{-1}}\int_{C(\bm z)}Q(\acute{T},\acute{\bm Z})\Phi(t,q,\bm z)W(t,\bm y)\,dt,
\eneq
where the integral converges for any $(q,\bm z)\in \widetilde{\mathbb C^*}\times\Omega$. In particular, we have
\beq\label{int2}
\Psi^m(q,\bm y,\bm z)=\frac{e^{\pi\sqrt{-1}\sum_{i=1}^nz_i}}{2\pi\sqrt{-1}}\int_{C(\bm z)}e^{2\pi\sqrt{-1}mt}e^{-\pi\sqrt{-1}nt}q^t\prod_{a=1}^n\Gamma(z_a-t)\prod_{j=1}^{n-1}(y_j-t)\,dt.
\eneq
\end{lemma}

\begin{oss}
These formulae differ from the corresponding ones in \cite{tarvar}. See also Remarks \ref{diff1} and \ref{diff2}. 
\end{oss}

\subsection{Coxeter element, and elements $\gamma_n,\delta_{n,{\rm odd}},\delta_{n,\rm even}\in\mathcal B_n$}\label{coxetersec} The Coxeter element of $\mathcal B_n$ is the braid
\beq
C:=\tau_1\tau_2\dots\tau_{n-1}\in\mathcal B_n.
\eneq
For any $n\geq 3$, let  
\begin{empheq}[left={\ell_n:=}\empheqlbrace]{align*} 
      n-1,&\quad\text{for $n$ odd},\\
      \\
      n-2,&\quad\text{for $n$ even}.
        \end{empheq}
Set $\gamma_2:=1$, and for $n\geq 3$,
\[\beta_k:=\tau_k\tau_{k+1}\dots\tau_{n-1},\quad \gamma_n:=\beta_{\ell_n}\beta_{\ell_{n}-2}\dots\beta_2.
\]
Define also 
\[\delta_{n,\rm odd}=\tau_1\tau_3\dots\tau_{n-2},\quad \delta_{n,\rm even}=\tau_2\tau_4\dots\tau_{n-1},\quad \text{for $n$ odd},
\]
\[\delta_{n,\rm odd}=\tau_1\tau_3\dots\tau_{n-1},\quad \delta_{n,\rm even}=\tau_2\tau_4\dots\tau_{n-2},\quad \text{for $n$ even}.
\]
The elements $C,\ \gamma_n,\ \delta_{n,{\rm odd}},\ \delta_{n,\rm even}$ satisfy the following relation.
\begin{lemma}[{\cite[Lemma 6.3]{tarvar}}]\label{lemidbr1}
We have the following identity in $\mathcal B_n$:
\beq\label{idbr1}
\delta_{n,\rm even}\ \delta_{n,\rm odd}\ \gamma_n=\gamma_n\ C.
\eneq 
\end{lemma}
\subsection{Exceptional bases $Q_k,Q_k',Q_k'',\widetilde{Q}_k,\widetilde{Q}_k',\widetilde{Q}_k''$}\label{qbas} For any $k\in\mathbb Z$, we define the basis $Q_k$ of solutions of the joint system \eqref{eqde}, \eqref{qkz} to be the basis
\beq Q_k:=(\Psi^{k+n-1},\dots,\Psi^{k+1},\Psi^k).
\eneq 

\begin{lemma}\label{qexcbas}
The basis $Q_k$ is an exceptional basis of $\mathcal S_n$. Via the isomorphism
$\quad$ $\theta\colon K_0^{\mathbb T}(\mathbb P^{n-1})_{\mathbb C}\to\mathcal S_n$, 
it is identified with the  exceptional basis 
\[([\mathcal O(-k-n+1)],\dots,[\mathcal O(-k-1)],[\mathcal O(-k)]),
\]
of $K_0^{\mathbb T}(\mathbb P^{n-1})_{\mathbb C}$,
obtained from the Beilinson basis $([\mathcal O(i)])_{i=0}^{n-1}$ 
by twisting it with $\otimes [\mathcal O(-k-n+1)]$.
\end{lemma}
\proof
It follows from Corollary \ref{corsn}.
\endproof
For any $k\in\mathbb Z$, we define the exceptional bases $Q_k'$ and $Q_k''$ through the mutations
\beq Q_k':=\gamma_n Q_k,\quad Q_k'':=\delta_{n,\rm odd}Q'_k.
\eneq

\begin{prop}[{\cite[Lemma 6.6, Corollary 7.2]{tarvar}}]
The basis $Q_{k}$ and $Q_{k-1}$ are related by the so-called \emph{modified Coxeter map}: this means that $Q_{k-1}$ is obtained from $CQ_k$  by multiplying its last element by $(-1)^{n+1}s_n(\bm Z^{-1})$. Moreover, the basis $Q_{k-1}'$ is obtained from the basis $\delta_{n,\rm even}Q''_k$ by multiplying its last vector by $(-1)^{n+1}s_n(\bm Z^{-1})$.
\end{prop}

\begin{oss}\label{diff3}Our bases $Q_k,Q_k',Q_k''$ have the same elements as
 the bases $Q_{k-1},Q_{k-1}',Q_{k-1}''$ of \cite{tarvar}, but ordered in the opposite way, see Remark \ref{chi-A}.
\end{oss}

Introduce three more families of exceptional bases of $\mathcal S_n$, denoted by $\widetilde{Q}_k,\widetilde{Q}_k',\widetilde{Q}_k''$.  For $k\in \mathbb Z$ define
\beq
 \widetilde Q_k:=C^{-k}\widetilde Q_0\,,
\qquad
\widetilde Q_k':=\gamma_n \widetilde Q_k,\qquad \widetilde Q_k'':=\delta_{n,\rm odd}\widetilde Q'_k.
\eneq
We have
\[\widetilde Q_{k-1}'=\delta_{n,\rm even}\widetilde Q_k''
\]
by formula  \eqref{idbr1}. The diagram
\beq\label{diagQ}
\xymatrix{
\widetilde Q_k'\ar@{|->}[rr]^{\delta_{n,\rm odd}}&&\widetilde Q_k''\ar@{|->}[rr]^{\delta_{n,\rm even}}&&\widetilde Q_{k-1}'\\
\widetilde Q_k\ar@{|->}[rrrr]_C\ar@{|->}[u]^{\gamma_n}&&&&\widetilde Q_{k-1}\ar@{|->}[u]_{\gamma_n}
}
\eneq
is commutative by Lemma \ref{lemidbr1}.

\begin{oss}
Notice that the pre-image of a basis $\widetilde{Q}_k$ via 
the isomorphism $\theta\colon K_0^{\mathbb T}(\mathbb P^{n-1})_{\mathbb C}\to\mathcal S_n$ 
is a  foundation $\frak E_k$ of the helix generated by the Beilinson exceptional collection $(\mathcal O(-n+1),\dots,\mathcal O(-1),\mathcal O)$. In particular, $\frak E_n$ is the adjacent foundation to the left of $\frak E_0$: by Proposition \ref{serreexc}, the objects of $\frak E_n$ are obtained by applying the Serre functor to objects of $\frak E_0$.
\end{oss}

For any $\ell,m\in\mathbb Z$ such that $0\leq m-\ell \leq n$, denote
\beq\label{psiml}\Psi^m(\ell):=\Psi^m-s_1(\bm Z)\Psi^{m-1}+\dots+(-1)^{m-\ell}s_{m-\ell}(\bm Z)\Psi^\ell.
\eneq 

The explicit formulae for $Q_k'$ and $Q''_k$ follow from
  \cite[Section 6.3]{tarvar}, Remarks \ref{diff2} and \ref{diff3}.

If $n=2h+1$, we have
\begin{itemize}
\item the basis $Q_k'$ is the basis in which the solutions $\Psi^k,\dots,\Psi^{k+h}$ stay at the positions $2h+1,2h-1,\dots,1$, and the solutions $\Psi^{k+2h}(k+1),\Psi^{k+2h-1}(k+2),\dots,$ $\Psi^{k+h+1}(k+h)$ stay at the positions $2h,2h-2,\dots,6,4,2$, respectively;
\item the basis $Q_k''$ is the basis in which the solutions $\Psi^k,\dots,\Psi^{k+h}$ stay at the positions $2h,2h-2,\dots,6,4,2,1$, and the solutions $\Psi^{k+2h}(k),\Psi^{k+2h-1}(k+1),\dots,$ $\Psi^{k+h+1}(k+h-1)$ stay at the positions $2h+1,\dots,5,3$.
\end{itemize}

If $n=2h$, we have
\begin{itemize}
\item the basis $Q_k'$ is the basis in which the solutions $\Psi^k,\dots,\Psi^{k+h}$ stay at the positions $2h,2h-2,\dots,4,2,1$, and the solutions $\Psi^{k+2h-1}(k+1),\Psi^{k+2h-2}(k+2),\dots,$ $\Psi^{k+h+1}(k+h-1)$ stay at the positions $2h-1,\dots,7,5,3$;
\item the basis $Q_k''$ is the basis in which the solutions $\Psi^k,\dots,\Psi^{k+h-1}$ stay at the positions $2h-1,\dots,5,3,1$, and the solutions $\Psi^{k+2h-1}(k),\Psi^{k+2h-2}(k+1),\dots,\quad$
 $\Psi^{k+h}(k+h-1)$ stay at the positions $2h,\dots,6,4,2$.
\end{itemize}

The bases $\widetilde{Q}_k',\widetilde{Q}_k''$ can be obtained from the bases $Q_k',Q_k''$ by application of the following rule.

\begin{regola}
\label{regola}

The basis $\widetilde{Q}_k'$ (resp., $\widetilde{Q}_k''$) is obtained from
the basis $Q_k'$ (resp., $Q_k''$) by substituting any solution $\Psi^m(\ell)$ with
\[\left((-1)^{n+1}s_n(\bm Z)\right)^a\Psi^m(\ell),
\]where $a\in\mathbb Z$ is such that
\[0\leq m+an\leq n-1.
\]
\end{regola}

\begin{es}Let $n=5$ and $k=-1$. We have 
\begin{align*}
Q_{-1}'=&\ (\Psi^1,\Psi^2(1),\Psi^0,\Psi^3(0),\Psi^{-1}),\\
\widetilde{Q}_{-1}'=&\ (\Psi^1,\Psi^2(1),\Psi^0,\Psi^3(0),s_5(\bm Z)\Psi^{-1}),\\
Q_{-1}''=&\ (\Psi^1,\Psi^0,\Psi^2(0),\Psi^{-1},\Psi^3(-1)),\\
\widetilde{Q}_{-1}''=&\ (\Psi^1,\Psi^0,\Psi^2(0),s_5(\bm Z)\Psi^{-1},\Psi^3(-1)).
\end{align*}
\end{es}

\begin{theorem}\label{exobj}
Via the isomorphism $\theta\colon K_0^{\mathbb T}(\mathbb P^{n-1})_{\mathbb C}\to\mathcal S_n$, the solution $\Psi^m(\ell)$ corresponds 
to the $K$-class of the exceptional object \[
\bigwedge\nolimits^{ m-\ell}\mathcal T(-m):=\left(\bigwedge\nolimits^{ m-\ell}\mathcal T\right)\otimes \mathcal O(-m),\]
placed in degree $\ell-m$. Here $\mathcal T$ denotes the tangent sheaf of $\mathbb P^{n-1}$ 
 with its natural $\mathbb T$-equivariant structure. 
\end{theorem}

\proof Let $V=\mathbb C^n$ be the diagonal representation on $\mathbb T$ described in Section \ref{torusaction}.
Consider the Euler exact sequence, together with its exterior powers 
\begin{equation}\label{euler}\xymatrix@R=1pt{
0\ar[r]&\mathcal O\ar[r]&V\otimes \mathcal O(1)\ar[r]&\mathcal T\ar[r]&0,\\
0\ar[r]&\mathcal T\ar[r]&\bigwedge\nolimits^2V\otimes \mathcal O(2)\ar[r]&\bigwedge\nolimits^2\mathcal T\ar[r]&0,\\
&\vdots&\vdots&\vdots\\
0\ar[r]&\bigwedge\nolimits^{h-1}\mathcal T\ar[r]&\bigwedge\nolimits^hV\otimes \mathcal O(h)\ar[r]&\bigwedge\nolimits^h\mathcal T\ar[r]&0,\\
&\vdots&\vdots&\vdots\\
0\ar[r]&\bigwedge\nolimits^{n-2}\mathcal T\ar[r]&\bigwedge\nolimits^{n-1}V\otimes \mathcal O(n-1)\ar[r]&\mathcal O(n)\ar[r]&0.\\
}
\end{equation}Each morphism in \eqref{euler} is $\mathbb T$-equivariant. In equivariant $K$-theory we have
\[
\left[\bigwedge\nolimits^h\mathcal T\right]=\left[\bigwedge\nolimits^hV\otimes \mathcal O(h)\right]-\left[\bigwedge\nolimits^{h-1}\mathcal T\right]=s_h(\bm Z)\left[\mathcal O(h)\right]-\left[\bigwedge\nolimits^{h-1}\mathcal T\right],
\]
for $h=1,\dots, n$. By induction, we obtain
\[
\left[\bigwedge\nolimits^h\mathcal T\right]=\pm\sum_{j=0}^h(-1)^js_j(\bm Z)\left[\mathcal O(j)\right],
\]
where the sign is $+$ for even $h$, and $-$ for odd $h$.
The result follows from identity \eqref{psiml}.
\endproof

\begin{cor}
\label{corexcol}
For any $k\in\mathbb Z$, via the isomorphism $\theta\colon K_0^{\mathbb T}(\mathbb P^{n-1})_{\mathbb C}\to\mathcal S_n$, the bases $Q_k',Q_k''$ correspond to the $K$-classes of the following $\mathbb T$-full exceptional collections: \begin{enumerate}
\item If $n$ is \emph{odd}, the basis $Q_k'$ corresponds to
{\footnotesize
\begin{align*}
\mathcal O\left(-k-\frac{n-1}{2}\right),\ &\mathcal T\left(-k-\frac{n-1}{2}-1\right),\ \mathcal O\left(-k-\frac{n-1}{2}+1\right),\ \bigwedge\nolimits^3\mathcal T\left(-k-\frac{n-1}{2}-2\right),\dots,\\
\dots&\ ,\ \bigwedge\nolimits^{n-4}\mathcal T\left(-k-n+2\right),\ \mathcal O(-k-1),\ \bigwedge\nolimits^{n-2}\mathcal T\left(-k-n+1\right),\ \mathcal O(-k),
\end{align*}
}
and the basis $Q_k''$ corresponds to
{\footnotesize
\begin{align*}
\mathcal O\left(-k-\frac{n-1}{2}\right),\ &\mathcal O\left(-k-\frac{n-1}{2}+1\right),\ \bigwedge\nolimits^2\mathcal T\left(-k-\frac{n-1}{2}-1\right),\ \mathcal O\left(-k-\frac{n-1}{2}+2\right),\dots,\\
\dots\ &,\ \mathcal O(-k-1),\ \bigwedge\nolimits^{n-3}\mathcal T\left(-k-n+2\right),\ \mathcal O(-k),\ \bigwedge\nolimits^{n-1}\mathcal T\left(-k-n+1\right).
\end{align*}
}
\item If $n$ is \emph{even}, the basis $Q_k'$ corresponds to
{\footnotesize
\begin{align*}
\mathcal O\left(-k-\frac{n}{2}\right),\ &\mathcal O\left(-k-\frac{n}{2}+1\right),\ \bigwedge\nolimits^2\mathcal T\left(-k-\frac{n}{2}-1\right),\ \mathcal O\left(-k-\frac{n}{2}+2\right),\dots,\\
\dots\ &,\ \bigwedge\nolimits^{n-4}\mathcal T\left(-k-n+2\right),\ \mathcal O(-k-1),\ \bigwedge\nolimits^{n-2}\mathcal T\left(-k-n+1\right),\ \mathcal O(-k),
\end{align*}
}and the basis $Q_k''$ corresponds to
{\footnotesize
\begin{align*}
\mathcal O\left(-k-\frac{n}{2}+1\right),\ &\mathcal T\left(-k-\frac{n}{2}\right),\ \mathcal O\left(-k-\frac{n}{2}+2\right),\ \bigwedge\nolimits^3\mathcal T\left(-k-\frac{n}{2}-1\right),\dots,\\
\dots&\ ,\ \mathcal O(-k-1),\ \bigwedge\nolimits^{n-3}\mathcal T\left(-k-n+2\right),\ \mathcal O(-k),\ \bigwedge\nolimits^{n-1}\mathcal T\left(-k-n+1\right).
\end{align*}
}
\end{enumerate}
In these exceptional collections, each of the objects 
$\mathcal O(m)$  sits in degree $0$ and each of the objects
  $\bigwedge^h\mathcal T(m)$ sits in degree $-h$.
\end{cor}

\proof
It follows from Theorem \ref{exobj} and the description of the bases $Q_k',Q_k''$ given above.
\endproof

\begin{cor}
\label{corexcol-t}

The objects corresponding to the elements of the 
bases $\widetilde{Q}_k',\widetilde{Q}_k''$ are obtained 
from the objects corresponding to the elements of the 
bases $Q_k',Q_k''$ by twisting their $\mathbb T$-equivariant structures.

More precisely, for $a\in\mathbb Z$ define the $\mathbb T$-characters 
\beq 
\underbrace{\bigwedge\nolimits ^n V\otimes\dots\otimes \bigwedge\nolimits ^n V}_{a\text{ times}},\quad \text{if}\
\ a\geq 0,
\eneq
\beq
 \underbrace{\bigwedge\nolimits ^n V^*\otimes\dots\otimes \bigwedge\nolimits ^n V^*}_{-a\text{ times}},\quad 
 \quad \text{if}\
\ a< 0,
\eneq
where $V\cong\mathbb C^n$ is the diagonal representation of $\mathbb T$.
Given $m\in \mathbb Z$ define $a\in \mathbb Z$ from 
\beq
0\leq m+an\leq n-1.
\eneq

Then the $\mathbb T$-equivariant structure of any
object $\mathcal O(-m)$ or $\bigwedge^{m-\ell}\mathcal T(-m)$ corresponding to bases
$Q_k',Q_k''$  must be tensored with the corresponding character defined above.
\end{cor}

\proof
It follows from Corollary \ref{corexcol} and Rule \ref{regola}.
\endproof

\subsection{Asymptotic expansion of  bases $Q_k'$ and $Q_k''$ in sectors $\mathcal V_k'$ and $\mathcal V_k''$}
\label{asymexps}

 Following  \cite{tarvar},  introduce the
 coordinates $(r,\phi)$ on the universal cover $\widetilde{\mathbb C^*}$ of the punctured $s$-plane $\mathbb C^*$:
 \beq\label{eqrhophi}
q=s^n,\quad s=re^{-2\pi\sqrt{-1}\phi},\quad r>0,\ \phi\in\mathbb R.
\eneq

\begin{lemma}[{\cite[Lemma 5.1]{tarvar}}]\label{lemmapsim}
For $m\in\mathbb Z$ and $\phi\in\mathbb R$ such that
\beq
\frac{m}{n}-1<\phi<\frac{m}{n},
\eneq
we have the asymptotic expansion as $s\to\infty$
\beq
\Psi^m(s^n,\bm y,\bm z)=\frac{(2\pi)^\frac{n-1}{2}}{\sqrt{n}}e^{\pi\sqrt{-1}\sum_{i=1}^nz_i}\left(e^{-\sqrt{-1}\pi}\zeta_n^ms\right)^{\sum_{i=1}^nz_i+\frac{n-1}{2}}e^{ns\zeta_n^m}\left(1+O\left(\frac{1}{s}\right)\right),
\eneq
\[\zeta_n:=\exp\left(\frac{2\pi\sqrt{-1}}{n}\right),
\]
where  $\arg\left(e^{-\sqrt{-1}\pi}\zeta_n^ms\right)=2\pi \frac{m}{n}-\pi-2\pi\phi$, so that $\left|\arg\left(e^{-\sqrt{-1}\pi}\zeta_n^ms\right)\right|<\pi$.
\end{lemma}

Consider the sectors
\beq\label{stsec1}\mathcal V_k':=\left\{s\in\widetilde{\mathbb C^*}\colon \frac{k}{n}-\frac{1}{2}-\frac{1}{2n}<\phi<\frac{k}{n}\right\},
\eneq
\beq\label{stsec2}
\mathcal V_k'':=\left\{s\in\widetilde{\mathbb C^*}\colon \frac{k}{n}-\frac{1}{2}-\frac{1}{n}<\phi<\frac{k}{n}-\frac{1}{2n}\right\},
\eneq
for $k\in\mathbb Z$.

Let us recall the main result of \cite{tarvar} concerning the asymptotic expansion of the bases $Q_k', Q_k''$.

\begin{theorem}[{\cite[Theorem 7.1]{tarvar}}]\label{maintarvar}
The elements of the basis of solutions $Q_k'$ (resp. $Q_k''$) can be reordered to a basis $(I_m(s^n,\bm y,\bm z))_{m=0}^{n-1}$ with asymptotic expansion
\beq
I_m(s^n,\bm y,\bm z)=\frac{(2\pi)^\frac{n-1}{2}}{\sqrt{n}}e^{\pi\sqrt{-1}\sum_{i=1}^nz_i}\left(e^{-\sqrt{-1}\pi}\zeta_n^ms\right)^{\sum_{i=1}^nz_i+\frac{n-1}{2}}e^{ns\zeta_n^m}\left(1+O\left(\frac{1}{s}\right)\right),
\eneq
for $s\to\infty$ and $s\in\mathcal V_k'$ (resp $s\in\mathcal V_k''$).
Here, for defining $\left(e^{-\sqrt{-1}\pi}\zeta_n^ms\right)^{\sum_{i=1}^nz_i+\frac{n-1}{2}}$, the following choice of the branch of $\log\left(e^{-\sqrt{-1}\pi}\zeta_n^ms\right)$ is done: for every $m$ the argument of $\left(e^{-\sqrt{-1}\pi}\zeta_n^ms\right)$ is chosen so that $|\arg(e^{-\sqrt{-1}\pi}\zeta_n^ms)|<\pi$ when $\phi$ tends 
\begin{enumerate}
\item to $\frac{k}{n}$ inside $\left(\frac{2k-n-1}{2n},\frac{k}{n}\right)$, for the case of $Q_k'$,
\item to $\frac{2k-1}{2n}$ inside $\left(\frac{2k-n-2}{2n},\frac{2k-1}{2n}\right)$, for the case of $Q_k''$.
\end{enumerate}
In both cases, the argument of $\left(e^{-\sqrt{-1}\pi}\zeta_n^ms\right)$ is continuous for $\phi$ in the intervals above.
\end{theorem}

In terms of the bases $\widetilde{Q}_k'$ and $\widetilde{Q}_k''$ we can recast this result as follows.

\begin{cor}\label{corqtilde}
The elements of the basis of solutions $\widetilde{Q}_k'$ (resp. $\widetilde{Q}_k''$) can be reordered to a basis $(I_m(s^n,\bm y,\bm z))_{m=0}^{n-1}$ with asymptotic expansion
\beq
I_m(s^n,\bm y,\bm z)=\frac{(2\pi)^\frac{n-1}{2}}{\sqrt{n}}e^{-\pi\sqrt{-1}\frac{n-1}{2}}\left(\zeta_n^ms\right)^{\sum_{i=1}^nz_i+\frac{n-1}{2}}e^{ns\zeta_n^m}\left(1+O\left(\frac{1}{s}\right)\right),
\eneq
for $s\to\infty$ and $s\in\mathcal V_k'$ (resp $s\in\mathcal V_k''$).
Here for defining $\left(\zeta_n^ms\right)^{\sum_{i=1}^nz_i+\frac{n-1}{2}}$, the principal determination of the argument of $\left(\zeta_n^ms\right)$ is chosen, i.e. $\arg\left(\zeta_n^ms\right)=2\pi\frac{m}{n}-2\pi\phi$.
\end{cor}

\proof
From Lemma \ref{lemmapsim} and Theorem \ref{maintarvar}, we have that, for any $m\in\mathbb Z$, the element $\Psi^m(\ell)$ of $Q_k'$ admits the following expansion on $\mathcal V_k'$:
\begin{align*}
\Psi^m(\ell)(s^n,\bm y,\bm z)&=\frac{(2\pi)^\frac{n-1}{2}}{\sqrt{n}}e^{\pi\sqrt{-1}s_1(\bm z)}\left(e^{-\sqrt{-1}\pi}\zeta_n^ms\right)^{s_1(\bm z)+\frac{n-1}{2}}e^{ns\zeta_n^m}\left(1+O\left(\frac{1}{s}\right)\right)\\
=&\frac{(2\pi)^\frac{n-1}{2}}{\sqrt{n}}e^{\pi\sqrt{-1}s_1(\bm z)}\left(e^{-\sqrt{-1}\pi}\zeta_n^{m+an}s\right)^{s_1(\bm z)+\frac{n-1}{2}}\cdot\\
&\cdot(\zeta_n^{an})^{-s_1(\bm z)-\frac{n-1}{2}}e^{ns\zeta_n^m}\left(1+O\left(\frac{1}{s}\right)\right),
\end{align*}
where $a\in\mathbb Z$ is such that
\[0\leq m+an\leq n-1.
\]
Thus, we have the following asymptotic expansion for $s\to \infty$ and $s\in\mathcal V_k'$:
\begin{align*}(\zeta_n^{an})^{s_1(\bm z)+\frac{n-1}{2}}\Psi^m(\ell)&(s^n,\bm y,\bm z)=\\
&\frac{(2\pi)^\frac{n-1}{2}}{\sqrt{n}}e^{-\pi\sqrt{-1}\frac{n-1}{2}}\left(\zeta_n^{m+an}s\right)^{s_1(\bm z)+\frac{n-1}{2}}e^{ns\zeta_n^m}\left(1+O\left(\frac{1}{s}\right)\right).
\end{align*}Notice that
\[(\zeta_n^{an})^{s_1(\bm z)+\frac{n-1}{2}}\Psi^m(\ell)=\left((-1)^{n-1}s_n(\bm{\acute{Z}})\right)^a\Psi^m(\ell)
\]is the element of $\widetilde{Q}_k'$ obtained the element  $\Psi^m(\ell)$ of $Q_k'$ by applying the Rule \ref{regola}. This proves the statement for $\widetilde{Q}_k'$. The same argument applies for $\widetilde{Q}_k''$.
\endproof

\section{{\cyr B}-classes and {\cyr B}-Theorem}\label{sec9}

In this Section we prove the \textcyr{B}-Theorem, one of the main results of this paper.

\subsection{Morphism \textcyr{B}}

\begin{defi}
Define the morphism of complex vector spaces $$\textnormal{\textcyr{B}}\colon K_0^{\mathbb T}(\mathbb P^{n-1})_{\mathbb C}\to  H^\Omega_{\mathbb T}(\mathbb P^{n-1})$$ by
\beq
\label{ccmor} 
\textnormal{\textcyr{B}}(F):=\widehat\Gamma^+_{\mathbb P^{n-1}}
\cdot
\exp\left(\pi\sqrt{-1}c_1(\mathbb P^{n-1})\right) \cdot
{\rm Ch}_{\mathbb T}(F).
\eneq
See Section \ref{charclass} for the definition of the characteristic classes in the r.h.s. of \eqref{ccmor}.
\end{defi}

\subsection{\textcyr{B}-Theorem}
Consider the space $\clgr{S}\ _n$ of solutions $I(q,\bm z)$ of the equivariant 
quantum differential equation \eqref{eqde} that are holomorphic 
wrt $\bm z$ in $\Omega$. The space $\clgr{S}\ _n$ is a module over $\mathcal O_\Omega$. Since elements of $\mathcal S _n$ can be seen as element of $\clgr{S}\ _n$, the isomorphism $\theta\colon K_0^\mathbb T(\mathbb P^{n-1})\to\mathcal S_n$ of Corollary \ref{corsn} induces a map
\[ \theta\colon K_0^\mathbb T(\mathbb P^{n-1})\to\clgr{S}\ _n.
\]
The restriction $\mathcal S^o$ of the topological-enumerative morphism,
 defined in \eqref{nn}, allows us to associate an element of $\clgr{S}\ _n$ 
 to any element of $H^\bullet_\mathbb T(\mathbb P^{n-1},\mathbb C)$.
  By extension of scalars, there is an induced morphism
\[
\mathcal S^o\colon H^\Omega_\mathbb T(\mathbb P^{n-1})
\to \clgr{S}\ _n,\quad \alpha\mapsto \mathcal S^o(q)\alpha.
\]

\begin{theorem}[\textcyr{B}-Theorem]\label{bteo}The following diagram is commutative:
\[\xymatrix{
K_0^\mathbb T(\mathbb P^{n-1})_{\mathbb C}\ar[rr]^{\textnormal{\textcyr{B}}}\ar[dr]_{\theta}&& H^\Omega_\mathbb T(\mathbb P^{n-1})\ar[dl]^{\mathcal S^o}\\
&{\clgr{S}\ _n}&}
\]

In other words, if $F\in K_0^{\mathbb T}(\mathbb P^{n-1})_{\mathbb C}$ and $\theta(F)\in\mathcal S_n$ is the corresponding solution to the joint system of equations \eqref{eqde} and \eqref{qkz}, then the meromorphic functions $h_{X,j}(\bm z)$, with $j=1,\dots, n$, defined by the identity
\beq
\theta(F)=\sum_{j=1}^{n}h_{X,j}(\bm z)\Psi_{{\rm top,}j},
\eneq
are the components of the equivariant cohomology class \textnormal{\textcyr{B}}$(F)$ wrt the basis $(x_\alpha)_{\alpha=0}^{n-1}$, i.e.
\beq
\text{\textnormal{\textcyr{B}}}(F)=\sum_{j=1}^nh_{X,j}(\bm z)x_{j-1}.
\eneq
\end{theorem}

\begin{oss}

The relation between $\theta(F)$ and the equivariant-topological solution is the equivariant 
version of part 3.b of \cite[Conjecture 5.2]{CDG1} for $\mathbb P^{n-1}$,
see also \cite{dubro0,KKP,gamma1}. 
 Notice that this is also a refinement of Gamma Theorem of \cite{tar-var}.
\end{oss}

\proof
We prove the statement of the 
theorem for a basis in $ K_0^{\mathbb T}(\mathbb P^{n-1})_{\mathbb C}$. Then
the result for an arbitrary element $F\in K_0^{\mathbb T}(\mathbb P^{n-1})_{\mathbb C}$
follows by linearity and Lemma \ref{lemmach}.

 For $k\in\mathbb Z$ consider the
  basis $([\mathcal O(-k-n+1)],\dots,[\mathcal O(-k-1)],[\mathcal O(-k)])$ 
  in $ K_0^{\mathbb T}(\mathbb P^{n-1})_{\mathbb C}$. Its $\theta$-image in $\mathcal S_n$ 
 is the basis $\left(\Psi^{k+n-1-m}\right)_{m=0}^{n-1}$. 
 Let $Y_{{\rm TV},k}=\left([Y_{{\rm TV},k}]^\lambda_m\right)_{\lambda,m}$
 be the matrix defined by
\[
\Psi^{k+n-1-m}=\sum_\lambda[Y_{{\rm TV},k}]^\lambda_m\ x_\lambda,\quad m=0,\dots,n-1.
\]
For $\bm z\in\Omega$, the matrix $Y_{{\rm TV},k}$ is a fundamental system of solutions 
of system \eqref{eqpn}. 
The matrix $C_{{\rm TV},k}$ connecting the basis $Y_{{\rm TV},k}$ with the topological-enumerative 
solution $Y_{\rm top}$, 
\[Y_{{\rm TV},k}=Y_{\rm top}\cdot C_{{\rm TV},k},
\]
equals
\[
C_{{\rm TV},k}=C\cdot \operatorname{diag}\left(\exp(2\pi\sqrt{-1}(k+n-1-m)z_{m+1})\right)_{m=0}^{n-1},
\]
where $C$ is given by \eqref{conn1}. This shows that $C_{{\rm TV},k}$ is the matrix of the morphism \textcyr{B} with respect the bases $([\mathcal O(-k-n+1)],\dots,[\mathcal O(-k-1)],[\mathcal O(-k)])$ and $(x_\alpha)_{\alpha=0}^{n-1}$. This concludes the proof.
\endproof

\section{Formal solutions of the system of $qDE$ and $qKZ$ equations}\label{sec10}
\subsection{Matrix form of $qDE$ and $qKZ$ difference equations}
The sections 
\[
(q,\bm z)\mapsto x_{\alpha}|_{q,\bm z},
\qquad\alpha=0,\dots, n-1,
 \]
 of the bundle ${\rm pr}^*H$, introduced in Section \ref{secqkzdiff},
 define a trivialization of ${\rm pr}^*H$. This trivialization, allows us write 
 the joint system of the  $qDE$ and $qKZ$ difference equations \eqref{eqde}, \eqref{qkz}
 in matrix form. 

For a basis $I_1(q,\bm z),\dots, I_n(q,\bm z)$ of solutions  to  the  joint system
 \eqref{eqde}, \eqref{qkz}, introduce a matrix 
  $Y(q,\bm z)=(Y^\alpha_m(q,\bm z))_{\alpha,m}$, with $\alpha=0,\dots, n-1$,
   $m=1,\dots, n$,  by the formula:
\beq
I_m(q,\bm z)=\sum_{\alpha=0}^{n-1}Y^\alpha_m(q,\bm z)x_\alpha|_{q,\bm z},\quad m=1,\dots, n.
\eneq
Then $Y(q,\bm z)$ is a fundamental system of solutions of the joint system of equations
\beq\label{mqde}
\frac{d}{dq}Y(q,\bm z)=\mathcal A(q,\bm {z})Y(q,\bm z),
\eneq
\beq\label{mqkz}
Y(q,z_1,\dots, z_i-1,\dots, z_n)=\clgr{K}\ _i(q,\bm z)Y(q,\bm z),\quad
\eneq
where $\mathcal A(q,\bm z)$ is the matrix \eqref{eqpn} attached to the operator $x*_{q,\bm z}$ wrt the basis $(x_\alpha)_{\alpha}$, and the matrix $\clgr{K}\ _i(q,\bm z)$ is the matrix attached to the isomorphism \eqref{qkziso} wrt the basis $(x_\alpha)_\alpha$. 

\begin{oss}\label{qkzmatrgbas}
The sections 
\[
(q,\bm z)\mapsto g_j|_{q,\bm z},
\qquad j=1,\dots,n,
\]
 define another trivialization of ${\rm pr}^*H$. In this trivialization, the $qKZ$ difference equations are
\beq
\label{hh}
\widehat{Y}(q,z_1,\dots,z_i-1,\dots, z_n)=K_i(q,\bm z)\widehat Y(q,\bm z),\quad i=1,\dots,n,
\eneq
where $K_i(q,\bm z)$ are the matrices of $qKZ$-operators \eqref{qkzop} wrt the basis $(g_j)_j$
and the matrix $\widehat{Y}(q,\bm z)$ is defined by 
\[
 I_m(q,\bm z):=\sum_{j=1}^{n}\widehat Y^j_m(q,\bm z)g_j|_{q,\bm z},\quad m=1,\dots, n.
\]

Notice the difference between \eqref{mqkz} and (\ref{hh}).
\end{oss}

In Section \ref{sec6} we studied equation \eqref{mqde} only. Now we will study the joint system of equations
\eqref{mqde} and \eqref{mqkz}. 
As a result of this Section and Section 10 we will deduce the following theorem.
(For its precise statement see in Theorem \ref{gaugeGeqde} and Corollary \ref{cormathKj}.)

\begin{theorem}\label{formred}
Consider the joint system \eqref{mqde},
\eqref{mqkz} of the
$qDE$ and $qKZ$ equations for $\mathbb P^{n-1}$. This system is equivalent at $q=\infty$, up to change of variable $q=s^n$, to the system
\begin{align}
\label{redsyst1}
\frac{dZ}{ds}&=UZ,\\
\label{redsyst2}
Z(s,z_1,\dots,z_j-1,\dots,z_n)&=\mathcal K_jZ(s,\bm z),\quad j=1,\dots,n,
\end{align}
where 
\begin{align*}
U&:={\rm diag}(n\zeta_n^0,\dots,n\zeta_n^{n-1}),\\
\mathcal K_j&:={\rm diag}\left(\zeta_n^0,\zeta_n^{-1},\dots,\zeta_n^{-(n-1)}\right),\quad j=1,\dots,n,
\\
\zeta_n&:=\exp\left(\frac{2\pi\sqrt{-1}}{n}\right).
\end{align*}
\end{theorem}

The theorem says that after a formal transformation, 
the system of $qDE$ and $qKZ$ equations becomes a system with constant coefficients and separated variables. Moreover, 
the system splits into the direct sum of systems of rank one.

System \eqref{redsyst1}, \eqref{redsyst2} admits the basis of solutions
\[Z_i(s,\bm z)=\exp\left(n\zeta_n^{i-1}s+\frac{2(i-1)\pi\sqrt{-1}}{n}\sum_{a=1}^n z_a\right)
\begin{pmatrix}
0\\
\vdots\\
1_i\\
\vdots\\
0
\end{pmatrix} \qquad i=1,\dots, n.
\]
 All solutions of  system \eqref{redsyst1}, \eqref{redsyst2} are linear combinations of these basis solutions 
 with coefficients 1-periodic in $z_1,\dots,z_n$. 

The formal transformation which realizes the reduction to  system  \eqref{redsyst1}, \eqref{redsyst2}, will be described
 in the following subsections.

\subsection{Shearing transformation}The singularity
at  $q=\infty$ of the differential system \eqref{eqpn} is irregular of Poincar\'{e} rank $1$. 
It is known \cite{wasow,BJL79a,div1,div2} that \eqref{eqpn} admits a formal solution of the form
\beq
Y_{\rm form}(q,\bm z)=\Phi\left(q^{\frac{1}{\nu}},\bm z\right)q^{\Lambda(\bm z)}\exp\left(P\left(q^{\frac{1}{\nu}},\bm z\right)\right),
\eneq
where 
\begin{itemize}
\item $\nu\in\mathbb N$ is the \emph{degree of ramification} of the singularity,
\item $\Phi$ is an $n\times n$ matrix-valued  formal power series in 
${q^{-\frac{1}{\nu}}}$ of the form
\[
\Phi\left(q^{\frac{1}{\nu}},\bm z\right)=
\sum_{j=0}^\infty\Phi_j(\bm z)q^{-\frac{j}{\nu}},\quad \det\Phi_0(\bm z)\neq 0,
\]
\item $\Lambda$ is an $n\times n$-matrix depending only on $\bm z$ (the 
\emph{exponent of formal monodromy}),
\item $P={\rm diag}(p_1,\dots, p_n)$ where each $p_j(q^{\frac{1}{\nu}},\bm z)$ is a polynomial in $q^{\frac{1}{\nu}}$ of the form
\[p_j(q^{\frac{1}{\nu}},\bm z)=\sum_{\ell=1}^{N_j}p_{j\ell}(\bm z)q^{\frac{\ell}{\nu}},\quad N_j\geq 1.
\]
\end{itemize}
To find the formal solution $Y_{\rm form}$, we
 perform the gauge transformation of \eqref{eqpn} defined by 
\beq
\label{shearing}
Y(q,\bm z)=\mathcal H(q)\cdot \widetilde{T}(q,\bm z),\quad \mathcal H(q):={\rm diag}\left(1,q^{-\frac{1}{n}},\dots,q^{-\frac{n-1}{n}}\right),
\eneq
called the
\emph{shearing transformation}, see \cite[Section 19]{wasow}. The function $\widetilde T$ satisfies the 
 differential equation
\beq
\label{sheqpn}
\frac{d\widetilde T}{dq}=\mathcal A_{sh}(q,\bm z)\widetilde T,\quad \mathcal A_{sh}=\mathcal H^{-1}\cdot\mathcal A\cdot \mathcal H-\mathcal H^{-1}\frac{d\mathcal H}{dq}.
\eneq
Explicitly, the entries of $\mathcal A_{sh}$ are given by
\begin{align*}
(\mathcal A_{sh})^\alpha_\beta=&\ q^{\frac{1-n}{n}}\delta_{\alpha-\beta,1}+q^{\frac{1-n}{n}}\delta_{\alpha,1}\delta_{\beta,n}\\
&+\sum_{j=1}^n(-1)^{j+1}s_j(\bm z)\delta_{\alpha+j,n+1}\delta_{\beta,n}q^{\frac{1-n-j}{n}}+\delta_{\alpha\beta}\frac{\beta-1}{n}q^{-1},
\end{align*}
for $\alpha,\beta=1,\dots, n$. With the change of variable $q=s^n$, the function $T(s,\bm z):=\widetilde T(s^n,\bm z)$ is a solution of the equation
\beq
\label{eqshea}
\frac{d}{ds} T(s,\bm z)=\mathcal B(s,\bm z) T(s,\bm z),\quad\mathcal B(s,\bm z):=ns^{n-1}\mathcal A_{sh}(s^n,\bm z).
\eneq

\begin{lemma}\label{lemmacoeffb}
We have the following expansion for the coefficient $\mathcal B(s,\bm z)$:
\beq
\mathcal B(s,\bm z)=\mathcal B_0+\frac{1}{s}\mathcal B_1(\bm z)+\sum_{j=2}^n\frac{1}{s^j}\mathcal B_j(\bm z),
\eneq
where
\begin{align}\label{b0}
\mathcal B_0:=&\begin{pmatrix}
0&&\dots&0& n\\
n&0&\dots&0&0\\
0&n&\dots&0&0\\
&&\ddots&&\vdots\\
&&&n&0
\end{pmatrix},
\end{align}
\begin{align}
\label{b1}
\mathcal B_1(\bm z):=&
\begin{pmatrix}
0&&&&&\\
&1&&&&\\
&&2&&&\\
&&&\ddots&&\\
&&&&n-2&\\
&&&&&n-1+ns_1(\bm z)
\end{pmatrix},
\end{align}
\begin{align}
\label{b2}
\mathcal B_j(\bm z):=&\begin{pmatrix}
0&&0&0\\
\vdots&&\vdots&\vdots\\
0&&0&0\\
0&&0&(-1)^{j+1}ns_j(\bm z)\\
0&&0&0\\
\vdots&&\vdots&\vdots\\
0&&0&0
\end{pmatrix},\quad j=2,\dots, n.
\end{align}
\qed
\end{lemma}

By the shearing transformation \eqref{shearing} and the change of variable $q=s^n$, we have reduced the equivariant quantum differential equation of $\mathbb P^{n-1}$ to the equation
\beq\label{redeqpn}
\frac{d T}{ds}=\mathcal B(s,\bm z) T.
\eneq
Equation \eqref{redeqpn} has
 an irregular singularity at $s=\infty$ of Poincar\'{e} rank 1. 
 The essential difference between differential systems \eqref{eqpn} 
 and \eqref{redeqpn}: the matrix  $\mathcal B(s,\bm z)$ 
 has the leading term $\mathcal B_0$ with distinct eigenvalues
\beq
u_k:=n\zeta_n^{k-1},\quad \zeta_n:=\exp\left(\frac{2\pi\sqrt{-1}}{n}\right),\quad k=1,\dots, n,
\eneq
while
the matrix $\mathcal A(q,\bm z)$ has the nilpotent leading term $\mathcal A_0(\bm z)$.

\subsection{The $\mathcal E$-matrix} Let $(e_1,\dots, e_n)$ be the standard basis of $\mathbb C^n$. Let\footnote{Here the subscript ``cl'' stands for \emph{classical}. The matrix $\eta_{\rm cl}$, indeed, appears in the study of the quantum cohomology of $\mathbb P^n$ as the classical Poincar\'e metric. See Remark \ref{remfrobpn}.} $\eta_{\rm cl}$ 
 be the bilinear form on  $\mathbb C^n$ 
 with  matrix 
\beq \label{nepoi}
(\eta_{\rm cl})_{\alpha\beta}=\delta_{n+1,\alpha+\beta}
\eneq 
wrt the standard basis.

For fixed $\bm z\in\mathbb C^n$, consider the $\mathbb C$-linear endomorphisms $\mathcal B_0,\mathcal B_1(\bm z)\in{\rm End}(\mathbb C^n)$ defined, in the standard basis, by the matrices $\mathcal B_0$ and $\mathcal B_1(\bm z)$ of equations \eqref{b0}, \eqref{b1}.  
Introduce the matrix $\mathcal E\in GL(n,\mathbb C)$, 
\beq \left(\mathcal E\right)_{i\alpha }:=\frac{1}{\sqrt{n}}\exp\left(\frac{(i-1)(2\alpha-1)\sqrt{-1}\pi}{n}\right),\quad \alpha,i=1,\dots, n,
\eneq
with inverse 
\beq\label{Pinv} (\mathcal E^{-1})_{\alpha i}=\frac{1}{\sqrt{n}}\exp\left(\frac{(i-1)(1-2\alpha)\sqrt{-1}\pi}{n}\right),\quad \alpha,i=1,\dots, n.
\eneq

\begin{lemma}\label{fbasis}
Define the basis $(f_1,\dots, f_n)$ of $\mathbb C^n$ by
\[f_j:=\sum_{\alpha=1}^n\left(\mathcal E^{-1}\right)_{\alpha j}e_\alpha,\quad j=1,\dots,n,
\]
then
\begin{enumerate}
\item The basis $(f_1,\dots, f_n)$ is orthonormal wrt the bilinear form $\eta_{\rm cl}$.
\item The basis $(f_1,\dots, f_n)$ is an eigenbasis of the operator $\mathcal B_0$.
\item For any fixed $\bm z\in\mathbb C^n$, 
\[\eta_{\rm cl}\left(\mathcal B_1(\bm z)f_i,f_i\right)=s_1(\bm z)+\frac{n-1}{2},\quad i=1,\dots,n.
\]
\end{enumerate}
\end{lemma}

\proof
The statements are equivalent to the identities 
\beq
(\mathcal E^{-1})^T\ \eta_{\rm cl}\ \mathcal E^{-1}=\mathbbm1,
\eneq
\beq\label{propP1}
\mathcal E\ \mathcal B_0\ \mathcal E^{-1}={\rm diag}(u_1,\dots,u_n),
\eneq
\beq\label{propP2}
\mathcal E\ \mathcal B_1(\bm z)\ \mathcal E^{-1}=\left(s_1(\bm z)+\frac{n-1}{2}\right)\cdot \mathbbm1+B^{\rm od}(\bm z),
\eneq
where $B^{\rm od}$ is an off-diagonal matrix, i.e. $(B^{\rm od})_{ii}=0$. A straightforward computation shows the validity of these identities.
\endproof

\begin{oss}\label{remfrobpn}
The 
matrices $\mathcal B_0$, $\eta_{\rm cl}$ and $\mathcal E$ 
appear in the study of the quantum cohomology of $\mathbb P^{n-1}$ seen as a Frobenius manifold, see Appendix \ref{qdeDubr} for details.
\end{oss}

\subsection{Formal reduction of the system of $qDE$ and $qKZ$ equations}
$\quad$

\begin{theorem}\label{gaugeGeqde}There exists a unique $n\times n$-matrix $G(s,\bm z)$, of the form \beq\label{formgauge}
G(s,\bm z)=\mathcal H(s^n)\mathcal E^{-1} F(s,\bm z)s^{\Lambda(\bm z)},
\eneq
with
\begin{align}
\mathcal H(s^n)&={\rm diag}(1,s^{-1},\dots, s^{-(n-1)}),\\
\label{seriesF} F(s,\bm z)&
=\mathbbm 1+\sum_{k=1}^\infty\frac{\widetilde{F}_k(\bm z)}{s^k},
\quad \widetilde{F}_k(\bm z)\text{ polynomials,}\\
\Lambda(\bm z) &=\left(s_1(\bm z)+\frac{n-1}{2}\right)\cdot \mathbbm1,
\end{align}
such that the  transformation
\beq\label{gaugeG}
Y(s^n,\bm z)=G(s,\bm z)Z(s,\bm z)
\eneq 
transforms
the joint system \eqref{mqde}, \eqref{mqkz} of qDE and qKZ equations to the system 
\beq\label{gaugeeqqde}
\frac{dZ}{ds}=UZ,\quad U={\rm diag}(n\zeta_n^0,\dots,n \zeta_n^{n-1}),
\eneq
\beq\label{gaugeqkz}
Z(s,z_1,\dots,z_j-1,\dots, z_n)=\mathcal K_j(\bm z)Z(s,\bm z),\quad j=1,\dots,n,
\eneq
where the matrices $\mathcal K_j(\bm z)$ are diagonal and polynomial in $\bm z$.
\end{theorem}

\proof
The theorem follows from Theorem \ref{teoappb} of Appendix \ref{redjoisyst}, after shearing transformation \eqref{shearing} and change of variables $q=s^n$. Notice that Assumption (1)-(4) of Theorem \ref{teoappb} are satisfied: see Lemma \ref{lemmacoeffb}, Lemma \ref{fbasis}, the expression of $qKZ$-operators $K_i$'s in the $g$-bases \eqref{qkzop} and Remark \ref{qkzmatrgbas}.

The functions $\widetilde{F}_k$ are polynomial in $\bm z$: this
 follows from the procedure described in the proof of Theorem \ref{teoappb} and the fact that
\begin{itemize}
\item the matrices $\mathcal B_0(\bm z),\dots, \mathcal B_n(\bm z)$ are polynomial in $\bm z$,
\item the matrices $U$ and $\mathcal E$ do not depend on $\bm z$.
\end{itemize}
The matrices $\clgr{K}\ _j$'s and $\mathcal K_j$'s are related by the identity
\beq\label{modqkzop}
\mathcal K_j(\bm z):=s F(s,z_1,\dots, z_j-1,\dots,z_n)^{-1} \mathcal E \mathcal H(s^n)^{-1} \clgr{K}\ _j(s^n,\bm z) \mathcal H(s^n) \mathcal E^{-1} F(s,\bm z),
\eneq
for $j=1,\dots, n$.
This implies that the matrices $\mathcal K_j(\bm z)$ are polynomial in $\bm z$.
\endproof

\begin{cor}\label{cormodqkz}
The following identity holds true
\beq
\mathcal K_j(\bm z)=\underset{s=0}{\rm Res}\left(\mathcal E\ \mathcal H(s)^{-1} \clgr{K}\ _j(s^n,\bm z) \mathcal H(s)\ \mathcal E^{-1}\right),
\eneq
for any $j=1,\dots, n$.
\end{cor}

\proof
Since the series $ F(s,\bm z)$ has the form \eqref{seriesF}, from \eqref{modqkzop} we deduce that
\begin{align*}\mathcal K_j(\bm z)=&s\left(\mathbbm1 +O\left(\frac{1}{s}\right)\right)\left(\mathcal E\ \mathcal H(s)^{-1} \clgr{K}\ _j(s^n,\bm z) \mathcal H(s)\ \mathcal E^{-1}\right)\left(\mathbbm1 +O\left(\frac{1}{s}\right)\right)\\
=&s \left(\mathcal E\ \mathcal H(s)^{-1} \clgr{K}\ _j(s^n,\bm z) \mathcal H(s)\ \mathcal E^{-1}\right)\left(\mathbbm1 +O\left(\frac{1}{s}\right)\right).
\end{align*}
Hence, 
\[\mathcal E\ \mathcal H(s)^{-1} \clgr{K}\ _j(s^n,\bm z) \mathcal H(s)\ \mathcal E^{-1}=\frac{\mathcal K_j(\bm z)}{s}+O\left(\frac{1}{s^2}\right).
\]
\endproof

\subsection{Formal solutions of the system of $qDE$ and $qKZ$ equations at $q=\infty$}
Consider the system
\begin{align}
\label{sist1-1}
\frac{dZ}{ds}&=UZ,\\
\label{sist2-1}
Z(s,z_1,\dots,z_j-1,\dots, z_n)&=\mathcal K_j(\bm z)Z(s,\bm z),\quad j=1,\dots,n,
\end{align}
 in Theorem \ref{gaugeGeqde}.

\begin{lemma}
\label{lemmasolc}
Let $C(\bm z)$ be a meromorphic $n\times n$-matrix-valued function on $\mathbb C^n$, regular on $\Omega$ and with non-vanishing determinant. The following conditions are equivalent:
\begin{enumerate}
\item The matrix $C(\bm z)$ is a fundamental system of solutions of equation \eqref{sist2-1} only, over the ring of $1$-periodic functions in $\bm z$;
\item The matrix $Z(s,\bm z):=\exp(sU)C(\bm z)$ is a fundamental system of solutions of the joint system of equations \eqref{sist1-1}, \eqref{sist2-1}, over the ring of $1$-periodic functions in $\bm z$.
\end{enumerate}
\end{lemma}

\proof 
We have
\beq
\label{KU}
 \left[\mathcal K_j(\bm z),\exp(sU)\right]=0,
\eneq
 since both $\exp(sU)$ and $K_j(\bm z)$ are diagonal.

(1) imples (2):
 Let $v(s,\bm z)$
be  a column vector satisfying both \eqref{sist1-1} and \eqref{sist2-1}. Then there exists a unique column vector $c_1(\bm z)$ such that $v(s,\bm z)=\exp(sU)c_1(\bm z)$, since the columns of $\exp(sU)$ 
give a $\mathbb C$-basis of solutions of \eqref{sist1-1}. The vector $c_1(\bm z)$ is a solution of \eqref{sist2-1}
by (\ref{KU}).
Hence  $c_1(\bm z)=C(\bm z)c_2(\bm z)$ for a unique column
 vector $c_2(\bm z)$, which is 1-periodic with respect ot $z_1,\dots,z_n$. 
 This shows that $\exp(sU)C(\bm z)$ is a system of fundamental solutions.

(2) implies (1): Equation (\ref{KU})  easily implies that  $C(\bm z)$ is a solution of the equation \eqref{sist2-1} only, and moreover $C(\bm z)$ is a fundamental solution.
\endproof

\begin{theorem}
\label{teofor} 
Let $C(\bm z)$ be an
$n\times n$-diagonal-matrix-valued function, meromorphic on $\mathbb C^n$, regular on $\Omega$ and
with non-vanishing entries on $\Omega$.
The following conditions are equivalent:

\begin{enumerate}
\item The matrix $C(\bm z)$ is a fundamental system of solutions of the difference equations \eqref{sist2-1}, over the ring of $1$-periodic functions in $\bm z$;
\item There exist meromorphic $n\times n$-matrix valued functions $(F_k(\bm z))_{k=1}^\infty$ 
regular on $\Omega$, such that the matrix 
\beq
\label{formalsolutionC}
Y_{\rm form}(s^n,\bm z)=\mathcal H(s^n)\mathcal E^{-1}{F}(s,\bm z)s^{\Lambda(\bm z)}e^{Us},
\eneq
\beq\label{solforCz}
{F}(s,\bm z)=C(\bm z)+\sum_{j=1}^\infty\frac{{F}_j(\bm z)}{s^j},
\eneq
is a formal solution of the
joint system of the
$qDE$ and $qKZ$ equations for $\mathbb P^{n-1}$, 
\begin{align}
\label{orsist1}\frac{d}{ds}Y(s^n,\bm z)&=ns^{n-1}\mathcal A(s^n,\bm z)Y(s^n,\bm z),\\
\label{orsist2}Y(s^n,z_1,\dots,z_j-1,&\dots,z_n)= \clgr{K}\ _j(s^n,\bm z)Y(s^n,\bm z).
\end{align}
\end{enumerate}
Moreover, if such a formal solution exists, then it is unique.

\end{theorem}

\proof
We have 
\beq\label{commC}
[s^{\Lambda(\bm z)}e^{Us},C(\bm z)]=0,
\eneq
since both $s^{\Lambda(\bm z)}e^{Us}$ and $C(\bm z)$ are diagonal.

Let \beq
G(s,\bm z)=\mathcal H(s^n)\mathcal E^{-1}\!\!\left(\mathbbm 1+\sum_{k=1}^\infty\frac{\widetilde{F}_k(\bm z)}{s^k}\right)\!\!s^{\Lambda(\bm z)}
\eneq
be as in Theorem \ref{gaugeGeqde}.

(1) implies (2): By Lemma \ref{lemmasolc}, the matrix $Z(s,\bm z):=\exp(sU)C(\bm z)$ is a solution of the joint system \eqref{sist1-1}, \eqref{sist2-1}.  By Theorem \ref{gaugeGeqde}, the matrix 
\beq
Y_{\rm form}(s^n,\bm z):=G(s,\bm z)Z(s,\bm z),
\eneq
is a formal solution of the joint system \eqref{orsist1}, \eqref{orsist2}. By \eqref{commC}, $Y_{\rm form}(s^n,\bm z)$ can be re-written in the form \eqref{formalsolutionC}, with $F_k(\bm z):=\widetilde{F}_k(\bm z)C(\bm z)$ for $k\in\mathbb N^*$.

(2) implies (1): By \eqref{commC}, we have
\beq \label{form1}
Y_{\rm form}(s^n,\bm z)=\widetilde{G}(s,\bm z)e^{Us}C(\bm z),
\eneq 
\beq
\widetilde{G}(s,\bm z):=\mathcal H(s^n)\mathcal E^{-1}\!\!\left(\mathbbm 1+\sum_{k=1}^\infty\frac{{F}_k(\bm z)C(\bm z)^{-1}}{s^k}\right)\!\!s^{\Lambda(\bm z)}.
\eneq
Thus, the gauge transformation $Y(s^n,\bm z)=\widetilde{G}(s,\bm z)Z(s,\bm z)$ transforms the $qDE$ \eqref{orsist1} into \eqref{sist1-1}. Hence, it automatically transforms the joint system \eqref{orsist1}, \eqref{orsist2} into \eqref{sist1-1}, \eqref{sist2-1}, see the proof of Theorem \ref{teoappb}. 

The function $Z(s,\bm z):=\exp{(sU)}C(\bm z)$ is a solution of the joint system \eqref{sist1-1}, \eqref{sist2-1}. By Lemma \ref{lemmasolc}, one concludes.

The uniqueness of the formal solution follows from Theorem \ref{gaugeGeqde}: we have $\widetilde{G}(s,\bm z)=G(s,\bm z)$.
\endproof

\begin{es}
Let us consider the case of $\mathbb P^1$. The original system of $qDE$ and $qKZ$ equations is the following
\begin{align*}
\frac{d}{dq}Y(q,\bm z)&=\frac{1}{q}\begin{pmatrix}
0&q-z_1z_2\\
1&z_1+z_2
\end{pmatrix}Y(q,\bm z),\\
Y(q,z_1-1,z_2)&=\frac{1}{q}\left(
\begin{array}{cc}
 -{z_2} & q-{z_1 z_2} \\
 {1} & {z_1}\\
\end{array}
\right)Y(q,\bm z),\\
Y(q,z_1,z_2-1)&=\frac{1}{q}\left(
\begin{array}{cc}
 -{z_1} & q-{z_1z_2} \\
 {1} & {z_2} \\
\end{array}
\right)Y(q,\bm z).
\end{align*}
Through a formal gauge transformation $Y(s^2,\bm z)=G(s,\bm z)Z(s,\bm z)$, the system above can be reduced to the system
\begin{align*}
\frac{d}{ds}Z(s,\bm z)&=\begin{pmatrix}
2&0\\
0&-2
\end{pmatrix}Z(s,\bm z),\\
Z(s,z_1-1,z_2)&=\begin{pmatrix}
1&0\\
0&-1
\end{pmatrix}Z(s,\bm z),\\
Z(s,z_1,z_2-1)&=\begin{pmatrix}
1&0\\
0&-1
\end{pmatrix}Z(s,\bm z).
\end{align*}
The formal gauge $G(s,\bm z)$ is given by
\[G(s,\bm z)=\mathcal H(s^2)\mathcal E^{-1}\left(\mathbbm1+\frac{F_1(\bm z)}{s}+\frac{F_2(\bm z)}{s^2}+\dots\right)s^{\Lambda(\bm z)},
\]where
\[\mathcal H(s^2)={\rm diag}\left(1,\frac{1}{s}\right),\quad \mathcal E^{-1}=\left(
\begin{array}{cc}
 \frac{1}{\sqrt{2}} & -\frac{i}{\sqrt{2}} \\
 \frac{1}{\sqrt{2}} & \frac{i}{\sqrt{2}} \\
\end{array}
\right),\quad \Lambda(\bm z)=\left(\frac{1}{2}+z_1+z_2\right)\mathbbm 1,
\] and the coefficients can be computed recursively as in the proof of Theorem \ref{teoappb}. Here we give just the first coefficient $F_1$: 
\[ F_1=\left(
\begin{array}{cc}
 F_{1,1}& F_{1,2}\\
 F_{1,2}&- F_{1,1}
\end{array}\right),
\]where
\[ F_{1,1}=s_2(\bm z)-\frac{1}{16} \left(2 s_1(\bm z)+1\right){}^2,\quad  F_{1,2}=-\frac{\sqrt{-1}}{8}  \left(2 s_1(\bm z)+1\right).
\]
Notice that Corollary \ref{cormodqkz} allow us to compute directly the coefficients $\mathcal K_1,\mathcal K_2$:
\[\mathcal E\ \mathcal H(s)^{-1} \clgr{K}\ _1(s^2,\bm z)\ \mathcal H(s)\ \mathcal E^{-1}=\left(
\begin{array}{cc}
 \frac{z_1-z_2}{2 s^2}+\frac{1}{s}-\frac{z_1 z_2}{2 s^3} & \frac{\sqrt{-1} \left(z_1+z_2\right)}{2 s^2}-\frac{\sqrt{-1} z_1 z_2}{2 s^3} \\
 -\frac{\sqrt{-1} z_1 z_2}{2 s^3}-\frac{\sqrt{-1} \left(z_1+z_2\right)}{2 s^2} & \frac{z_1-z_2}{2 s^2}+\frac{z_1 z_2}{2 s^3}-\frac{1}{s} \\
\end{array}
\right),
\]
\begin{align*}\mathcal E\ \mathcal H(s)^{-1} &\clgr{K}\ _2(s^2,\bm z)\ \mathcal H(s)\ \mathcal E^{-1}=\\
=&\left(
\begin{array}{cc}
 \frac{-z_1+z_2}{2 s^2}+\frac{1}{s}-\frac{z_1 z_2}{2 s^3} & \frac{\sqrt{-1} \left(z_1+z_2\right)}{2 s^2}-\frac{\sqrt{-1} z_1 z_2}{2 s^3} \\
 -\frac{\sqrt{-1} z_1 z_2}{2 s^3}-\frac{\sqrt{-1} \left(z_1+z_2\right)}{2 s^2} & \frac{-z_1+z_2}{2 s^2}+\frac{z_1 z_2}{2 s^3}-\frac{1}{s} \\
\end{array}
\right).
\end{align*}
By Corollary \ref{cormodqkz}, we obtain that
\[\mathcal K_1(\bm z)=\mathcal K_2(\bm z)=\left(
\begin{array}{cc}
1&0\\
0&-1
\end{array}\right).
\]Notice, in particular, that both $\mathcal K_1$ and $\mathcal K_2$ are equal  constant matrices: in Corollary \ref{cormathKj} we will prove that this is the
 general property valid for all projective spaces $\mathbb P^{n-1}$.
\end{es}

\section{Stokes bases of the system of $qDE$ and $qKZ$ equations}\label{sec11}
\subsection{Stokes rays, Stokes sectors}
The solution $Y_{\rm form}(s^n,\bm z)$ described in Theorem \ref{teofor} is  formal: the series $F(s,\bm z)$ is typically divergent. Nevertheless, $Y_{\rm form}(s^n,\bm z)$ contains information about genuine solutions of the   differential system \eqref{eqpn}. The formal solution prescribes indeed the asymptotics of genuine fundamental solutions of \eqref{eqpn}.

As in Section \ref{asymexps}, we use the coordinates $(r,\phi)$ for the universal cover $\widetilde{\mathbb C^*}$ of the punctured $s$-plane $\mathbb C^*$, see equation \eqref{eqrhophi}.

\begin{defi}
We call  a \emph{Stokes ray} any ray in the universal cover $\widetilde{\mathbb C^*}$ of the $s$-plane, defined by the equation
\beq\label{stokesrays}
\phi=\frac{k}{2n},\quad k\in\mathbb Z.
\eneq
We will denote the ray in \eqref{stokesrays} by $R_k$.
\end{defi}

The meaning of Stokes rays is explained by the following lemma.
\begin{lemma}
A number $\phi\in\mathbb R$ is of the form $\phi=k/2n$ for some $k\in\mathbb Z$, if and only if there are integers $m_1,m_2$ such that ${\rm Re}(\zeta_n^{m_1}s)={\rm Re}(\zeta_n^{m_2}s)$ and $m_1\not\equiv m_2\ ({\rm mod}\ n)$. \qed
\end{lemma}

\begin{defi}
We call \emph{Stokes sector} any open sector in $\widetilde{\mathbb C^*}$ which contains exactly $n$ consecutive Stokes rays $R_k,\dots, R_{k+n}$ for some $k\in\mathbb Z$.
\end{defi}

\begin{lemma}
Any open sector $\mathcal V\subset\widetilde{\mathbb C^*}$ of width $\pi+\delta$, i.e. of the form
\beq
\mathcal V=\left\{s\in\widetilde{\mathbb C^*}\colon a-\frac{1}{2}-\delta<\phi<a\right\},\quad a\in\mathbb R,
\eneq is a Stokes sector for $\delta>0$ sufficiently small.
\end{lemma}

The following theorem follows from the general theory of differential equations.

\begin{theorem}[\cite{wasow,BJL79b,fed,sibook,painkapaev}]\label{teostok}
Let $\mathcal V\subseteq\widetilde{\mathbb C^*}$ be a Stokes sector, and let $Y_{\rm form}(s^n,\bm z)$ denote the unique solution described in Theorem \ref{teofor}. There exists a unique fundamental solution $Y(s^n,\bm z)$ of the  differential system \eqref{eqpn} satisfying the asymptotic condition 
\beq\label{asym} Y(s^n,\bm z)\sim Y_{\rm form}(s^n,\bm z),\quad s\to\infty,\quad s\in\mathcal V,
\eneq uniformly on compact subsets of $\Omega$.  The asymptotic expansion \eqref{asym} can actually be extended to a sector wider than $\mathcal  V$, up to the nearest Stokes rays.
\end{theorem}

\begin{oss}
In the notations of Theorem \ref{teofor}, the precise meaning of the asymptotic relation \eqref{asym} is the following:  
\[
\forall K\Subset \Omega,
\
 \forall h\in\mathbb N,
 \
  \forall \overline{\mathcal V}\subsetneq \mathcal V,
  \
   \exists C_{K,h,\overline{\mathcal V}}>0\colon \text{ if }s\in\overline{\mathcal V}\setminus\left\{0 \right\}\text{ then }\]
   \[
    \sup_{\bm z\in K}\left\| \mathcal E\cdot \mathcal H(s)^{-1}\cdot Y(s,\bm z)\cdot\exp(-s U)s^{-\Lambda(\bm z)}-\sum_{m=0}^{h-1}\frac{F_m(\bm z)}{s^m}\right\|
    <\frac{C_{K,h,\overline{\mathcal V}}}{|s|^h}.
\]
Here $\overline{\mathcal V}$ denotes any unbounded closed sector of $\widetilde{\mathbb C^*}$ with vertex at $0$, and $F_0(\bm z)=C(\bm z)$.  Here, for defining $s^{\Lambda(\bm z)}$, the principal branch of $\log s$ is chosen.
\end{oss}

\begin{lemma}\label{stoksec}
The sectors $\mathcal V_k',\mathcal V_k''$, defined by \eqref{stsec1}-\eqref{stsec2}, are maximal Stokes sectors wrt to the inclusion, i.e.
\begin{enumerate} 
\item they are Stokes sectors,
\item any Stokes sector $\mathcal V$ is contained in one (and only one) $\mathcal V_k'$ or $\mathcal V_k''$.
\end{enumerate}
\qed
\end{lemma}

\subsection{Stokes bases and Stokes matrices}
Consider a basis $(I_i(s^n,\bm y,\bm z))_{i=1}^n$ of solutions of the joint system \eqref{eqde} and \eqref{qkz}, and denote by $Y(s^n,\bm z)$ the corresponding matrix-valued function defined by
\beq
I_\beta(s^n,\bm y,\bm z)=\sum_{\alpha=0}^{n-1}Y(s^n,\bm z)^\alpha_\beta\  x_{\alpha},\quad \beta=1,\dots, n.
\eneq
The function $Y(s^n,\bm z)$ is a fundamental system of solution of the joint system \eqref{mqde} and \eqref{mqkz}. Denote by $Y_{\rm form}(s^n,\bm z)$ the unique formal solution associated to a diagonal matrix $C(\bm z)$ described in Theorem \ref{teofor}.
\begin{defi}
We say that a basis $(I_i(s^n,\bm y,\bm z))_{i=1}^n$ of solutions of the joint system \eqref{eqde} and \eqref{qkz} is a \emph{Stokes basis with normalization $C(\bm z)$} on a sector $\mathcal V$ if it can be reordered in such a way that the corresponding matrix-valued solution $Y(s^n,\bm z)$ satisfies the asymptotic expansion
\beq
Y(s^n,\bm z)\sim Y_{\rm form}(s^n,\bm z),\quad s\to\infty,\quad s\in\mathcal V,
\eneq uniformly on compact subsets of $\Omega$. The matrix $Y(s^n,\bm z)$ is called\footnote{Here we introduce this convenient terminology, though not standard in the literature of ordinary differential equations.} the \emph{Stokes fundamental solution with normalization $C(\bm z)$} on $\mathcal V$ of the joint system \eqref{mqde} and \eqref{mqkz}.
\end{defi}

\begin{oss}By Theorem \ref{teostok}, if $\mathcal V$ is a Stokes sector, and $C(\bm z)$ is a fixed normalization, two Stokes bases on $\mathcal V$ differ only for the order of their objects. Thus, by abuse of language, we will refer to the $\frak S_n$-orbit of Stokes bases on $\mathcal V$ as \emph{the} Stokes basis on $\mathcal V$. Furthermore, if $\mathcal V\subseteq \mathcal V_k'$ (or $\mathcal V_k''$, resp.) then the Stokes basis on $\mathcal V$ is actually the Stokes basis on $\mathcal V_k'$ (or $\mathcal V_k''$, resp.), by Theorem \ref{teostok}.\newline
\end{oss}

Notice that if $\mathcal V$ is a Stokes sector, then also $e^{\pi\sqrt{-1}}\mathcal V$ and $e^{2\pi\sqrt{-1}}\mathcal V$ are Stokes sectors.

\begin{defi}\label{stokesmatr}
Let $Y(s^n,\bm z)$ be the Stokes fundamental solution (with normalization $C(\bm z)$) of system \eqref{eqpn} on the Stokes sector $\mathcal V$. Let $Y_1(s^n,\bm z)$, and $Y_2(s^n,\bm z)$, be the Stokes solutions on $e^{\pi\sqrt{-1}}\mathcal V$ and on $e^{2\pi\sqrt{-1}}\mathcal V$, respectively. Define the \emph{ Stokes matrices} attached to $\mathcal V$ and $C(\bm z)$ as the matrices $\mathbb S_1,\mathbb S_2$ (depending on $\bm z\in\Omega$) for which we have
\beq
Y_1(s^n,\bm z)=Y(s^n,\bm z)\mathbb S_1,\quad Y_2(s^n,\bm z)=Y_1(s^n,\bm z)\mathbb S_2,\quad s\in\widetilde{\mathbb C^*},\quad \bm z\in\Omega.
\eneq
\end{defi}

\subsection{Properties of Stokes matrices and lexicographical order}Let $\mathcal V$ be a Stokes sector, and let $(I_i(s^n,\bm y,\bm z))_{i=1}^n$ be the Stokes basis of the joint system \eqref{eqde} and \eqref{mqkz} on $\mathcal V$ with normalization $C(\bm z)$. Each element $I_i$ corresponds to one eigenvalue $u_j=n\zeta_n^j$, $j=0,\dots,n-1$. Any ordering of the eigenvalues $u_j$'s (i.e. any permutation of the diagonal entries of $U$) corresponds to an ordering of the elements $I_i$'s. Correspondingly, the Stokes matrices $\mathbb S_1$ and $\mathbb S_2$ attached to $\mathcal V$ transform by conjugation by a permutation matrix.

\begin{prop}Denote by $\mathbb S_1,\mathbb S_2$ the Stokes matrices computed wrt the Stokes sector $\mathcal V$.
There exists a unique order of the entries of $U$ such that for all $\bm z\in\Omega$ the matrix $\mathbb S_1$ (resp. $\mathbb S_2$) is upper triangular (resp. lower triangular) with ones along the diagonal. 
\end{prop}

\proof The reader may consult \cite{wasow,BJL79b,painkapaev,div1,div2}. See also \cite{CDG0,CDG1}.
\endproof

The order which realizes the upper triangular form of $\mathbb S_1$ (and consequently the lower triangular form of $\mathbb S_2$) is unique, since $u_i\neq u_j$ for $i\neq j$, and it will be called the \emph{lexicographical order} wrt the Stokes sector $\mathcal V$.

\begin{prop}\label{propstok}
In the  notations of Definition \ref{stokesmatr}, the following identities hold true for any $s\in\widetilde{\mathbb C^*}$ and $\bm z\in\Omega$:
\begin{enumerate}
\item $Y_{2}\left((e^{2\pi\sqrt{-1}}s)^n,\bm z\right)=Y(s^n,\bm z)\cdot \exp\left(2\pi\sqrt{-1}\Lambda (\bm z)\right)$,
\item $Y_2(s^n,\bm z)=Y(s^n,\bm z)\cdot \mathbb S_1\mathbb S_2$,
\item $Y\left((e^{2\pi\sqrt{-1}}s)^n,\bm z\right)=Y(s^n,\bm z)\cdot \exp\left(2\pi\sqrt{-1}\Lambda (\bm z)\right)\cdot \left(\mathbb S_1\mathbb S_2\right)^{-1}$.
\end{enumerate}
Here $\Lambda(\bm z)$ is the exponent of formal monodromy, i.e.
\[\Lambda(\bm z)=\Lambda(\bm z):=\left(s_1(\bm z)+\frac{n-1}{2}\right)\cdot \mathbbm1.
\]
\end{prop}

\proof
For (1), notice that
\[Y_{2}\left((e^{2\pi\sqrt{-1}}s)^n,\bm z\right)\cdot \exp\left(-2\pi\sqrt{-1}\Lambda (\bm z)\right)
\]is a solution of \eqref{eqpn} with asymptotic expansion $Y_{\rm form}(s^n,\bm z)$ on the Stokes sector $\mathcal V$. Hence it must coincide with $Y(s^n,\bm z)$. Point (2) is a direct consequence of the definition of Stokes matrices. Point (3) follows from points (1) and (2).
\endproof

\subsection{Stokes bases $\widetilde{Q}_k'$ and $\widetilde{Q}_k''$}$\quad $

\begin{prop}\label{corstok}
The basis $\widetilde{Q}_k'$ (resp. $\widetilde{Q}_k''$) is a Stokes basis on $\mathcal V_k'$ (resp. $\mathcal V_k''$) with normalization
\beq C(\bm z)=(2\pi)^{\frac{n-1}{2}}e^{-\pi\sqrt{-1}\frac{n-1}{2}}{\rm diag}\left(e^{\frac{m\pi\sqrt{-1}}{n}}(\zeta_n^m)^{s_1(\bm z)+\frac{n-1}{2}}\right)_{m=0}^{n-1}.
\eneq
\end{prop}

\proof
It follows from Corollary \ref{corqtilde}, and formula \eqref{Pinv} for $\mathcal E^{-1}$.
\endproof

\begin{cor}\label{cormathKj}
The operators $\mathcal K_j(\bm z)$, with $j=1,\dots, n$ are all equal and independent of $\bm z$. Indeed, we have 
\[\mathcal K_j={\rm diag}\left(\zeta_n^{-m}\right)_{m=0}^{n-1},\quad j=1,\dots, n.
\]
\end{cor}

\proof
It follows from Theorem \ref{teofor} and the explicit computation $$C(z_1,\dots,z_j-1,\dots,z_n)C(\bm z)^{-1}.$$
\endproof

\subsection{Stokes bases as $\mathbb T$-full exceptional collections}$\quad$

\begin{theorem}\label{stokexcbas}
Via the isomorphism $\theta\colon K_0^{\mathbb T}(\mathbb P^{n-1})_{\mathbb C}\to\mathcal S_n$, Stokes bases on Stokes sectors of the $qDE$ of $\mathbb P^{n-1}$ correspond to $K$-classes of $\mathbb T$-full exceptional collections in $\mathcal D^b_{\mathbb T}(\mathbb P^{n-1})$. 
\end{theorem}

\proof Stokes bases correspond to $\mathbb T$-full exceptional collections of Corollary \ref{corexcol}.
\endproof

In particular, the Stokes basis $\widetilde{Q}_{1-n}''$ corresponds (up to shifts) to the exceptional collection
\beq\label{excol1}\left(\mathcal O\left(\frac{n}{2}\right),\bigwedge\nolimits^1\mathcal T\left(\frac{n}{2}-1\right),\mathcal O\left(\frac{n}{2}+1\right),\bigwedge\nolimits^3\mathcal T\left(\frac{n}{2}-2\right),\dots,\mathcal O(n-1),\bigwedge\nolimits^{n-1}\mathcal T\right)
\eneq
for $n$ even, and
{\footnotesize\beq\label{excol2}\left(\mathcal O\left(\frac{n-1}{2}\right),\mathcal O\left(\frac{n+1}{2}\right),\bigwedge\nolimits^2\mathcal T\left(\frac{n-3}{2}\right),\mathcal O\left(\frac{n+3}{2}\right),\bigwedge\nolimits^4\mathcal T\left(\frac{n-5}{2}\right),\dots,\mathcal O\left(n-1\right),\bigwedge\nolimits^{n-1}\mathcal T\right)
\eneq} for $n$ odd. All other Stokes bases, and corresponding exceptional collections, are obtained by application of a braid of the form
\[\dots\delta_{n,\rm odd}\delta_{n,\rm even}\delta_{n,\rm odd}\delta_{n,\rm even},
\]or
\[\dots\delta_{n,\rm even}^{-1}\delta_{n,\rm odd}^{-1}\delta_{n,\rm even}^{-1}\delta_{n,\rm odd}^{-1}.
\]

\begin{oss}
Exceptional collections \eqref{excol1} and \eqref{excol2} are the natural equivariant lift in $\mathcal D^b_\mathbb T(\mathbb P^{n-1})$ of the exceptional collections of \cite[Corollary 6.11]{CDG1}. Also in this non-equivariant case, these collections are identified with Stokes bases of the $qDE$ of $\mathbb P^{n-1}$ in suitable Stokes sectors, see \cite[Section 6]{CDG1} for details. 
\end{oss}

\begin{oss}\label{excsym}
All the objects on the $\mathbb T$-full exceptional collections attached to Stokes bases are equipped with their natural $\mathbb T$-equivariant structure, restriction of the natural $GL(n,\mathbb C)$-equivariant structure. Under the presentation \eqref{Kpn}, their $K$-theoretical classes in $K_0^\mathbb T(\mathbb P^{n-1})$ are symmetric polynomials wrt the equivariant parameters $\bm Z$.
\end{oss}

\section{Stokes matrices as Gram matrices of exceptional collections}\label{sec12}

\subsection{Musical notation for braids} 
\begin{figure}[ht!]
\centering
\def\svgscale{.5}
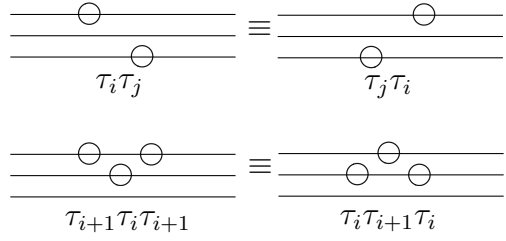
\caption{Braid relations in musical notation.}
\label{braid1}
\end{figure}

We introduce a notation for braids
in $\mathcal B_n$. Elements of $\mathcal B_n$ will be represented as notes on a musical $(n-1)$-line staff. The lines are enumerated from the bottom (1-st line) to the top ($(n-1)$-th line). The generator $\tau_i$ is represented as a hollow oval note head on the $i$-th line. 
The relations defining the braid group $\mathcal B_n$ translate into the diagrammatic rules described in Figure \ref{braid1}.

\subsection{An identity in $\mathcal B_n$}$\quad $
For $n\geq 2$, define the braids 
\beq
\sigma_{n,{\rm even}}:=\underbrace{\delta_{n,{\rm even}}\delta_{n,{\rm odd}}\delta_{n,{\rm even}}\delta_{n,{\rm odd}}\dots}_{n\text{ factors}},
\eneq
\beq
\sigma_{n,{\rm odd}}:=\underbrace{\delta_{n,{\rm odd}}\delta_{n,{\rm even}}\delta_{n,{\rm odd}}\delta_{n,{\rm even}}\dots}_{n\text{ factors}},
\eneq
where $\delta_{n,{\rm even}}$ and $\delta_{n,{\rm odd}}$ are defined in Section \ref{coxetersec}. Set $\sigma_{1,{\rm even}}:=1$, $\sigma_{1,{\rm odd}}=1$.

\begin{lemma}\label{lemmaidentity0}
For any $n\geq 2$, the following identities hold true in $\mathcal B_n$:
\begin{align}
\label{sigma1}
\sigma_{n,{\rm odd}}&=\sigma_{n-1,{\rm odd}}(\tau_{n-1}\tau_{n-2}\dots \tau_{1}),\\
\label{sigma2}
\sigma_{n,{\rm even}}&=\sigma_{n-1,{\rm even}}(\tau_{n-1}\tau_{n-2}\dots \tau_{1}).
\end{align}
\end{lemma}

\proof
We prove \eqref{sigma1} by induction on $n$. For $n=2$, the statement is obvious, being
\[\delta_{2,\rm odd}=\tau_1,\quad\delta_{2, \rm even}=1.
\]

The musical diagram corresponding to the braid $\delta_{n,{\rm odd}}\delta_{n,{\rm even}}\delta_{n,{\rm odd}}\delta_{n,{\rm even}}\dots$ is the following, according to the parity of $n$.

\begin{figure}[ht!]
\centering
\def\svgscale{.6}
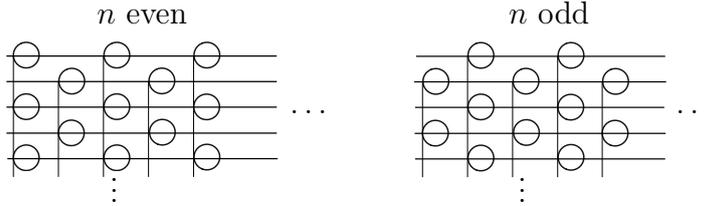
\caption{Diagrammatic notation for the braid $\delta_{n,{\rm odd}}\delta_{n,{\rm even}}\delta_{n,{\rm odd}}\dots$ according to the parity of $n$.}
\label{braid2}
\end{figure}

We collected with a stem the notes corresponding to a same factor $\delta_{n,{\rm odd}}$ (or $\delta_{n,{\rm even}}$). By commutativity, the order of the notes in any factor $\delta_{n,{\rm odd}}$ (or $\delta_{n,{\rm even}}$) can be modified at will, and for this reason we simply collect them with a vertical stem.
In the top line we have 
\begin{itemize}
\item $\frac{n}{2}$ notes, if $n$ is even,
\item $\frac{n-1}{2}$ notes, if $n$ is odd.
\end{itemize}
We call \emph{top factors} those factors $\delta_{n,{\rm odd}}$'s (or $\delta_{n,{\rm even}}$'s) which contains the notes on the top line. In other words, the top factors are
\begin{enumerate}
\item the factors $\delta_{n,{\rm odd}}$'s for $n$ even,
\item the factors $\delta_{n,{\rm even}}$'s for $n$ odd.
\end{enumerate}

The factorization \eqref{sigma1} can be reached by filling the empty spaces between two notes in the last factor $\delta_{n,{\rm odd}}$ (or $\delta_{n,{\rm even}}$), from the bottom to the top line.
We perform this in several steps:
\begin{enumerate}
\item Label by $A_0$ the first (from the left) elementary braid on the $(n-1)$-th line. By a chain of elementary moves, the braid $A_0$ can be moved on the $(n-2)$-th line, towards the right, and can be collected with the next top factor, as described in the following Figure \ref{braid3}. 

\begin{figure}[ht!]
\centering
\def\svgscale{.6}
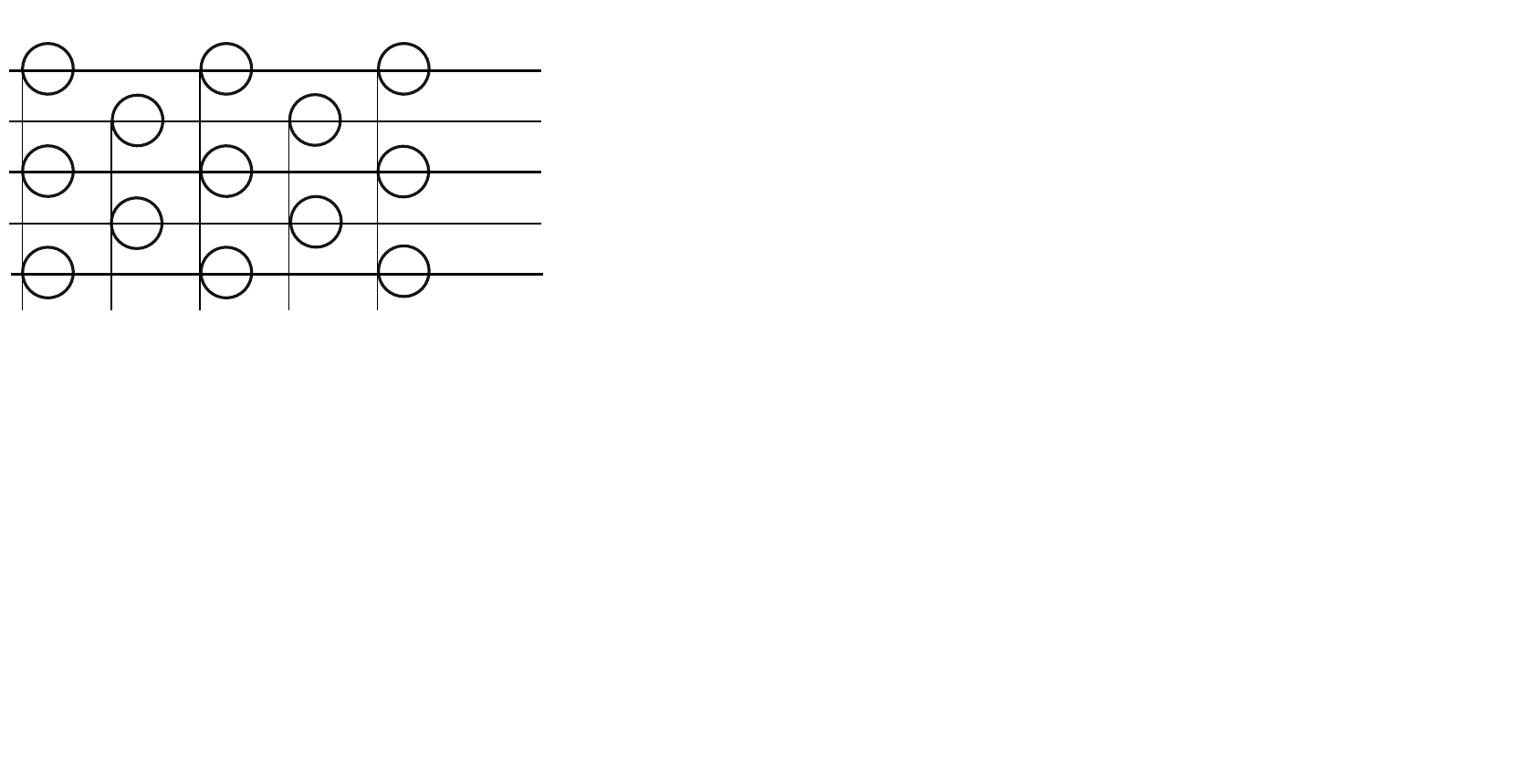
\caption{}
\label{braid3}
\end{figure}

We call $A_1$ the new note obtained from $A_0$. In this way, this factor is ``overcharged'' of notes (i.e. it contains notes $\tau_i$ and $\tau_{i+1}$ for some $i$), and we have an inclination of the stem, the order of the elementary braids being not anymore arbitrary.
\item By the braid relations, the braid $A_1$ can be moved on the $(n-3)$-th line, towards the right, and can be collected with the next factor (not a top factor), as described in the following picture. 

\begin{figure}[ht!]
\centering
\def\svgscale{.6}
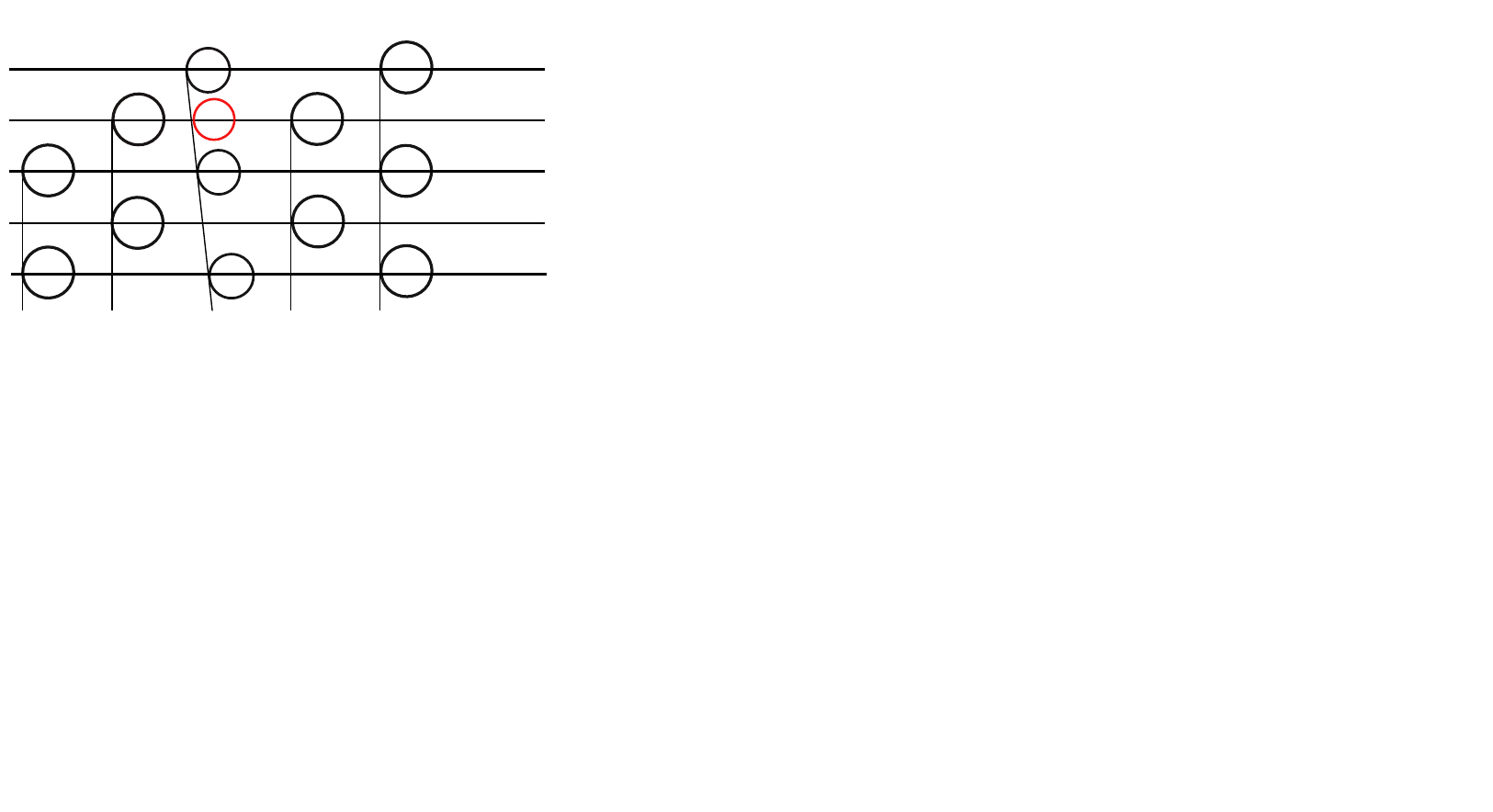
\caption{}
\label{braid4}
\end{figure}

We call $A_2$ the new note obtained from $A_1$. Also in this case, we have an inclination of the stem.
\item Starting from $A_{j}$ on the $(n-1-j)$-th line, iterate the procedure of point (2) in order to produce a new braid $A_{j+1}$ in the line $(n-2-j)$-th line, by overcharging the next factor.
\item Stop when the final braid $A_{j+1}$ fills the empty space on the 
\begin{itemize}
\item 1-st line if $n$ is even,
\item 2-nd line if $n$ is odd.
\end{itemize}
\item Iterate points (1),(2),(3),(4) and stop when the final braid $A_{j+1}$ fills the first empty space from the bottom line to the top.
\end{enumerate}

By applying the procedure above, the factorization \eqref{sigma1} is reached. The argument for \eqref{sigma2} is similar.
\endproof

\begin{es}
Consider $n=7$. The factorization \eqref{sigma2} is obtained with the moves described in Figure \ref{braid5}. For simplicity, we remove all the stems of the notes.

\begin{figure}[ht!]
\centering
\def\svgscale{.6}
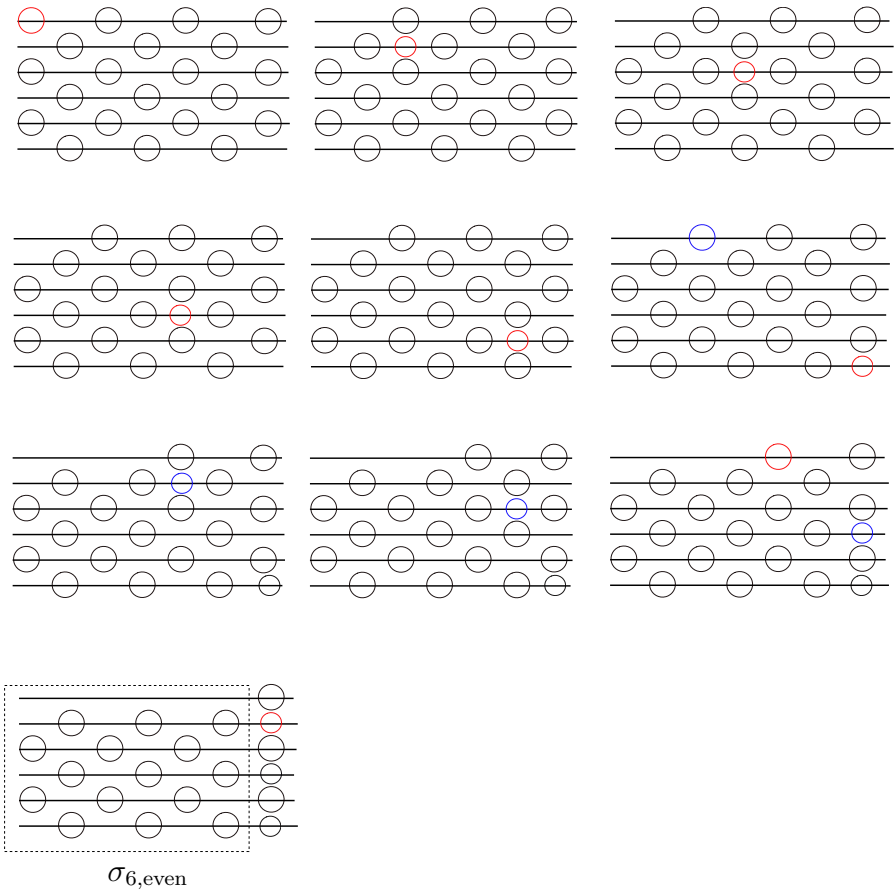
\caption{These are the moves described in the proof of Lemma \ref{lemmaidentity0} in order to obtain the factorization \eqref{sigma2} for $n=7$.}
\label{braid5}
\end{figure}

\end{es}

\begin{cor}\label{lemmaidentity}For $n\geq 2$, we have
\beq
\sigma_{n,{\rm even}}=\beta,\quad
\sigma_{n,{\rm odd}}=\beta,
\eneq
where $\beta$ is given by equation \eqref{braiddual}, i.e.
\[\beta:=\tau_1(\tau_{2}\tau_{1})\dots(\tau_{n-2}\dots \tau_{1})(\tau_{n-1}\tau_{n-2}\dots \tau_{1}).
\]
\qed
\end{cor}

\subsection{Stokes matrices as Gram matrices}$\quad$

\begin{lemma}\label{lemstoksec}
The following identities among Stokes sectors hold true:

\begin{empheq}[left={e^{\pi\sqrt{-1}}\mathcal V'_k=}\empheqlbrace]{align*} 
      &\mathcal V'_{k-\frac{n}{2}},\quad\text{if }n\equiv 0\ ({\rm mod}\ 2),\\
      \\
    & \mathcal V''_{k-\frac{n-1}{2}},\quad\text{if }n\equiv 1\ ({\rm mod}\ 2),
        \end{empheq}
\begin{empheq}[left={e^{\pi\sqrt{-1}}\mathcal V''_k=}\empheqlbrace]{align*} 
&\mathcal V''_{k-\frac{n}{2}},\quad\text{if }n\equiv 0\ ({\rm mod}\ 2),\\
\\
 &\mathcal V'_{k-\frac{n-1}{2}},\quad\text{if }n\equiv 1\ ({\rm mod}\ 2).
        \end{empheq}
\end{lemma}

\proof It is readily obtained from the definition of $\mathcal V_k'$ and $\mathcal V_k''$.
\endproof

\begin{cor}\label{cortrecce}
For any $k\in\mathbb Z$, the Stokes basis  on $e^{\pi\sqrt{-1}}\mathcal V'_k$ is obtained  by acting (on the left) on the Stokes basis $\widetilde{Q}'_k$ with the braid
\beq
\underbrace{\dots\delta_{n,{\rm odd}}\delta_{n,{\rm even}}\delta_{n,{\rm odd}}\delta_{n,{\rm even}}\delta_{n,{\rm odd}}}_{n\text{ factors}}.
\eneq
For any $k\in\mathbb Z$, the Stokes basis on $e^{\pi\sqrt{-1}}\mathcal V''_k$ is obtained  by acting (on the left) on the Stokes basis $\widetilde{Q}''_k$ with the braid
\beq
\underbrace{\dots\delta_{n,{\rm even}}\delta_{n,{\rm odd}}\delta_{n,{\rm even}}\delta_{n,{\rm odd}}\delta_{n,{\rm even}}}_{n\text{ factors}}.
\eneq
\end{cor}

\proof
It follows from the definition of $\widetilde{Q}'_k$ and $\widetilde{Q}''_k$ (see also diagram \eqref{diagQ}), from Proposition \ref{corstok} and Lemma \ref{lemstoksec}.
\endproof

\begin{cor}
Let $\mathcal V$ be a Stokes sector of the $qDE$ of $\mathbb P^{n-1}$, and let $\frak E$ be the $\mathbb T$-full exceptional collection corresponding to the Stokes basis on $\mathcal V$. The exceptional collection $\frak E'$ corresponding to the Stokes sector $e^{2\pi\sqrt{-1}}\mathcal V$ is a foundation of the helix generated by $\frak E$. More precisely, $\frak E'$ is the adjacent foundation on the right of $\frak E$. The collection $\frak E'$ is obtained by application of the inverse Serre functor to objects of $\frak E$: 
\beq
(E_1,\dots,E_n)\mapsto \left(E_1\otimes (\omega^{\mathbb T}_{\mathbb P^{n-1}})^{-1}[-n+1],\dots,E_n\otimes (\omega^{\mathbb T}_{\mathbb P^{n-1}})^{-1}[-n+1]\right),
\eneq
where $\omega^{\mathbb T}_{\mathbb P^{n-1}}$ denotes the $\mathbb T$-equivariant canonical sheaf of $\mathbb P^{n-1}$.
\end{cor}

\proof
By Corollary \ref{cortrecce}, $\frak E'$ is obtained from $\frak E$ by mutation either with the braid 
\beq
\underbrace{\dots\delta_{n,{\rm odd}}\delta_{n,{\rm even}}\delta_{n,{\rm odd}}\delta_{n,{\rm even}}}_{2n\text{ factors}},
\eneq
or with the braid 
\beq
\underbrace{\dots\delta_{n,{\rm even}}\delta_{n,{\rm odd}}\delta_{n,{\rm even}}\delta_{n,{\rm odd}}}_{2n\text{ factors}}.
\eneq
In both cases, by Corollary \ref{lemmaidentity}, the resulting braid is $\beta^2$, where
\[\beta:=\tau_1(\tau_{2}\tau_{1})\dots(\tau_{n-2}\dots \tau_{1})(\tau_{n-1}\tau_{n-2}\dots \tau_{1}).
\]It is well-known that $\beta^2=(\tau_1\dots\tau_{n-1})^{n}$ (see e.g. Theorem 1.24 of \cite{braids}). The result follows from Proposition \ref{serreexc}.
\endproof

\begin{theorem}\label{teostokgram}
Let $\mathbb S_1,\mathbb S_2$ the Stokes matrices computed wrt a Stokes sector $\mathcal V$ in lexicographical order. Let $\varepsilon$ be the exceptional basis of $K_0^{\mathbb T}(\mathbb P^{n-1})_{\mathbb C}$ associated with the Stokes basis on $\mathcal V$ via the isomorphism $\theta\colon K_0^\mathbb T(\mathbb P^{n-1})_{\mathbb C}\to \mathcal S_n$ defined in \eqref{isoth}, and let $\mathcal G$ be the Gram matrix of $\chi^\mathbb T$ wrt $\varepsilon$. Let $J$ be the anti-diagonal matrix
\[J_{\alpha\beta}:=\delta_{\alpha+\beta,n+1},\quad \alpha,\beta=1,\dots, n.
\]
\begin{enumerate}
\item The Stokes matrix $\mathbb S_1$ is equal to the Gram matrix of  $\chi^\mathbb T$ wrt the left dual exceptional basis $^\vee\varepsilon$, i.e.
\beq
\mathbb S_1=J\left(\mathcal G^\dag\right)^{-1}J.
\eneq
\item The matrix $J\mathbb S_2J$ is equal to the Gram matrix of $\chi^\mathbb T$ wrt the exceptional basis $\varepsilon$, i.e.
\beq
\mathbb S_2=J\mathcal G J.
\eneq
\end{enumerate}
\end{theorem}

\proof
By Lemma \ref{stoksec}, $\mathcal V$ is contained in one (and only one) $\mathcal V'_k$ or $\mathcal V''_k$. By Corollary \ref{cortrecce}, the Stokes basis on $e^{\sqrt{-1}\pi}\mathcal V$ is obtained from the Stokes basis on $\mathcal V$ by applying either the braid 
\[\underbrace{\delta_{n,{\rm odd}}\delta_{n,{\rm even}}\delta_{n,{\rm odd}}\delta_{n,{\rm even}}\dots}_{n\text{ factors}},
\]or the braid
\[\underbrace{\delta_{n,{\rm even}}\delta_{n,{\rm odd}}\delta_{n,{\rm even}}\delta_{n,{\rm odd}}\dots}_{n\text{ factors}}.
\]
Consequently, by Corollary \ref{lemmaidentity}, the exceptional basis associated to the Stokes basis on $e^{\sqrt{-1}\pi}\mathcal V$ is $^\vee\varepsilon$. Both points (1) and (2) then follow from Proposition \ref{propconndual}, more precisely from the second identity \eqref{dual1}.
\endproof

\begin{cor}
Let $\mathbb S_1$ and $\mathbb S_2$ be the Stokes matrices computed wrt a Stokes sector $\mathcal V$. We have
\beq\label{ssdaga}
\mathbb S_2=\left(\mathbb S_1^\dag\right)^{-1}.
\eneq
\qed
\end{cor}

\begin{oss}
In theory of Frobenius manifolds \cite{dubro1,dubro2,dubro0,CDG}, the Stokes matrices of the associated isomonodromic system of differential equations  
satisfy an analogous identity, in which the $\dag$-operator is replaced by transposition, see \cite[Theorem 4.3]{dubro2}.
\end{oss}

\begin{cor}
Let $\mathbb S_1$ and $\mathbb S_2$ be the Stokes matrices computed wrt a Stokes sector $\mathcal V$. Both $\mathbb S_1$ and $\mathbb S_2$ have entries in the ring of symmetric Laurent polynomials with integer coefficients, i.e. $\mathbb S_1,\mathbb S_2\in M_n\left(\mathbb Z[\bm Z^{\pm 1}]^{\frak S_n}\right)$. 
\end{cor}

\proof
By Lemma \ref{chiO}, the Gram matrix associated with the Beilinson exceptional collection is with integer symmetric Laurent polynomial entries. The braid group action preserves this property.
\endproof

\begin{cor}
Let $\mathbb S$ be a Stokes matrix of the differential system \eqref{eqpn} computed wrt a Stokes sector $\mathcal V$. Then
\beq\label{dioph2}
\det\left(\lambda\cdot\mathbbm 1-\mathbb S^\dag\mathbb S^{-1}\right)=\sum_{j=0}^n(-1)^j\lambda^{n-j}s_j\left((-1)^{n+1}\frac{Z_1^n}{s_n(\bm Z)},\dots,(-1)^{n+1}\frac{Z_n^n}{s_n(\bm Z)}\right).
\eneq
\end{cor}

\proof
The Corollary immediately follows from Theorem \ref{propdioph} and Theorem \ref{teostokgram}. For an alternative proof (purely analytical), notice that without loss of generality we can assume that $\mathbb S$ is the matrix $\mathbb S_1$ computed wrt $\mathcal V$.
From point (3) of Proposition \ref{propstok}, and equation \eqref{ssdaga} we deduce that
\[\exp(2\pi\sqrt{-1}\Lambda(\bm z))\mathbb S^\dag\mathbb S^{-1}=M_0(\bm z)^n,
\]where $M_0(\bm z)$ is the monodromy operator of differential system \eqref{eqpn}. From Corollary \ref{monoeig} we deduce the constraint
\beq
\det\left(\lambda\cdot\mathbbm 1-(-1)^{n+1}s_n(\bm Z)\mathbb S^\dag\mathbb S^{-1}\right)=\sum_{j=0}^n(-1)^j\lambda^{n-j}s_j(\bm Z^n),
\eneq
which is easily seen to be equivalent to equation \eqref{dioph2}.
\endproof

\section{Specialization of the $qDE$ at roots of unity}
\label{sec13}

\subsection{Specialization of equivariant $K$-theory at roots of unity} 
Fix the equivariant parameters in $K_0^{\mathbb T}(\mathbb P^{n-1})$  by setting
\beq\label{locuspoli}
Z_m=\zeta_n^{m-1},\quad m=1,\dots, n.
\eneq 
Denote  by $K_{\bm \zeta}$ this specialization of $K_0^{\mathbb T}(\mathbb P^{n-1})$.

\begin{theorem}[{\cite[Theorem 1.1]{poli} }]$\quad$\begin{enumerate}
\item The Grothendieck-Euler-Poincar\'{e} pairing $\chi^{\mathbb T}$ specializes to an Hermitian positive definite form $\chi_\zeta$ on $K_{\bm \zeta}$. 
\item If $E$ is an exceptional object in $\mathcal D^b(\mathbb P^{n-1})$ equipped with a $\mathbb T$-equivariant structure, then the class $[E]$ in $K_{\bm \zeta}$ has length 1 wrt the Hermitian form $\chi_{\bm \zeta}$.
\item If $(E_1,E_2)$ is an exceptional pair in $\mathcal D^b(\mathbb P^{n-1})$, with both $E_1$ and $E_2$ equipped with a $\mathbb T$-structure, then the classes $[E_1],[E_2]$ in $K_{\bm \zeta}$ are orthogonal wrt $\chi_{\bm \zeta}$.
\item If $(E_1,\dots, E_n)$ is a full exceptional collection in $\mathcal D^b(\mathbb P^{n-1})$, with each $E_i$ equipped with a $\mathbb T$-structure, then each unit vector in $K_{\bm \zeta}$ is of the form $\pm\zeta_n^k[E_i]$ for some $i$ and some $k$.
\item The action of the braid group on the set of orthonormal exceptional bases of $K_{\bm \zeta}$ 
reduces to the action by permutations of basis vectors.
\end{enumerate}
\end{theorem}

\subsection{Identities for Stirling numbers}
The Stirling numbers of the first kind ${n\brack k}$ are defined recursively by 
\beq\label{stirl1}
{n+1\brack k}=n{n\brack k}+{n\brack k-1},
\eneq
for $k>0$, with the initial conditions 
\beq
{0\brack 0}=1,\quad {0\brack n}={n\brack 0}=0,\quad n>0.
\eneq
The Stirling numbers of the second kind ${n\brace k}$ are defined recursively by
\beq\label{stirl2}
{n+1\brace k}=k{n\brace k}+{n\brace k-1},
\eneq
for $k>0$, with the initial conditions 
\beq
{0\brace 0}=1,\quad {0\brace n}={n\brace 0}=0,\quad n>0.
\eneq
The Stirling numbers of the first and  second kind are related by the identity
\beq\label{stirl12}
\sum_{j\geq 0}(-1)^{n-j}{n\brace j}{j\brack k}=\delta_{nk}.
\eneq

\begin{lemma}Let $n\geq 2$ and $1\leq k\leq n$. We have
\beq\label{s1}
s_k\left(0,\frac{1}{n},\frac{2}{n},\dots,\frac{n-1}{n}\right)=\frac{1}{n^k}{n\brack n-k},
\eneq
\beq \label{s2}
m_k\left(0,\frac{1}{n},\frac{2}{n},\dots,\frac{n-1}{n}\right)=\frac{1}{n^k}{n+k-1\brace n-1}.
\eneq
\end{lemma}

\proof It is sufficient to prove the identities 
\beq\label{ss1} s_k(0,1,\dots,n-1)={n\brack n-k},
\eneq
\beq\label{ss2}
m_k(0,\dots,n-1)={n+k-1\brace n-1}.
\eneq
They are proved by induction on $n$. For $n=2$ both equations \eqref{ss1} and \eqref{ss2} hold true. 
Recall the  following recurrence equations: for $k\geq 2$
\begin{align}\label{st1}
s_k(z_1,\dots,z_n)&=s_k(z_1,\dots,z_{n-1})+z_ns_{k-1}(z_1,\dots,z_{n-1})\\
\label{st2}
m_k(z_1,\dots,z_n)&=m_k(z_1,\dots,z_{n-1})+z_nm_{k-1}(z_1,\dots,z_{n}).
\end{align}
Now equation \eqref{ss1} follows from \eqref{st1} and  \eqref{stirl1}, equation  \eqref{ss2} follows from \eqref{st2} and  \eqref{stirl2}.
\endproof

\begin{lemma}
If $\vartheta_s:=s\frac{d}{ds}$, then
\beq\label{diffe1}
\vartheta_s^n=\sum_{j=1}^n{n\brace j}s^j\frac{d^j}{ds^j},
\eneq
\beq\label{diffe2}
s^n\frac{d^n}{ds^n}=\sum_{j=1}^n(-1)^{n-j}{n\brack j}\vartheta_s^j.
\eneq
\end{lemma}

\proof
Identity \eqref{diffe1} is easily proved by induction on $n$. Identity \eqref{diffe2} follows from \eqref{diffe1} and 
 \eqref{stirl12}.
\endproof

\subsection{Scalar equivariant quantum differential equation at roots of unity}
Consider the specialization of the equivariant parameters $\bm z$ in
$H^\bullet_\mathbb T(\mathbb P^{n-1},\mathbb C)$ defined by the equations
\beq
\exp(2\pi\sqrt{-1}z_m)=\zeta_n^{m-1},\quad m=1,\dots,n.
\eneq
These equations define the locus
\beq
\mathcal P:=\left\{\bm z\in\mathbb C^n\colon 	\bm z=\left(k_1,k_2+\frac{1}{n},\dots,k_n+\frac{n-1}{n}\right),\quad \bm k\in\mathbb Z^n\right\}.
\eneq
We have a distinguished point $\bm z_o\in\mathcal P$, 
\beq\label{zo}
\bm z_o:=\left(0,\frac{1}{n},\dots,\frac{n-1}{n}\right).
\eneq

\begin{theorem}
At $\bm z=\bm z_o$ the scalar equivariant quantum differential equation \eqref{scaleqpn} of $\mathbb P^{n-1}$ for the function $\phi(q)$ reduces to the linear differential equation with constant coefficients,
\beq\label{bellaeq}
\frac{d^n}{ds^n}\varphi(s)=n^n\varphi(s),\quad \varphi(s):=	\phi(s^n).
\eneq
\end{theorem}

\proof By the change of variable $q=s^n$, equation \eqref{scaleqpn} reduces to
\beq
\left[\frac{1}{n^n}\vartheta_s^n+\sum_{j=1}^{n-1}(-1)^{n-j}s_{n-j}(\bm z)\frac{1}{n^j}\vartheta_s^j-\left(s^n+(-1)^{n-1}s_n(\bm z)\right)\right]\phi(s^n)=0.
\eneq
If $\bm z=\bm z_o$, then  the equation reduces to
\beq
\left[\frac{1}{n^n}\vartheta_s^n+\sum_{j=1}^{n-1}(-1)^{n-j}\frac{1}{n^n}{n\brack j}\vartheta_s^j-s^n\right]\phi(s^n)=0,
\eneq
by identity \eqref{s1}.
Using identity \eqref{diffe2}, we obtains the equation
\beq
\left(s^n\frac{d^n}{ds^n}-n^ns^n\right)\phi(s^n)=0.
\eneq
\endproof

Equation \eqref{bellaeq} has two natural bases of solutions:
\begin{enumerate}
\item the basis $(f_m(s))_{m=0}^{n-1}$\,,
\beq
f_m(s):=\exp(n\zeta_n^ms),\quad m=0,\dots,n-1;
\eneq
\item the basis $(g_m(s))_{m=0}^{n-1}$\,,
\beq
g_m(s):=\sum_{k=0}^\infty\frac{(ns)^{m+kn}}{(m+kn)!},\quad m=0,\dots,n-1.
\eneq
\end{enumerate}
The functions $g_m(s)$ are real-valued for $s\in\mathbb R$ and define a partition of the exponential function $e^{ns}$, 
\beq
\sum_{m=0}^{n-1}g_m(s)=e^{ns}.
\eneq

\begin{lemma}\label{lemmasolfg}
The cyclic group $\mathbb Z/n\mathbb Z$ acts on the space of solutions of equation \eqref{bellaeq} via the transformations $T_k\colon s\mapsto \zeta_n^ks$, $k=1,\dots, n$. The basis $(f_m(s))_{m=0}^{n-1}$ is cyclically permuted by this action, while the basis $(g_m(s))_{m=0}^{n-1}$ is an eigenbasis.
\qed
\end{lemma}

Introduce the matrices

\beq
\hat{Y}_f(s^n):=\left(\hat{Y}_f(s^n)^h_m\right)_{h,m=0,\dots, n-1},\quad \hat{Y}_f(s^n)^h_m:=\frac{1}{n^h}\vartheta_s^hf_m(s),
\eneq

\beq
\hat{Y}_g(s^n):=\left(\hat{Y}_g(s^n)^h_m\right)_{h,m=0,\dots, n-1},\quad \hat{Y}_g(s^n)^h_m:=\frac{1}{n^h}\vartheta_s^hg_m(s).
\eneq
Both $\hat{Y}_f(s^n)$ and $\hat{Y}_g(s^n)$ are solutions of the differential equation \eqref{eqpndual}, specialized at $\bm z=\bm z_o$.

 \begin{prop}
 The matrix-valued function $\eta(\bm z_o)^{-1}\hat{Y}_f(s^n)$ is a fundamental system of solutions of \eqref{eqpn} of the form
 \beq\eta(\bm z_o)^{-1}\hat{Y}_f(s^n)=G(s)\exp(sU), 
 \eneq where  $U=\text{diag}(n\zeta^0_n,\dots,n\zeta^{n-1}_n)$ and 
  $G(s)$ is a polynomial in $s$ of degree $n-1$. Hence, $\eta(\bm z_o)^{-1}\hat{Y}_f(s^n)$ is a Stokes basis of the equivariant $qDE$ \eqref{eqpn} at $\bm z=\bm z_o$ for any Stokes sector. In particular, the corresponding formal series $F(s,\bm z_o)$ of the form \eqref{solforCz} is actually convergent.
 \end{prop}
 
 \proof
 The matrix $\eta(\bm z_o)^{-1}\hat{Y}_f(s^n)\eta(\bm z_o)$ is a fundamental system of solutions of \eqref{eqpn} by the discussion in Section \ref{secscaleqpn}. Hence $\eta(\bm z_o)^{-1}\hat{Y}_f(s^n)$ also is 
 a fundamental system of solutions. The matrix $G(s)$ is given by
 \beq
 G(s)=\eta(\bm z_o)^{-1}\cdot L(s),\quad L(s)=\left(L(s)^h_m\right)_{h,m=0}^{n-1},\quad L(s)^h_m:=\zeta_n^h s^h.
 \eneq
 Thus the series $F(s,\bm z_o)$ is given by
 \beq
 F(s,\bm z_o)=\mathcal D_q\mathcal H(s^n)^{-1}\eta(\bm z_o)^{-1}L(s)s^{1-n},
 \eneq
 and is  convergent. The normalization of the Stokes basis can be readily computed from this formula.
 \endproof

 \begin{prop}
 The matrix-valued function $\eta(\bm z_o)^{-1}\hat{Y}_g(s^n)$ is a fundamental system of solutions of \eqref{eqpn} of the form
 \beq
 \eta(\bm z_o)^{-1}\hat{Y}_g(s^n)=Y_o(s^n,\bm z_o)\cdot C,
 \eneq
 where the matrix $Y_o(s^n,\bm z)$ is the Levelt solution defined in Corollary \ref{solstand} and $C$ is a diagonal matrix.
 \end{prop}

 \proof
 The proposition follows from Lemma \ref{lemmasolfg} and Corollary \ref{monoeig}. We leave to the reader the explicit computation of the matrix $C$.
 \endproof

\begin{theorem}
The following conditions are equivalent:
\begin{enumerate}
\item $\bm z'\in\mathcal P$;
\item the formal gauge transformation $G(s,\bm z')$ of Theorem \ref{gaugeGeqde} is convergent;
\item the Stokes phenomenon of the  differential system \eqref{eqpn} specialized at $\bm z=\bm z'$ is trivial, i.e. all the Stokes matrices $\mathbb S(\bm z')$ for all Stokes sectors
are the identity matrix;
\item the monodromy matrix $M_0(\bm z')$ of the equivariant quantum differential equation \eqref{eqpn} specialized at $\bm z=\bm z'\in\Omega$ has order $n$.
\end{enumerate}
\end{theorem}

\proof
We prove that $(1)\Rightarrow (2)\Rightarrow (4)\Rightarrow (1)$. 

\smallskip
Assume $\bm z'=\left(k_1,k_2+\frac{1}{n},\dots, k_n+\frac{n-1}{n}\right)$ for some $\bm k\in\mathbb Z^n$. Then, we have to show that the series $F(s,\bm z)$ in \eqref{formgauge} is convergent for $\bm z=\bm z'$. From the identity \eqref{modqkzop} we deduce that
\[F(s,z_1,\dots, z_j-1,\dots,z_n)=W_j(s,\bm z) F(s,\bm z)\mathcal K_j^{-1},
\]
where
\[W_j(s,\bm z):=s\mathcal D_{\rm q} \mathcal H(s^n)^{-1} K_j(s^n,\bm z) \mathcal H(s^n) \mathcal D_{\rm q}^{-1},\quad j=1,\dots,n.
\]
Hence, we have
\[F(s,\bm z')=\left(\prod_{j=1}^nW_j(s,\bm z_o)^{-k_j}\right)F(s,\bm z_o)\left(\prod_{j=1}^n\mathcal K_j^{-k_j}\right),
\]
and the convergence of $F(s,\bm z')$ follows from the convergence of $F(s,\bm z_o)$.

If (2) holds then $Y(s^n,\bm z')=G(s,\bm z')e^{sU}$ is a solution of system \eqref{eqpn} at $\bm z=\bm z'$, and the transformation $s\mapsto \zeta_ns$ cyclically permutes the diagonal entries of $U$.
Thus (4) holds true.

If (4) holds true, then from Corollary \ref{monoeig} we deduce that $\exp(2\pi\sqrt{-1}z'_j)$ is a $n$-th root of unity, i.e. $\bm z'\in\mathcal P$.

The equivalence of (2) and (3) is obvious.
\endproof

\appendix
\section{Formal reduction of the joint system }
\label{redjoisyst} 

Consider a joint system of differential and difference equations
\begin{align}
\label{jseq1}\frac{d}{ds}X(s,\bm z)&=A(s,\bm z)X(s,\bm z),\\
\label{jseq2}X(s,z_1,\dots, z_i-1,\dots, z_n)&=P_i(s,\bm z)X(s,\bm z),\quad i=1,\dots,n,
\end{align}
where $A, P_i$ are meromorphic $m\times m$-matrix valued functions of $(s,\bm z)\in \mathbb C\times \mathbb C^n$.

Assume that equations \eqref{jseq1}, \eqref{jseq2} are compatible,
\begin{align}\label{comp1}
\frac{d}{ds}P_i(s,\bm z)=&A(s,z_1,\dots, z_i-1,\dots,z_n)P_i(s,\bm z)-P_i(s,\bm z)A(s,\bm z),
\quad i=1,\dots,n,
\\
\label{comp2}
P_i(s,z_1,\dots, z_i,\dots&,z_j-1,\dots,z_n)P_j(s,\bm z)=P_j(s,z_1,\dots, z_i-1,\dots,z_j,\dots,z_n)P_i(s,\bm z), 
\end{align}
for all $i, j$.
 Assume that 
\begin{enumerate}
\item 
the matrices $A(s,\bm z), P_i(s,\bm z)$ have the following convergent power series expansions
\beq
A(s,\bm z)=\sum_{k=0}^\infty A_k(\bm z)\frac{1}{s^k},\quad P_i(s,\bm z)=\sum_{k=0}^\infty P_k(\bm z)\frac{1}{s^k},\quad |s|>\rho,
\eneq
where the matrices $A_k(\bm z), P_k(\bm z)$ are holomorphic functions of $\bm z$ and $\rho>0$;
\item the matrix $A_0(\bm z)$ is diagonalizable, 
\beq
D(\bm z)\cdot A_0(\bm z)\cdot D(\bm z)^{-1}=U(\bm z),\quad U(\bm z):={\rm diag}(u_1(\bm z),\dots, u_n(\bm z)),
\eneq
with $D(\bm z)$ a holomorphic matrix; 
\item the matrix  $U(\bm z)$ is {1-periodic}, 
\beq
U(z_1,\dots, z_i-1,\dots, z_n)=U(\bm z),\quad i=1,\dots,n;
\eneq
\item the eigenvalues $u_i(\bm z)$'s are pairwise distinct for values of $\bm z$ in
 an open subset $W\subset \mathbb C^n$.
\end{enumerate} 
Introduce the diagonal matrix $\Lambda_1(\bm z)$ by the formula
\beq
\label{L1}
[\Lambda_1(\bm z)]_{ij}:=\left[D(\bm z)\cdot A_1(\bm z)\cdot D(\bm z)^{-1}\right]_{ij}\delta_{ij}.
\eneq

Assumptions (1-4) imply that for all $\bm z\in W$ the point $s=\infty$ is an irregular singularity (of Poincar\'e rank 1) of the differential equation \eqref{jseq1}.

\begin{theorem}\label{teoappb}

If the joint  system of equations  \eqref{jseq1}, \eqref{jseq2} satisfies assumptions $(1$-$4)$,
then there exists a {unique} 
$m\times m$-matrix $G(s,\bm z)$ of the form 
\beq
G(s,\bm z)=D(\bm z)^{-1}F(s,\bm z)s^{\Lambda_1(\bm z)},
\eneq
where $F(s,\bm z)$ is a formal power series
\beq
F(s,\bm z) =\mathbbm 1+\sum_{k=1}^\infty F_k(\bm z)\frac{1}{s^k},
\eneq
with $F_k(\bm z)$  regular $m\times m$-matrix-valued functions 
on $W$, such that the change of variables  $X(s,\bm z)=G(s,\bm z)Z(s,\bm z)$
 transforms system \eqref{jseq1}, \eqref{jseq2} into the joint system
\begin{align}
\label{tjseq1}
\frac{d}{ds}Z(s,\bm z)&=U(\bm z)Z(s,\bm z),\\
\label{tjseq2}
Z(s,z_1,\dots, z_i-1,\dots, z_n)&=\mathcal P_i(\bm z)Z(s,\bm z),\quad i=1,\dots,n,
\end{align}
where $\mathcal P_i$ are  diagonal matrices.
\end{theorem}

\proof
First we prove that there exists a unique formal transformation
$X(s,\bm z)$ 
\linebreak
$=G(s,\bm z) Z(s,\bm z)$,
which transforms equation \eqref{jseq1} into equation \eqref{tjseq1}.
Then we prove that this transformation automatically
transforms equations \eqref{jseq2} into equations \eqref{tjseq2}
with diagonal $\mathcal P_i$'s.

If a transformation $X(s,\bm z)=G(s,\bm z)Z(s,\bm z)$,
 \beq
G(s,\bm z)=D(\bm z)^{-1}F(s,\bm z)s^{\Lambda_1(\bm z)},
\qquad
F(s,\bm z)=\mathbbm 1+\sum_{k=1}^\infty F_k(\bm z)\frac{1}{s^k},
\eneq
transforms equation \eqref{jseq1} to the equation
\beq
\frac{d}{ds}Z(s,\bm z)=U(\bm z)Z(s,\bm z),
\eneq
  then  $Z(s,\bm z)$ is a solution of the equation
\beq
\frac{d}{ds}Z(s,\bm z)=\left(G(s,\bm z)^{-1}A(s,\bm z) G(s,\bm z)-G(s,\bm z)^{-1}\frac{d}{ds}G(s,\bm z)\right)Z(s,\bm z).
\eneq
Thus,  $F(s,\bm z)$ satisfies the equation
\begin{align*}
U(\bm z)F(s,\bm z)&+\left(\sum_{k=1}^\infty \frac{D(\bm z)A_k(\bm z)D(\bm z)^{-1}}{s^k}\right)F(s,\bm z)\\
&=\frac{d}{ds}F(s,\bm z)+\frac{1}{s}F(s,\bm z)\Lambda_1(\bm z)+F(s,\bm z)U(\bm z),
\end{align*}
which  gives a system of equations for the coefficients $F_k(\bm z)$.

Denote $\hat A_j:=DA_jD^{-1}$.  The matrix $\hat A_1^{\rm od}:=\hat A_1-\Lambda_1$ is an off-diagonal matrix. 
The first equation is
\beq
UF_1+\hat A_1^{\rm od}=F_1U.
\eneq
For $\alpha\neq \beta$, 
\beq
F_1(\bm z)_{\alpha\beta}=\frac{1}{u_{\beta}(\bm z)-u_\alpha(\bm z)}\left(\hat A_1^{\rm od}(\bm z)\right)_{\alpha\beta}.
\eneq
The second equation is 
\beq
\label{doe}
UF_2+\hat A_2+\hat A_1F_1=-F_1+F_1\Lambda+F_2U.
\eneq
We find the diagonal entries of $F_1(\bm z)$ from the diagonal part of equations (\ref{doe}).  We
compute the off-diagonal entries $F_2(\bm z)_{\alpha\beta}$ with $\alpha\neq \beta$, by the formula:
\beq
F_2(\bm z)_{\alpha\beta}=\frac{1}{u_\beta(\bm z)-u_\alpha(\bm z)}\left(F_1(\bm z)-F_1(\bm z)\Lambda(\bm z)+\hat A_2(\bm z)+\hat A_1(\bm z)F_1(\bm z)\right).
\eneq
After $k$ steps of this procedure, we will determine all the coefficients $F_1,\dots, F_{k-1}$ and
all the off-diagonal entries  of $F_k$. Then the $k+1$st equation
\[
-kF_k+[F_{k+1},U]=-F_k\Lambda+\hat A_{k+1}+\sum_{h+\ell=k+1}\hat{A_h}F_\ell
\]
determines uniquely the diagonal entries of $F_k$ and the  off-diagonal entries 
of $F_{k+1}$ and we may continue  this procedure. 

This procedure shows that  the desired series
 $F(s,\bm z)$ does exist and is  unique.
 
 \smallskip
The gauge transformation $X(s,\bm z)=G(s,\bm z)Z(s,\bm z)$ transforms the joint system \eqref{jseq1}, \eqref{jseq2} into the joint system
\begin{align}
\label{nsyst1}
\frac{dZ}{ds}&=U(\bm z)Z,\\
\label{nsyst2}
Z(s,z_1,\dots, z_i-1,\dots,z_n)&=\mathcal P_i(s,\bm z)Z(s,\bm z),\quad i=1,\dots,n,
\end{align}
where we set
\beq\label{modPi}
\mathcal P_i(s,\bm z):=G(s,z_1,\dots, z_i-1,\dots, z_n)^{-1}P_i(s,\bm z)G(s,\bm z)
\eneq
for $i=1,\dots,n$. 
We claim that the matrices $\mathcal P_i$ are diagonal and  independent of $s$.

The compatibility conditions of \eqref{nsyst1} and \eqref{nsyst2} imply that
\beq
\frac{d}{ds}\mathcal P_i=[U,\mathcal P_i],\quad i=1,\dots, n.
\eneq
Thus the entries of $\mathcal P_i$ are of the form 
\beq
\label{expPi}
\mathcal P_i(s,\bm z)_{\alpha\beta}=f_{\alpha\beta}(\bm z)\exp((u_\alpha(\bm z)-u_\beta(\bm z))s),
\eneq
where $\alpha,\beta=1,\dots, n$, and $f_{\alpha\beta}(\bm z)$ are functions of $\bm z$. 

Also we know that all the entries of the right-hand side of \eqref{modPi} are formal power series of the form 
\[s^{m(\bm z)}\sum_{n=0}^\infty\frac{a_n(\bm z)}{s^n},
\]for suitable functions $m(\bm z),a_n(\bm z)$. 
This shows  that the operator $\mathcal P_i$ can be of the form  \eqref{expPi} if and only if 
\[
f_{\alpha\beta}(\bm z)=0,\quad \alpha\neq \beta.
\]
This concludes the proof.
\endproof

\section{Relation of $qDE$ to Dubrovin's equation for $QH^\bullet(\mathbb P^{n-1})$}\label{qdeDubr}
Denote by
\beq
\iota^*\colon H^\bullet_{\mathbb T}(\mathbb P^{n-1},\mathbb C)\to H^\bullet(\mathbb P^{n-1},\mathbb C),\quad f(x,\bm z)\mapsto f(x,0),
\eneq
the \emph{non-equivariant limit morphism}\footnote{Recall that this morphism is induced in cohomology by the inclusion $\iota\colon\mathbb P^{n-1}\to\mathbb P^{n-1}_{\mathbb T}:=\mathbb P^{n-1}\times_{\mathbb T}E\mathbb T$. See \cite[Section 2]{atbot}.}.
It maps the $\mathbb C[\bm z]$-basis $(x_{\alpha})_{\alpha=0}^{n-1}$ of $H^\bullet_{\mathbb T}(\mathbb P^{n-1},\mathbb C)$ to the $\mathbb C$-basis $(\iota^* x_{\alpha})_{\alpha=0}^{n-1}$ of $H^\bullet(\mathbb P^{n-1},\mathbb C)$. Denote the dual coordinates on $H^\bullet(\mathbb P^{n-1},\mathbb C)$ by the notation $\bm t:=(t^0,\dots, t^{n-1})$. Consider the non-equivariant limit $F_0^{\mathbb P^{n-1}}\in\mathbb C[\![\bm t]\!]$ of \eqref{GWpot}, 
\beq
F_0^{\mathbb P^{n-1}}(\bm t):=\sum_{m=0}^\infty\sum_{d=0}^\infty\sum_{\alpha_1,\dots\alpha_m=0}^{n-1}\frac{t^{\alpha_1}\dots t^{\alpha_m}}{m!}\langle \iota^*x_{\alpha_1},\dots, \iota^*x_{\alpha_m}\rangle_{0,m,d}^{\mathbb P^{n-1}}.
\eneq
It is known that the Gromov-Witten potential $F_0^{\mathbb P^{n-1}}(\bm t)$ is convergent. The domain of convergence $M\subseteq H^\bullet(\mathbb P^{n-1},\mathbb C)$ of the Gromov-Witten potential $F_0^{\mathbb P^{n-1}}(\bm t)$ carries a \emph{Frobenius manifold structure} \cite{dubronapoli,dubro1,dubro0,dubro2,manin,hertling,sabbah}. Tangent spaces\footnote{Tangent spaces to $M$ are canonically identified with $H^\bullet(\mathbb P^{n-1},\mathbb C)$.} $T_pM$ are equipped with an associative, commutative algebra structure: the product $*_p\colon T_pM\times T_pM\to T_pM$ is compatible with the non-equivariant Poincar\'e metric $\eta_{\rm cl}:=\eta|_{\bm z=0}$,
\beq
\eta_{\rm cl}(\alpha*_p\beta,\gamma)=\eta_{\rm cl}(\alpha,\beta*_p\gamma),\quad \alpha,\beta,\gamma\in T_pM.
\eneq  
The metric $\eta_{\rm cl}$ is a non-degenerate pseudo-riemannian metric on $M$, whose Levi-Civita connection $\nabla$ is flat.
Consider the semisimple part $M_{ss}$ of $M$, namely the subset of points $p$ whose corresponding Frobenius algebra $T_pM$ is without nilpotents. Denote by $(\pi_1,\dots, \pi_n)$ the idempotent tangent vectors at $p\in M_{ss}$, and introduce the normalized frame $(f_1,\dots, f_n)$ by
\beq\label{normidemp}
f_i:=\eta_{\rm cl}(\pi_i,\pi_i)^{-\frac{1}{2}}\cdot \pi_i,\quad i=1,\dots,n,
\eneq
for arbirary choice of the square roots. Consider the \emph{Euler vector field} $E$ on $M$ defined by
\beq
E:=c_1(\mathbb P^{n-1})+\sum_{\alpha=0}^{n-1}\left(1-\frac{1}{2}{\rm deg}(\iota^*x_{\alpha})\right)t^\alpha\frac{\partial}{\partial t^\alpha}.
\eneq
Let $p\in M_{ss}$, and denote by $U(p)$ and $V(p)$ the matrices, wrt the frame $(f_1,\dots, f_n)$, of the morphisms
\beq\label{euleru}
\mathcal U(p)\colon T_pM\to T_pM,\quad v\mapsto E|_p*_pv,
\eneq
\beq\label{gradmu}
\mu(p) \colon T_pM\to T_pM,\quad v\mapsto \frac{3-n}{2}v-\left.\nabla_v E\right|_p.
\eneq
It is easily seen that they satisfy 
\beq
U(p)^T=U(p),\quad V(p)^T+V(p)=0.
\eneq

There is a local identification of $M_{ss}$ with the space of parameters of isomonodromic deformations of the ordinary differential equation
\beq\label{isomono}
\frac{d}{d\lambda}Y(\lambda,p)=\left(U(p)+\frac{1}{\lambda} V(p)\right)Y(\lambda,p),\quad \lambda\in\mathbb C^*,\ p\in M_{ss}
\eneq
for a $n\times n$-matrix valued function $Y$. Equation \eqref{isomono} is central in Dubrovin's study of Frobenius manifolds, see \cite[Lecture 3]{dubro1}, \cite[Lectures 3 and 4]{dubro2}, \cite{dubro0}. See also \cite{guzzetti1} and \cite[Section 6]{CDG1} for details on the monodromy and Stokes phenomenon of \eqref{isomono}.

The tangent space $T_0M\cong H^\bullet(\mathbb P^{n-1},\mathbb C)$ can be identified with $\mathbb C^n$ by fixing the frame $(1,x,\dots, x^{n-1})$, where $x$ denote the non-equivariant hyperplane class. Under this identification, 
\begin{itemize}
\item the Gram matrix of the non-equivariant Poincar\'{e} metric coincides with the matrix $\eta_{\rm cl}$ of equation \eqref{nepoi};
\item the matrix of the operator $\mathcal U(0)\colon T_0M\to T_0M$ of multiplication by the Euler vector field coincides with the matrix $\mathcal B_0$ of equation \eqref{b0};
\item the basis $(f_1,\dots, f_n)$ of Lemma \ref{fbasis} coincides with the orthonormalized idempotent frame \eqref{normidemp} at $T_0M$ for suitable choices of the square roots, see \cite[Section 6.1]{CDG1}.
\end{itemize}
In the standard notations of Dubrovin's theory of Frobenius manifolds, the matrix $\mathcal E$ is usually denoted by $\Psi$. Here we avoid this notation, the symbol $\Psi$ being already used for solutions of the joint systems \eqref{eqde}, \eqref{qkz}.

We close this Appendix by commenting the relation between equation \eqref{redeqpn} and the isomonodromic differential equation \eqref{isomono}.

\begin{prop}\label{isoqde}
For $\bm z=0$, equation \eqref{redeqpn} is
\beq
\label{redeq0}\frac{d}{ds} T(s)=\left(\mathcal B_0+\frac{1}{s}\mathcal B_1(0)\right) T(s).
\eneq If $ T(s)$ is a solution of \eqref{redeq0}, then  
\[Y(\lambda):=\lambda^{-\frac{n-1}{2}}\mathcal E\  T(\lambda)
\]is a solution of the equation \eqref{isomono} specialized at $p=0$.
\end{prop}

\proof
if $\bm z=0$, all the coefficients $\mathcal B_2,\dots,\mathcal B_n$ vanish, and the coefficient $\mathcal B_1$ takes the form
\beq
\mathcal B_1(0)={\rm diag}(0,1,2,\dots, n-1).
\eneq
We have 
\[\mathcal B_1(0)-\frac{n-1}{2}\mathbbm 1=\mu,
\]where $\mu$ is the matrix of the grading operator 
\eqref{gradmu}, written in coordinates wrt the basis $(\iota^*x_\alpha)_{\alpha=0}^{n-1}$.
As a consequence, the matrix $B^{\rm od}$ of formula \eqref{propP2} in the non-equivariant limit is given by
\[B^{\rm od}(0)=\mathcal E\ \left(\mathcal B_1(0)-\frac{n-1}{2}\mathbbm 1\right)\ \mathcal E^{-1}=\mathcal E\ \mu\ \mathcal E^{-1}=V,
\]where the matrix $V$ is the antisymmetric matrix $V(p)$ specialized at $p=0$.
The antisymmetry of $B^{\rm od}$ however is lost for $\bm z\neq 0$.
\endproof

\bibliographystyle{amsalpha}

\begin{thebibliography}{GRTV13}

\bibitem[Aig13]{aigner}
M.~Aigner.
\newblock {\em Markov's theorem and 100 years of the uniqueness conjecture}.
\newblock Springer, Cham, 2013.
\newblock A mathematical journey from irrational numbers to perfect matchings.

\bibitem[AB94]{ABRH}
D.V.\,Anosov and A.A.\,Bolibruch.
\newblock {\em The Riemann-Hilbert problem.}
\newblock Aspects of Mathematics, E22, Friedr. Vieweg \& Sohn, Braunschweig, 1994.

\bibitem[AB84]{atbot}
M.~F. Atiyah and R.~Bott.
\newblock The moment map and equivariant cohomology.
\newblock {\em Topology}, 23(1):1--28, 1984.



\bibitem[BJL79a]{BJL79a}
W.\,Balser, W.B.\,Jurkat, and D.A.\,Lutz. 
\newblock A General Theory of Invariants for Meromorphic Differential Equations; Part I, Formal Invariants. 
\newblock {\em Funkcialaj Evacioj}, 22, (1979) 197-221.

\bibitem[BJL79b]{BJL79b}
W.\,Balser, W.B.\,Jurkat, and D.A.\,Lutz.
\newblock Birkhoff Invariants and Stokes Multipliers for Meromorphic Linear Differential Equations.
\newblock {\em Journal Math. Analysis and Applications}, 71, 48-94, (1979).

\bibitem[BL94]{blunts}
J.~Bernstein and V.~Lunts.
\newblock {\em Equivariant sheaves and functors}, volume 1578 of {\em Lecture
  Notes in Mathematics}.
\newblock Springer-Verlag, Berlin, 1994.

\bibitem[Bon04]{bondalsympl}
A.~I. Bondal.
\newblock A symplectic groupoid of triangular bilinear forms and the braid
  group.
\newblock {\em Izv. Rremark. Akad. Nauk Ser. Mat.}, 68(4):19--74, 2004.

\bibitem[BO18]{bororl}
L.~Borisov and D.~Orlov.
\newblock {Equivariant Exceptional Collections on Smooth Toric Stacks}.
\newblock arXiv:1806.08281v2 [math.AG], 2018.

\bibitem[CG10]{chrissginz}
N.~Chriss and V.~Ginzburg.
\newblock {\em Representation theory and complex geometry}.
\newblock Modern Birkh{\"a}user Classics. Birkh{\"a}user Boston, Inc., Boston,
  MA, 2010.
\newblock Reprint of the 1997 edition.



\bibitem[CDG18]{CDG1}
G.~Cotti, B.~Dubrovin, and D.~Guzzetti.
\newblock Helix structures in quantum cohomology of fano varieties.
\newblock arXiv:1811.09235 [math.AG], 2018.

\bibitem[CDG19]{CDG0}
G.~Cotti, B.~Dubrovin, and D.~Guzzetti.
\newblock {Isomonodromy Deformations at an Irregular Singularity with
  Coalescing Eigenvalues}.
\newblock {\em Duke Math. J.}, 168(6):967--1108, 2019.
\newblock arXiv:1706.04808v1 [math.CA].

\bibitem[CDG20]{CDG}
G.~Cotti, B.~Dubrovin, and D.~Guzzetti.
\newblock Local moduli of semisimple {Frobenius} coalescent structures.
\newblock {\em SIGMA} 16 (2020), 040, 105 pages.
\newblock arXiv:1712.08575 [math.DG].

\bibitem[CV20]{cv-inprep}
G.~Cotti and A.~Varchenko.
\newblock {The $*$-Markov equation for Laurent polynomials}.
\newblock arXiv:2006.11753 [math.AG].

\bibitem[CK99]{cox}
D.A. Cox and S.~Katz.
\newblock {\em Mirror Symmetry and Algebraic Geometry}.
\newblock American Mathematical Society, 1999.

\bibitem[Dub92]{dubronapoli}
B.A. Dubrovin.
\newblock {Integrable systems in topological field theory}.
\newblock {\em Nucl. Phys. B}, 379:627--689, 1992.

\bibitem[Dub96]{dubro1}
B.A. Dubrovin.
\newblock {Geometry of Two-dimensional topological field theories}.
\newblock In M.~Francaviglia and S.~Greco, editors, {\em Integrable Systems and
  Quantum Groups}, volume Springer Lecture Notes in Math., pages 120--348,
  1996.

\bibitem[Dub98]{dubro0}
B.A. Dubrovin.
\newblock Geometry and analytic theory of {F}robenius manifolds.
\newblock In {\em Proceedings of the {I}nternational {C}ongress of
  {M}athematicians, {V}ol. {II} ({B}erlin, 1998)}, number Extra Vol. II, pages
  315--326, 1998.

\bibitem[Dub99]{dubro2}
B.A. Dubrovin.
\newblock {Painlev{\'e} Transcendents in two-dimensional topological field
  theories}.
\newblock In R.~Conte, editor, {\em The Painlev{\'e} property, One Century
  later}. Springer, 1999.

\bibitem[Ela09]{elagin}
A.~D. Elagin.
\newblock Semi-orthogonal decompositions for derived categories of equivariant
  coherent sheaves.
\newblock {\em Izv. Rremark. Akad. Nauk Ser. Mat.}, 73(5):37--66, 2009.

\bibitem[Fed87]{fed}
M.~V. Fedoryuk.
\newblock {\em {\cyr Asimptotika}: {\cyr integraly i ryady}}.
\newblock {\cyr Spravochnaya Matematicheskaya Biblioteka}. [Mathematical
  Reference Library]. ``Nauka'', Moscow, 1987.

\bibitem[FIKN06]{painkapaev}
A.S. Fokas, A.R. Its, A.A. Kapaev, and V.~Yu. Novokshenov.
\newblock {\em {Painlev{\'e} Transcendents - The Riemann-Hilbert Approach}}.
\newblock American Mathematical Society, 2006.

\bibitem[GGI16]{gamma1}
S.~Galkin, V.~Golyshev, and H.~Iritani.
\newblock {Gamma classes and quantum cohomology of Fano manifolds: Gamma
  conjectures}.
\newblock {\em Duke Math. J.}, 165(11):2005--2077, 2016.

\bibitem[GM03]{gelman}
S.I. Gelfand and Yu.I. Manin.
\newblock {\em Methods of homological algebra}.
\newblock Springer Monographs in Mathematics. Springer-Verlag, Berlin, second
  edition, 2003.

\bibitem[Giv96]{giv1}
A.~Givental.
\newblock {Equivariant Gromov-Witten invariants}.
\newblock {\em Int. Math. Res. Not.}, 13:613--663, 1996.

\bibitem[Giv98]{Giv}
A.~Givental.
\newblock {Elliptic Gromov-Witten invariants and the generalized mirror
  conjecture.}
\newblock In {\em {Integrable systems and Algebraic Geometry. Proceedings of
  the Taniguchi Symposium 1997 (M. H. Saito, Y. Shimizu, K. Ueno, eds.)}}.
  {World Scientific}, 1998.

\bibitem[GK95]{givkim}
A.~Givental and B.~Kim.
\newblock Quantum cohomology of flag manifolds and {T}oda lattices.
\newblock {\em Comm. Math. Phys.}, 168(3):609--641, 1995.

\bibitem[GRTV13]{GRTV}
V.~Gorbounov, R.~Rim\'{a}nyi, V.~Tarasov, and A.~Varchenko.
\newblock Quantum cohomology of the cotangent bundle of a flag variety as a
  {Y}angian {B}ethe algebra.
\newblock {\em J. Geom. Phys.}, 74:56--86, 2013.

\bibitem[GK04]{helix}
A.~L. Gorodentsev and S.~A. Kuleshov.
\newblock Helix theory.
\newblock {\em Mosc. Math. J.}, 4(2):377--440, 535, 2004.

\bibitem[Guz99]{guzzetti1}
D.~Guzzetti.
\newblock Stokes matrices and monodromy of the quantum cohomology of projective
  spaces.
\newblock {\em Comm. Math. Phys.}, 207(2):341--383, 1999.

\bibitem[Her02]{hertling}
C.~Hertling.
\newblock {\em {Frobenius Manifolds and Moduli Spaces for Singularities}}.
\newblock Cambridge Univ. Press, 2002.

\bibitem[Huy06]{huy}
D.~Huybrechts.
\newblock {\em Fourier-Mukai transforms in algebraic geometry}.
\newblock Oxford Univ. Press, 2006.

\bibitem[KT08]{braids}
C.~Kassel and V.~Turaev.
\newblock {\em Braid groups}, volume 247 of {\em Graduate Texts in
  Mathematics}.
\newblock Springer, New York, 2008.
\newblock With the graphical assistance of Olivier Dodane.

\bibitem[KKP08]{KKP}
L.~Katzarkov, M.~Kontsevich, and T.~Pantev.
\newblock Hodge theoretic aspects of mirror symmetry.
\newblock In {\em From {H}odge theory to integrability and {TQFT}
  tt*-geometry}, volume~78 of {\em Proc. Sympos. Pure Math.}, pages 87--174.
  Amer. Math. Soc., Providence, RI, 2008.

\bibitem[Kim96]{kim}
B.~Kim.
\newblock On equivariant quantum cohomology.
\newblock {\em Internat. Math. Res. Notices}, (17):841--851, 1996.

\bibitem[LLY97]{LLY}
B.H. Lian, K.~Liu, and S.T. Yau.
\newblock {Mirror principle I}.
\newblock arXiv:alg-geom/9712011v1, 1997.

\bibitem[LH09]{lipman}
J.~Lipman and M.~Hashimoto.
\newblock {\em Foundations of {G}rothendieck duality for diagrams of schemes},
  volume 1960 of {\em Lecture Notes in Mathematics}.
\newblock Springer-Verlag, Berlin, 2009.

\bibitem[LS17]{liush}
C.-C.~M. Liu and A.~Sheshmani.
\newblock Equivariant {G}romov-{W}itten invariants of algebraic {GKM}
  manifolds.
\newblock {\em SIGMA Symmetry Integrability Geom. Methods Appl.}, 13:Paper No.
  048, 21, 2017.

\bibitem[LR16]{div2}
M.~Loday-Richaud.
\newblock {\em Divergent series, summability and resurgence. {II}}, volume 2154
  of {\em Lecture Notes in Mathematics}.
\newblock Springer, [Cham], 2016.
\newblock Simple and multiple summability, With prefaces by Jean-Pierre Ramis,
  \'{E}ric Delabaere, Claude Mitschi and David Sauzin.

\bibitem[Man99]{manin}
Yu.I. Manin.
\newblock {\em {Frobenius manifolds, Quantum Cohomology, and Moduli Spaces}}.
\newblock Amer. Math. Soc., Providence, RI, 1999.

\bibitem[MO19]{mo}
D.~Maulik and A.~Okounkov.
\newblock Quantum groups and quantum cohomology.
\newblock {\em Ast\'{e}risque}, (408):ix+209, 2019.

\bibitem[Mih06]{mihalcea}
L.~Mihalcea.
\newblock Equivariant quantum {S}chubert calculus.
\newblock {\em Adv. Math.}, 203(1):1--33, 2006.

\bibitem[MS16]{div1}
C.~Mitschi and D.~Sauzin.
\newblock {\em Divergent series, summability and resurgence. {I}}, volume 2153
  of {\em Lecture Notes in Mathematics}.
\newblock Springer, [Cham], 2016.
\newblock Monodromy and resurgence, With a foreword by Jean-Pierre Ramis and a
  preface by \'{E}ric Delabaere, Mich{\`e}le Loday-Richaud, Claude Mitschi and
  David Sauzin.

\bibitem[Pol11]{poli}
A.~Polishchuk.
\newblock {$K$}-theoretic exceptional collections at roots of unity.
\newblock {\em J. K-Theory}, 7(1):169--201, 2011.

\bibitem[RTV15]{RTV}
R.~Rim\'{a}nyi, V.~Tarasov, and A.~Varchenko.
\newblock Partial flag varieties, stable envelopes, and weight functions.
\newblock {\em Quantum Topol.}, 6(2):333--364, 2015.

\bibitem[Sab08]{sabbah}
C.~Sabbah.
\newblock {\em Isomonodromic deformations and Frobenius manifolds: An
  introduction}.
\newblock Springer, 2008.

\bibitem[Sib90]{sibook}
Y.~Sibuya.
\newblock {\em Linear differential equations in the complex domain: problems of
  analytic continuation}, volume~82 of {\em Translations of Mathematical
  Monographs}.
\newblock American Mathematical Society, Providence, RI, 1990.
\newblock Translated from the Japanese by the author.

\bibitem[TV97a]{tarvar971}
V.~Tarasov and A.~Varchenko.
\newblock Geometry of {$q$}-hypergeometric functions as a bridge between
  {Y}angians and quantum affine algebras.
\newblock {\em Invent. Math.}, 128(3):501--588, 1997.

\bibitem[TV97b]{tarvar972}
V.~Tarasov and A.~Varchenko.
\newblock Geometry of {$q$}-hypergeometric functions, quantum affine algebras
  and elliptic quantum groups.
\newblock {\em Ast\'{e}risque}, (246):vi+135, 1997.

\bibitem[TV14]{tarvarhyp}
V.~Tarasov and A.~Varchenko.
\newblock Hypergeometric solutions of the quantum differential equation of the
  cotangent bundle of a partial flag variety.
\newblock {\em Cent. Eur. J. Math.}, 12(5):694--710, 2014.

\bibitem[TV19a]{tarvar}
V.~Tarasov and A.~Varchenko.
\newblock {Equivariant quantum differential equation, Stokes bases, and
  K-Theory for a Projective Space}.
\newblock arXiv:1901.02990v1, 2019.

\bibitem[TV19b]{tar-var}
V.~Tarasov and A.~Varchenko.
\newblock {$q$}-hypergeometric solutions of quantum differential equations,
  quantum {P}ieri rules, and gamma theorem.
\newblock {\em J. Geom. Phys.}, 142:179--212, 2019.

\bibitem[VV10]{varaves}
M.~Varagnolo and E.~Vasserot.
\newblock Double affine {H}ecke algebras and affine flag manifolds, {I}.
\newblock In {\em Affine flag manifolds and principal bundles}, Trends Math.,
  pages 233--289. Birkh\"{a}user/Springer Basel AG, Basel, 2010.

\bibitem[Was65]{wasow}
W.~Wasow.
\newblock {\em Asymptotic expansions for ordinary differential equations}.
\newblock Pure and Applied Mathematics, Vol. XIV. Interscience Publishers John
  Wiley \& Sons, Inc., New York-London-Sydney, 1965.

\bibitem[Wit90]{witten}
E.~Witten.
\newblock On the structure of the topological phase of two-dimensional gravity.
\newblock {\em Nucl. Phys. B}, 340:281--332, 1990.

\end{thebibliography}

\end{document}